\documentclass[11pt]{article}
\usepackage[utf8]{inputenc}
\usepackage[english]{babel}
\usepackage[titletoc]{appendix}
\usepackage[centertags]{amsmath}
\usepackage{amsthm, amssymb}
\usepackage{dsfont}
\usepackage{bbm}
\usepackage{graphicx}
\usepackage{pgf,tikz}
\usepackage{caption}
\usepackage[normalem]{ulem}
\usepackage{cancel}

\usepackage[pdftex,
            pdfauthor={Max Grieshammer},
            pdftitle={Measure representation of evolving genealogies},
            pdfproducer={pdfLaTex},
            pdfcreator={}]{hyperref}

\usepackage{fancyhdr}
\numberwithin{equation}{section}
\theoremstyle{plain}
\newtheorem{theorem}{Theorem}[section]
\newtheorem{proposition}[theorem]{Proposition}

\newtheorem{lemma}[theorem]{Lemma}
\newtheorem{assumption}[theorem]{Assumption}

\theoremstyle{definition}

\newtheorem{definition}[theorem]{Definition}

\theoremstyle{remark}
\newtheorem{remark}[theorem]{Remark}

\def\B{{\mathcal B}}
\def\M{\mathcal M}
\def\e{\varepsilon}
\def\L{\mathcal L}
\def\supp{\textrm{\normalfont supp}}

\newcommand{\norm}[1]{\big|\hspace{-0.05cm}\big|#1\big|\hspace{-0.05cm}\big|}

\def\mfx{\mathfrak{x}}

\def\SMR{\mathfrak{M}}
\def\TYPE{\mathbb{T}}							 							
\def\type{\mathfrak{t}}						 							
\def\MVX{\mathcal{X}}           	 							
\def\UM{\mathbb{U}}							 	 							
\def\MM{\mathbb{M}}							 	 							
\def\UMK{\mathbb{U}^{\TYPE}}	 									
\def\UMG{\mathbb{U}^{G}}	 											
\def\MMK{\mathbb{M}^{\TYPE}}	 									

\def\UMI{\mathbb{I}}
\def\UMIK{\mathbb{I}^{\TYPE}}
\def\mfu{\mathfrak{u}}						 							
\def\U{\mathcal{U}}								 							
\def\Cut{\Phi}							  	   							
\def\CutTwo{\hat \Phi}
\def\KERNEL{\kappa}															
\def\dmG{d_{\text{mGP}}}											  
\def\dPr{d_{\text{Pr}}}													
\def\dGP{d_{\text{GPr}}}											  
\def\rep{\mathfrak{r}}													
\def\h{h}
\def\num{n}												  
\def\F{\mathfrak{F}}												 		
\def\f{\mathfrak{f}}												  	
\def\ord{\mathfrak{o}}												 	
\def\DES{\mathcal{D}}

\def\king{\kappa} 
\def\King{\mathfrak{K}} 
\newcommand{\freq}[1]{|#1|_f}
\def\SEG{\mathbb E}											

\def\dTYPE{d^{\TYPE}} 
\def\dSK{d^{\text{SK}}} 
\def\dSeq{d^1} 
\def\pois{\rho} 
\def\MS{\mathfrak{S}}

\def\proj{\pi}
\def\eps{\varepsilon}
\def\isom{\tau}

\begin{document}

 \title{Measure representation of evolving genealogies}
 \author{Max Grieshammer\footnote{Institute for Mathematics, Friedrich-Alexander Universit\"at Erlangen-N\"urnberg, Germany;  	max.grieshammer@math.uni-erlangen.de, 
 MG was supported by DFG-Grant GR-876-17.1 of A. Greven}
 }
 
\maketitle
 
\thispagestyle{empty}

\begin{abstract}
We study evolving genealogies, i.e. processes that take values in the space of (marked) ultra-metric measure spaces and satisfy some sort of ``consistency'' condition. 
This condition is based on the observation that the genealogical distance of two individuals who do not have common ancestors up to a time $h$ in the past is completely determined by the genealogical distance of the respective ancestors at that time $h$ in the past. \par
Now the idea is to color all possible ancestors at time $h$ in the past and measure the relative number of their descendants. The resulting collection of measure-valued processes (the construction is possible for all $h$) is called a measure representation. \par 
As a main result we give a tightness criterion of evolving genealogies in terms of their measure representation.  
We then apply our theory to study a finite system scheme for tree-valued interacting Fleming-Viot processes. 
\end{abstract}
\medskip 

\noindent {\bf Keywords:} Evolving genealogies, marked (ultra-)metric measure spaces, tightness, measure-valued processes, finite-system scheme, Fleming-Viot model. 

\smallskip
\noindent {\bf AMS 2010 Subject Classification:}  Primary: 60G07, 60J25, 60B10; Secondary: 92D10, 92D25.\\

\newpage

\tableofcontents
\thispagestyle{empty}
\section{Introduction}

One objective in studying population models is the description of the genealogical tree of the population alive at the present time $T$. To do that there are several approaches. 
The classical approach is to study a process, which generates or describe the genealogy backward in time. For example, when we consider a Wright-Fisher population, then one can use the Kingman coalescent, see \cite{kingman1982coalescent}, to obtain genealogical information. Besides the coalescing models, there are other classical approaches like the lookdown construction of Donelly and Kurtz (see \cite{donnelly1996countable} and \cite{donnelly1999genealogical}). \par 
But in this paper we want to follow a different idea. Namely, the idea, that the genealogy of a population has an ``history'', i.e. it evolved as the population evolved. In other words, one uses a Markov dynamic on a suitable state space that generates the genealogy forward in time. This idea is based on the works of Depperschmidt, Greven, Pfaffelhuber and Winter (see \cite{GPW09}, \cite{GPW13}, \cite{DGP11} and \cite{DGP12} for the main references). \\

To be a bit more precise, we model the genealogical relationships in a population $X$, by the so called {\it genealogical distance} $r$, i.e.
$r(x,y)$ is the distance to the most resent common ancestor of the two individuals $x,y \in X$. Additionally, since we want to sample individuals from the population, we need a sampling measure $\tilde \mu$. 
Finally, since an individual $x$ carries some type $\type \in \TYPE$, with probability $\KERNEL(x,d\type)$, 
where we assume $\TYPE$ is a compact separable metric space, we formally get a measure $\mu(dx,d\type) = \KERNEL(x,d\type)\tilde \mu(dx) \in \mathcal M_1(X\times \TYPE)$. 
That means, we model genealogies of populations by the triples $(X,r,\mu)$, i.e. by a {\it marked metric measure spaces} or more precise, since the genealogical distance $r$ is an ultra-metric, 
by {\it marked ultra-metric measure spaces}. \par 
Since we are only interested in the genealogical structure of a population, a relabeling of the individuals (that preserves types and masses) should result in the same genealogical structure. 
Formally, two marked ultra-metric measure spaces $(X,r,\mu)$ and $(X',r',\mu')$ are said to be {\it equivalent}, if we can find an isometry 
$\tilde \varphi: \supp(\mu(\cdot \times \TYPE)) \to \supp(\mu'(\cdot \times \TYPE)))$ and $\varphi(x,\type):=(\tilde \varphi(x),\type)$ is measure-preserving. 
We denote by $[\cdot]$ the equivalence classes and define the space 
\begin{equation}
\UMK_1 := \{[X,r,\mu]: (X,r,\mu) \text{ is a marked ultra-metric measure space}\}.
\end{equation}

More general, we are working with finite measures $\mu(dx,d\type) = \KERNEL(x,d\type)\tilde \mu(dx)$, where $\KERNEL$ is a probability kernel and $\tilde \mu$ is a {\it finite} Borel-measure on $X$. 
The resulting space is denoted by $\UMK$ (see Section \ref{sec_basic_objects} for the details). Note that when $|\TYPE| = 1$, we can identify the space $\UMK$ with the space $\UM$ of (equivalent classes) of ultra-metric measure spaces (without marks). \par 
We note that ultra-metric space can be mapped isometrically to the set of leaves of a rooted $\mathbb R$-tree (see Remark 2.7 in \cite{DGP12} and Remark 2.2 in \cite{DGP11}), 
that is why we sometimes call an element in $\UMK$ a tree (and formally justifies the figures below, where we illustrate the genealogical distance by the genealogical tree). \\

\begin{figure}[ht]
\begin{tikzpicture}[scale = 0.37]
\draw  (0,0)-- (0,12);
\draw  (2,0)-- (2,12);
\draw  (4,0)-- (4,12);
\draw  (6,0)-- (6,12);
\draw [->] (8,-2) -- (8,14);
\draw [->] (6,10) -- (2,10);
\draw [->] (0,7) -- (4,7);
\draw [->] (6,4) -- (0,4);
\draw  (10,5)-- (10,4);
\draw  (12,5)-- (12,4);
\draw  (10,4)-- (12,4);
\draw  (11,4)-- (11,0);
\draw  (11,0)-- (14,0);
\draw  (14,0)-- (14,5);
\draw  (14,0)-- (16,0);
\draw  (16,0)-- (16,5);
\draw  (18,8)-- (18,7);
\draw  (18,7)-- (20,7);
\draw  (20,8)-- (20,7);
\draw  (19,7)-- (19,4);
\draw  (19,4)-- (22,4);
\draw  (22,4)-- (22,8);
\draw  (24,8)-- (24,0);
\draw  (24,0)-- (20.5,0);
\draw  (20.5,0)-- (20.5,4);
\draw  (26,11)-- (26,10);
\draw  (26,10)-- (28,10);
\draw  (28,10)-- (28,11);
\draw  (27,10)-- (27,7);
\draw  (27,7)-- (30,7);
\draw  (30,7)-- (30,11);
\draw  (28.5,7)-- (28.5,4);
\draw  (28.5,4)-- (32,4);
\draw  (32,4)-- (32,11);
\draw [dash pattern=on 7pt off 7pt] (0,5)-- (34,5);
\draw [dash pattern=on 7pt off 7pt] (0,8)-- (34,8);
\draw [dash pattern=on 7pt off 7pt] (0,11)-- (34,11);
\draw (8.25,12.25) node[anchor=north west] {$T_3$};
\draw (8.25,9.25) node[anchor=north west] {$T_2$};
\draw (8.25,6.25) node[anchor=north west] {$T_1$};
\draw (8,15) node[anchor=north west] {$time$};
\draw (8,11)-- ++(-5pt,-5pt) -- ++(10.0pt,10.0pt) ++(-10.0pt,0) -- ++(10.0pt,-10.0pt);
\draw (8,8)-- ++(-5pt,-5pt) -- ++(10.0pt,10.0pt) ++(-10.0pt,0) -- ++(10.0pt,-10.0pt);
\draw (8,5)-- ++(-5pt,-5pt) -- ++(10.0pt,10.0pt) ++(-10.0pt,0) -- ++(10.0pt,-10.0pt);
\end{tikzpicture}
\caption{\footnotesize Basic idea of evolving genealogies: On the left side we draw a population consisting of four individuals. 
The arrows indicate that at certain times one individual dies (the individual at the top of the arrow) and is replaced by an offspring of the individual at the bottom of the arrow. 
At the right side we draw the genealogical trees at certain times.}\label{fig.1}
\end{figure}
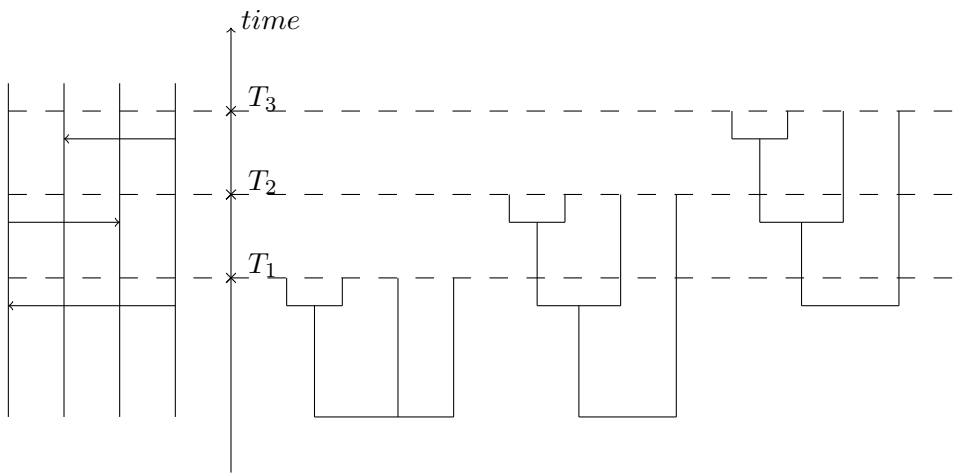

The aim of this paper is to give a strong tightness result for evolving genealogies, in terms of measure-valued processes. 
In order to get an idea, we first need to define the notion of evolving genealogy. To do so, we restrict at this point to the non-marked setting. 
Moreover, we do not want to go into technical details yet, i.e. we want to give a heuristic picture first. For the details and 
the general case see Section \ref{sec_main}.

\begin{figure}[ht]
\begin{tikzpicture}[scale = 0.45]
\draw [->] (8,-2) -- (8,14);
\draw  (10,5)-- (10,4);
\draw  (12,5)-- (12,4);
\draw  (10,4)-- (12,4);
\draw  (11,4)-- (11,0);
\draw  (11,0)-- (14,0);
\draw  (14,0)-- (14,5);
\draw  (14,0)-- (16,0);
\draw  (16,0)-- (16,5);
\draw [dash pattern=on 4pt off 4pt] (18,8)-- (18,7);
\draw [dash pattern=on 4pt off 4pt] (18,7)-- (20,7);
\draw [dash pattern=on 4pt off 4pt] (20,8)-- (20,7);
\draw  (19,4)-- (22,4);
\draw  (24,0)-- (20.5,0);
\draw  (20.5,0)-- (20.5,4);
\draw [dash pattern=on 4pt off 4pt] (26,11)-- (26,10);
\draw [dash pattern=on 4pt off 4pt] (26,10)-- (28,10);
\draw [dash pattern=on 4pt off 4pt] (28,10)-- (28,11);
\draw [dash pattern=on 4pt off 4pt] (27,7)-- (30,7);
\draw  (28.5,4)-- (32,4);
\draw (8.5,11.5) node {$T_3$};
\draw (8.5,8.5) node {$T_2$};
\draw (8.5,5.5) node {$T_1$};
\draw (9,14) node {$time$};
\draw [dash pattern=on 4pt off 4pt] (8,5)-- (34,5);
\draw  (19,4)-- (19,5);
\draw [dash pattern=on 4pt off 4pt] (19,5)-- (19,7);
\draw  (24,0)-- (24,5);
\draw [dash pattern=on 4pt off 4pt] (24,5)-- (24,8);
\draw  (22,4)-- (22,5);
\draw [dash pattern=on 4pt off 4pt] (22,5)-- (22,8);
\draw  (28.5,4)-- (28.5,5);
\draw  (32,4)-- (32,5);
\draw [dash pattern=on 4pt off 4pt] (28.5,5)-- (28.5,7);
\draw [dash pattern=on 4pt off 4pt] (32,5)-- (32,11);
\draw [dash pattern=on 4pt off 4pt] (30,7)-- (30,11);
\draw [dash pattern=on 4pt off 4pt] (27,7)-- (27,10);
\draw (10,5.75) node {$1$};
\draw (12,5.75) node {$1$};
\draw (14,5.75) node {$1$};
\draw (16,5.75) node {$1$};
\draw (19.5,5.75) node {$2$};
\draw (22.5,5.75) node {$1$};
\draw (24.5,5.75) node {$1$};
\draw (29,5.75) node {$3$};
\draw (32.5,5.75) node {$1$};
\draw (8,11)-- ++(-5pt,-5pt) -- ++(10.0pt,10.0pt) ++(-10.0pt,0) -- ++(10.0pt,-10.0pt);
\draw (8,8)-- ++(-5pt,-5pt) -- ++(10.0pt,10.0pt) ++(-10.0pt,0) -- ++(10.0pt,-10.0pt);
\draw (8,5)-- ++(-5pt,-5pt) -- ++(10.0pt,10.0pt) ++(-10.0pt,0) -- ++(10.0pt,-10.0pt);
\end{tikzpicture}
\caption{\footnotesize We fix a time $T_1$. The second truncated tree, for example, is obtained by taking the second tree from Figure \ref{fig.1}, cut it at height $T_1$ and give the new leaves the masses of the (old) leaves that ``descent'' from that new leaves.}\label{fig.2}
\end{figure}
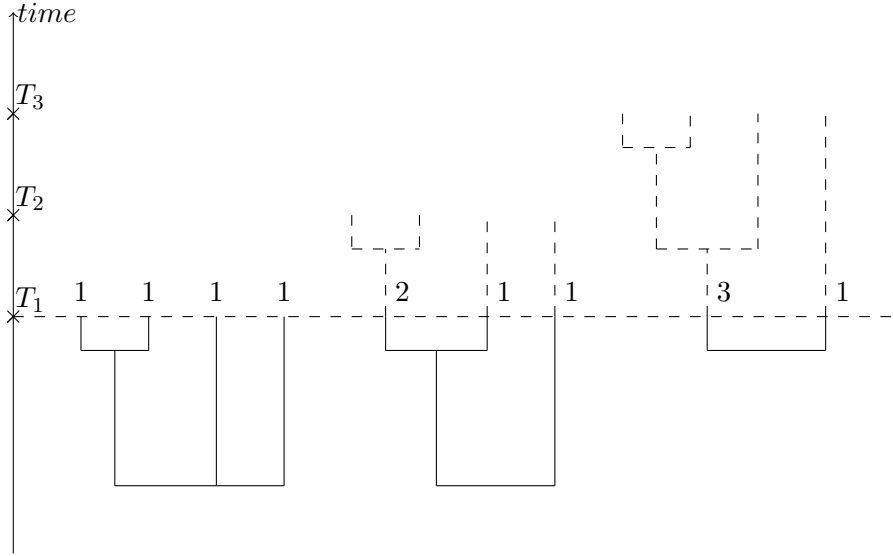

The prototype of an example we have in mind is shown in Figure \ref{fig.1}.  It is important to note that 
the trees on the right evolve in continuous time and we denote the path by $(\mfu_T)_{T \ge 0} \in D([0,\infty),\UM)$. 
Moreover, we observe that they are related in such a way that when we follow the lines backward from $T_3$ to $T_2$ in the left picture, then the third tree after time $T_2$ (measured backwards) 
is completely determined by the genealogical structure at time $T_2$. In other words, the genealogical structure of the individuals at the present time $T_3$ 
that do not have a common ancestor in the past $T_3-T_2$ years is completely determined by the ancestral 
tree at time $T_2$. 
In order to formalize this observation, we consider Figure \ref{fig.2}.

First observe that for $\mfu = [X,r,\mu] \in \UM$ and all $h > 0$ (this is the depth - measured from above - where we cut the tree), we can partition the space $X$ in at most countable 
(since $(X,r)$ is separable) many disjoint balls (since $r$ is an ultra-metric). When we interpret $r$ as a genealogical distance, the closed ball $\bar B(x,h)$ of radius $h$ around $x$ 
can be interpreted as a certain family of age $h$. In other words we can decompose the population in at most countable many disjoint families for all ages $h > 0$. 
We denote by $n(h)\in\mathbb N \cup\{\infty\}$ the number of families and pick a member $\rep_1^h,\rep_2^h,\ldots,\rep_{n(h)}^h$ of each family. By definition, the (genealogical) 
distance of $\rep_i^h$ and $\rep_j^h$ is larger than $h$. Now the truncated tree, which we call {\it $h$-trunk} is defined by
\begin{equation}
\begin{split}
\Cut^h(\mfu)&= \left[\{\rep_1^h,\rep_2^h,\ldots,\rep_{n(h)}^h\},r-h_{\neq},\sum_{i=1}^{n(h)}\mu(\bar B(\rep^h_i,h))\delta_{\rep^h_i}\right],
\end{split}
\end{equation}   
where $(r-h_{\neq})(x,y) = r(x,y)-h$ if $x \neq y$ and $0$ otherwise, and $\mu(\bar B(\rep^h_i,h))$ is the relative number (with respect to the measure $\mu$) of family members. \par 
Now, the notion evolving genealogy is defined in terms of the ``$T$ slices'' with respect to $\Cut$ given in Figure \ref{fig.2}:
We call $(\mfu_T)_{T \ge 0} = ([X_T,r_T,\mu_T])_{T \ge 0}$ a {\it (path of an) evolving genealogy}, if 
\begin{itemize}
\item[] for all $T \ge 0$, all $0 \le h' \le h$ and $\Cut^h(\mfu_{T+h}) = [[0,1],r_{T+h},\mu_{T+h}^h]$, there is an isometry $ \varphi: (\supp(\mu_{T+h}^h),r_{T+h}) \to (\supp(\mu_{T+h'}^{h'}),r_{T+h'}) $. 
\end{itemize}

Finally, we need the notion of measure representation of evolving genealogies, where we say that a collection $\{\mfx^T:\ T \ge 0\}$ is 
a {\it measure representation} of an evolving genealogy $(\mfu_T)_{T \ge 0}$ if (compare also Figure \ref{fig.3})
\begin{itemize} 
\item[(i)] $\mfx^T \in D([0,\infty),\mathcal M_f[0,1])$ for $T \ge 0$ with $\supp(\mfx_h^T) \subset \supp(\mfx_{h'}^T)$ for $0 \le h' \le h$. 
\item[(ii)] there is a metric $r_T$ such that $[[0,1],r_T,\mfx^T_h] = \Cut^h(\mfu_{T+h})$, for all $T \ge 0$ and all $h \ge 0$. 
\end{itemize}

\begin{figure}
\begin{tikzpicture}[scale = 0.52]
\draw (10,5)-- (10,4);
\draw (12,5)-- (12,4);
\draw (10,4)-- (12,4);
\draw (11,4)-- (11,0);
\draw (11,0)-- (14,0);
\draw (14,0)-- (14,5);
\draw (14,0)-- (16,0);
\draw (16,0)-- (16,5);
\draw (19,4)-- (22,4);
\draw (24,0)-- (20.5,0);
\draw (20.5,0)-- (20.5,4);
\draw (28.5,4)-- (32,4);
\draw (19,4)-- (19,5);
\draw (22,4)-- (22,5);
\draw (28.5,4)-- (28.5,5);
\draw (24,0)-- (24,5);
\draw (32,4)-- (32,5);
\draw (13,-1) node {$[X,r,\mathcal X_0^{T_1}]$};
\draw (13,6.5) node {$\mathcal X_0^{T_1} = \frac{1}{4}(\delta_1+\ldots + \delta_4)$};
\draw (22,-1) node {$[X,r,\mathcal X_{T_2-T_1}^{T_1}]$};
\draw (22,6.5) node {$\mathcal X_{T_2-T_1}^{T_1} = \frac{1}{2} \delta_1+\frac{1}{4} \delta_2 + \frac{1}{4} \delta_4$};
\draw (30.2,-1) node {$ [X,r,\mathcal X_{T_3-T_1}^{T_1}]$};
\draw (30.2,6.5) node {$\mathcal X_{T_3-T_1}^{T_1} = \frac{3}{4} \delta_1+ \frac{1}{4} \delta_4$};
\end{tikzpicture}
\caption{\footnotesize Idea of measure representations: Roughly speaking, a measure representation is a collection of paths of measures $\{\mfx^T:\ T \ge 0\}$, 
indexed by the time $T$ of the forward evolution, where $\mfx^{T}_h$ is the measure on the leaves of the truncated trees given in Figure \ref{fig.2}.}\label{fig.3}
\end{figure}

Before we continue, we define the corresponding random objects: We say that a stochastic process $\U := (\U_T)_{T \ge 0}$ with values in $\UMK$ is an evolving genealogy, if $\U(\omega)$ is an evolving genealogy for almost all $\omega$ and we call a collection $\{\mathcal X^T: \ T \ge 0\}$ of stochastic processes with values in $\mathcal M_f([0,1])$, a measure representation, if it satisfies $(\mathcal X^T)_{T \in I} \stackrel{d}{=} (\hat {\mathcal X}^T)_{T \in I} $, $I\subset [0,\infty),\ |I| < \infty$, where $(\hat X^T)_{T \in I} $ is a ``strong measure representation'' in the sense that for almost all $\omega$, $\{\hat {\mathcal X}^T(\omega): \ T \ge 0\}$ is a measure representation of $\U(\omega)$. \\

The reason why we study measure representations, is that they contain a lot of (in some situations all) information of the genealogical structure (see Section \ref{sec_main} for the details). \par 
For example, the fixation time of the measure representation is related to the time to the most recent common ancestor (see Remark \ref{rem.MCRA1}). This observation opens up a new opportunity to study questions such as the time to the most recent common ancestor using diffusion theory (see Remark \ref{rem.MCRA2}). \\

Besides examples of this kind, which follow directly from the construction, our main result is a tightness criterion: Let $\U^N$, $N \in \mathbb N$ be a sequence of evolving genealogies with measure representations $\{\mathcal X^{T,N}: T \ge 0\}$, then up to some weak assumptions on the processes it basically says that 
\begin{itemize}
\item[] if $\mathcal L(\mathcal X^{N,T})$ converges for all $T$ as a process in $D([0,\infty),\mathcal M_f([0,1]))$, where $\mathcal M_f([0,1])$ 
is equipped with the {\it weak atomic} topology (see Definition \ref{def.weak.atomic}), and the limit process $\mathcal X^{T}$ is purely atomic for all positive times $h > 0$, 
then $\mathcal L(\U^N)$ is tight as measures on $D([0,\infty),\UMK)$,
\end{itemize}
where
\begin{itemize}
\item[] the weak atomic topology was introduced by \cite{EKatomic}. Convergence in this topology is equivalent to weak convergence together with a condition on the convergence of the atoms of the measures. It is exactly the topology that is induced by convergence of 
the corresponding cumulative distribution functions in the Skorohod topology (see Proposition 6.2 in \cite{decomposition}).
\end{itemize}
At this point, we do not want to go into details concerning the topology and the additional assumptions (see Section \ref{sec_basic_objects}
and Section \ref{sec_main}), we just want to point out, that this is a really strong result in the sense that we only need to prove weak convergence 
of certain measure-valued processes in order to deduce tightness of evolving genealogies. \par 
Clearly, there are other questions, like how one can transfer the property of being an evolving genealogy to a limit process. We refer to Section \ref{sec_main} for the details and the answer to that question. \\

As an example of an evolving genealogy, we consider the tree-valued interacting Moran models.  Recall that  the tree-valued interacting Moran models describe (dynamically) the genealogical 
relationship between individuals that live on some geographical space $G$ and evolve as follows: Individuals located at the same site act like non-spatial Moran models, 
i.e. after some exponential time we uniformly pick two individuals and then one of the two individuals dies and is replaced by an offspring of the other. 
In addition individuals are allowed to migrate according 
to some homogeneous kernel $a(\cdot,\cdot)$ on $G$. \par 
If we now take the large population limit of this tree-valued process, we get the tree-valued interacting Fleming-Viot process $\U$.
As a main application of our theory on evolving genealogies, we will prove a {\it finite system scheme} result for this limit process. 

The finite system scheme clarifies what spatially infinite systems can tell us about large finite 
systems in a rigorous fashion. This has been proved for infinite particle systems, 
interacting diffusions and measure-valued processes (see \cite{MR1355055}, \cite{MR1050744} and \cite{DGV}).
What really happens on the level of genealogies is an important issue. \par 

To be precise, assume $\U^N$ are tree-valued interacting Fleming-Viot processes on some finite geographical spaces $G_N$ with suitable migration kernels $a^N$, where $G_N \uparrow G$ and $G$ is infinite. 
If $a^N$ converges to some transient migration kernel $a$ on $G$, then clearly, the mutual distance between two individuals grows to infinity when $N \rightarrow \infty$. 
The question is how fast do these distances grow and can we show some convergence result under a suitable rescaling. Here we will consider the case $G = \mathbb Z^d$, $d \ge 3$ 
and $G_N = (-N,N]^d \cap \mathbb Z^d$. In this situation, the answer is that the speed of divergence is proportional to $|G_N|$, the size of the finite spaces. 
If we rescale the distances and the ``sampling'' measures by the factor $1/|G_N|$ and speed up time by the factor $|G_N|$, we will prove, that these rescaled processes converge to a unique limit process, 
namely again a (non-spatial) tree-valued Fleming-Viot process with a certain resampling rate. \par 
The key observations for the proof are the following: One can show that  $\U^N$ is an evolving genealogy and the measure representation $\{\mathcal X^{N,T}:\ T\ge 0\}$  
can be chosen such that $\mathcal X^{N,T}$ is a system  of measure-valued interacting Fleming-Viot processes for all $T \ge 0$ and that the coupling of the processes $\mathcal X^{N,T}$ 
for different $T$ is determined by a spatial Kingman coalescent.  Since  a finite system scheme result for the measure-valued process (see \cite{DGV}) and for the spatial Kingman-coalescent (see \cite{LS}) is known,
we can apply our theory to get the corresponding result for the tree-valued Fleming-Viot process. \\

We close this section by remarking, that some of the results are based on ideas in \cite{Grieshammer}. But for the sake of completeness
and since we have a different perspective on the problem (i.e. some of the results are non-trivial reformulations), we include all proofs needed for the results in this paper. \\


\section{Marked metric measure spaces}\label{sec_basic_objects}
		
	Here we give the definition and basic properties of marked metric measure spaces and the Gromov-weak topology (see \cite{DGP11} and \cite{GPW09}). \par 
	First, recall that the support, $\supp(\mu)$, of a finite Borel measure $\mu$, 
	on some separable metric space $(X,d)$ is the smallest closed set $C$ with $\mu(X\backslash C) = 0$. Equivalently,  
	$\supp(\mu)$ is given by
	\begin{equation}\label{eq_supp}
		\supp(\mu) = \{x \in X\big|\ \forall \eps >0: \ \mu(B(x,\eps)) > 0\},
	\end{equation}
	where $B(x,\eps)$ is the open ball of radius $\eps$ around $x$.
	\begin{assumption}
	Throughout this paper we assume that
	\begin{equation}\label{eq.Type.metric}
		(\TYPE,d^\TYPE) \text{ is a separable compact metric space.}
	\end{equation}
	\end{assumption}

	\begin{definition}\label{def.mms} 
	  We call the triple $(X,r,\mu)$ 
	  \begin{enumerate}
	    \item a {\it $\TYPE$-marked metric measure space}, short mmm-space, if 
	      \begin{enumerate}
			      \item[(a)] $(X,r)$ is a complete separable metric space, where we assume that $X \subset \mathbb R$ (one needs this to
			      avoid set theoretic pathologies).
			      \item[(b)] $\mu \in \M_f(X\times \TYPE)$, i.e. $\mu$ is a finite measure on the Borel sets with respect to the product topology 
			      generated by $r$ and $d^\TYPE$. We will write $\mu(dx,d\type) = \KERNEL(x,d\type) \tilde \mu(dx)$ for some probability kernel 
			      $\KERNEL$ and $\tilde \mu:= \mu \circ \proj_X^{-1}$, where	$\proj_X:X\times\TYPE \to X$ is the projection. Note that we will
			always use $\tilde{\cdot}$ in order to indicate the projection to the first component. 
	      \end{enumerate} 
	    \item {\it ultra-metric}, if $r(x_1,x_2) \le r(x_1,x_3)\vee r(x_3,x_2)$, for $\tilde \mu$ almost all $x_1,x_2,x_3$.
	    \item {\it purely atomic}, if $\sum_{x \in X} \tilde \mu(\{x\}) = \tilde \mu(X)$.
	  \end{enumerate}		  
	  We say that two mmm-spaces $(X,r_X,\mu_X)$ and $(Y,r_Y,\mu_Y)$ are {\it equivalent} if there is a map
	  $\tilde \varphi: \supp(\tilde \mu_X) \rightarrow \supp(\tilde \mu_Y)$ with $r_X(x,y) = r_Y(\tilde \varphi(x),\tilde \varphi(y))$, 
	  $x,y \in \supp(\tilde \mu_X)$ and $\mu_Y = \mu_X\circ \varphi^{-1}$, where $\varphi(x,\type) := (\tilde \varphi(x),\type)$.\par 
	  This property defines an equivalence relation, and we denote by $[X,r,\mu]$ the equivalence 
	  class of a mmm-space $(X,r,\mu)$. Finally, we define the following sets: 
	  \begin{align}
	      \MMK &= \left\{[X,r,\mu]: \ (X,r,\mu) \textrm{ is a marked metric measure space}\right\}, \\
	      \UMK &= \left\{[X,r,\mu] \in \MMK: \ (X,r,\mu) \textrm{ is a ultra-metric measure space}\right\}.
	    \end{align}
	\end{definition}

	\begin{remark}
	  (1) Since we assume that $X\times \TYPE$ is Polish, we can always decompose $\mu(dx,d\type) = \KERNEL(x,d\type) \tilde \mu(dx)$ for some kernel $\KERNEL$.\par 
	  (2) We call the triple $(X,r,\mu)$, where $\mu \in \mathcal M_f(X)$, a {\it metric measure space}. As above we denote by  
	    $\MM$ and $\UM$  the spaces of equivalence classes induced by measure 	
	    preserving isometries. We note that if $|\TYPE| = 1$, we can identify the spaces $\MMK$ with $\MM$, $\UMK$ with $\UM$.
	\end{remark}

	In view of the above remark, we define for $\mfu = [X,r,\mu] \in \UMK$
	\begin{align}\label{eq.projMM}
	  \tilde \mfu := [X,r,\tilde \mu] \in \UM. 
	\end{align}

	\begin{definition}\label{def_dmd}Let $k \in \mathbb N_{\ge 2}$, $\mfu = [X,r,\mu]\in \UMK$ and set
	\begin{equation}\label{eq_distmap}
	 R^{k,(X,r)}: \left\{
	  \begin{array}{ll}
	      (X \times \TYPE)^{k} \rightarrow \mathbb R_+^{\binom{k}{2}} \times \TYPE^k,  \\
	      (x_{i},\type_{i})_{1\le i \le k} \mapsto  \left((r(x_{i},x_{j}))_{1\le i < j\le k},(\type_i)_{1\le i\le k}\right).
	  \end{array} \right.
	\end{equation} 
	We define the {\normalfont  distance matrix distribution of order $k$} by:
	\begin{equation}\label{eq_mdmd}
	  \begin{split}
	    \nu^{k,\mfu} :=  (R^{k,(X,r)})_*  \mu^{\otimes k} \in \mathcal{M}_f\left(\mathbb R_+^{\binom{k}{2}}\times \TYPE^k\right).
	  \end{split}
	\end{equation}
	For $k = 1$ we define (recall that $\KERNEL$ is a probability kernel)
	\begin{equation}
	  \nu^{1,\mfu} := \overline{\mfu} := \tilde \mu(X) = \sqrt{ \nu^{2,\mfu}(\mathbb R_+ \times \TYPE)}.
	\end{equation}
	\end{definition}

	\begin{remark} 
	  Note that $\nu^{k,\mfu}$ in the above definition does not depend on the representative $(X,r,\mu)$ of $\mfu$. In particular 
	  $\nu^{k,\mfu} $ is well defined.
	\end{remark}

	\begin{definition}\label{def_GS}
	  Let $\mfu,\mfu_1,\mfu_2, \ldots \in \UMK$. We say $\mfu_n \rightarrow \mfu$ for $n \rightarrow \infty$ in the {\it marked Gromov-weak
	  topology}, if
	  \begin{equation}
	    \nu^{k, \mfu_n} \stackrel{n \rightarrow \infty}{\Longrightarrow} \nu^{k, \mfu}
	  \end{equation}
	  in the weak topology on $\M_f\left( \mathbb R_+^{\binom{k}{2}}\times \TYPE^k\right)$ for all $k \in \mathbb N$, where $ \mathbb R_+^{\binom{k}{2}}\times \TYPE^k$ is equipped with the product topology.
	\end{definition}

	\begin{remark}\label{r.totalMassCont}
	  Assume that $\UMK$ is equipped with the marked Gromov-weak topology, then it is easy to see that the map $\sim: \UMK \to \UM, \ \mfu \mapsto \tilde \mfu$ is continuous. Furthermore, since $\overline{\mfu} = \sqrt{ \nu^{2,\tilde \mfu}(\mathbb R_+)}$, we also have $\mfu \mapsto \overline{\mfu}$ is continuous.
	\end{remark}

	\begin{proposition}\label{prop_M_polish}
	  $\UMK$ equipped with the Gromov-weak is Polish. 
	\end{proposition}

	\begin{proof}
	  This is Theorem 2 in \cite{DGP11}.
	\end{proof}
	
	\begin{remark}
	  An example of a complete metric is the so called {\it marked Gromov Prohorov distance} $\dmG$ (see \cite{DGP11}). We will give a definition of this distance in Section \ref{sec.GPM}.
	\end{remark}

\section{Main results: Measure representation of evolving genealogies}\label{sec_main}

		We start with the following lemma (see Lemma 3.1 in \cite{decomposition}). 

		\begin{lemma}\label{lem_representatives} 
		  Let $0< h$, $\mfu = [X,r,\mu] \in \UM$ and  $\bar B(x,h)$ be the closed ball of radius 
		  $\le h$ around $x \in X$. Then there is a $\num(h) \in \mathbb N \cup \{\infty\}$ and a family 
		  $\{\rep_i^\h: \ i \in \{1,2,\ldots,\num(h)\}\}$ of elements of $\supp(\mu)$ such that 
		  \begin{align}
		    &\mu\big(\bar B(\rep_i^h,h) \cap \bar B(\rep_j^h,h)\big) = 0, \qquad i \neq j \\
		    &\mu(X) = \sum_{i =1}^{\num(h)} \mu\big(\bar B(\rep_i^h,h)\big).
		  \end{align}
		  Moreover, if $0< \delta \le h$, then  there is a partition $\{I_i\}_{i \in 1,\ldots,n(h)}$ of $\{1,\ldots,\num(\delta)\}$ such
		  that
		  \begin{equation}
		    \mu(\bar B(\rep_i^h,h))= \sum_{j \in I_i}\mu(\bar B(\rep_j^\delta,\delta)),\qquad \forall i = 1,\ldots,\num(h).
		  \end{equation}
		\end{lemma}

		\begin{remark}\label{rem_pos} 
		  (1) By the definition of the support we get $ \mu(\bar B(\rep_i^h,h)) >0$ for all $i \in \{1,\ldots,\num(h)\}$. \par 
		  (2) The analogue of Lemma \ref{lem_representatives}  holds if we replace $\le h$ by $< h$. 
		\end{remark}
		
		\begin{definition}\label{def_Phi}
		  Let $\mfu = [X,r,\mu] \in \UMK$, $\mu(A \times B) = \int_A \KERNEL(x,B) \tilde \mu(dx)$ for some probability kernel $\KERNEL$
		  and set 
		  \begin{equation}\label{eq.muh}
		    \mu^h(A \times B):= \sum_{i \in \{1,\ldots,\num(h)\}}\int_{\bar B(\rep_i^h,h)} \KERNEL(x,B) \tilde \mu(dx) 
		    \delta_{\rep_i^h}(A)
		  \end{equation} 
		  for all measurable sets $A \subset \{\rep_i^h:\ i \in \{1,\ldots,\num(h)\}\},\ B \subset \TYPE$. We define: 
		  \begin{equation}
		    \CutTwo^h(\mfu) = \left[\{\rep_i^h:\ i \in \{1,\ldots,\num(h)\}\},\ r, \ \mu^h\right]
		  \end{equation}
		  and
		  \begin{equation}
		    \Cut^h(\mfu) = \left[\{\rep_i^h:\ i \in \{1,\ldots,\num(h)\}\},\ r- h_{\neq}, \ \mu^h\right],
		  \end{equation}
		  where $(r-h_{\neq})(x,y) = r(x,y) - h$ if $x \neq y$ and $0$ otherwise. We call  $\Cut^h$ the {\it $h$-trunk}.
		\end{definition}
		
		\begin{remark}
		 By Remark 5.1 in \cite{decomposition}, the above definition does not depend on the choice of representative $\{\rep_i^h:\ i \in \{1,\ldots,\num(h)\}$. 
		\end{remark}

		We are now ready to give the definition of evolving genealogies and measure representations

		\begin{definition}
		  We call $(\mfu_T)_{T \ge 0} = ([[0,1],r_T,\mu_T])_{T \ge 0} \in D([0,\infty),\UMK)$ a (path of an) {\it evolving 
		  genealogy}, if for all $T \ge 0$, $0 \le h' \le h$, there is an isometry $\tilde \isom_{h',h}^T: (\supp(\tilde \mu_{T+h}^h) ,r_{T+h'})\to (\supp(\tilde \mu_{T+h'}^{h'}),r_{T+h})$, where $\Cut^h(\mfu_{T+h}) = [[0,1],r_{T+h},\mu_{T+h}^h]$. We denote by $\SEG$ the space of evolving genealogies.
		\end{definition}

		\begin{remark}
		  (1) Note that a Polish space $(X,r)$ is Borel isomorphic to either $[0,1]$, $\mathbb Q \cap [0,1]$ or a finite space (see Corollary 6.8.8 in \cite{Bog}). 
		  Moreover, when $\mu$ is a Borel-measure on $(X,r)$, and $\lambda$ is the Lebesgue measure on $[0,1]$, then 
		  there is always a measure preserving map $[0,1] \to X$ (see Theorem 8.5.4 in \cite{Bog}).  \par 
		  (2) Note that in the case where $X \subset \mathbb R$ is countable, the metric $r_{\mathbb R \setminus X}$ induced by a Borel-isomorphism $X\setminus R \cup \{x^\ast\} \to [0,1]$ (with the standard metric), where $x^\ast \in X$,  is complete and the space equipped with that metric is separable. If $(X,r_X)$ is a complete separable metric space, then the metric 
			\begin{equation}
			r(x,y) := \left\{\begin{array}{ll}
			r_X(x,y),&\quad x,y \in X \\
			r_{\mathbb R \setminus Y}(x,y),&\quad x,y \in \mathbb R \setminus X \cup\{x^\ast\} \\
			r_{X}(x,x^\ast) + r_{\mathbb R \setminus Y}(x^\ast,y),&\quad x \in X,\ y \in \mathbb R \setminus X \cup\{x^\ast\} \\
			r_{X}(x^\ast,y) + r_{\mathbb R \setminus Y}(x,x^\ast),&\quad y \in X,\ x \in \mathbb R \setminus X \cup\{x^\ast\}
			\end{array}\right.
			\end{equation}
			extends $(X,r_X)$ to a complete separable metric space $(\mathbb R,r)$. 
			One consequence is, that we can always find a representative $([0,1],r',\mu') \in [X,r,\mu]$. Moreover, the Borel-sets generated by $r'$ on $\supp(\mu')$ are a subset of the Borel-sets on $[0,1]$, i.e. we will in the following always assume that the respective measures are Borel-measures on $[0,1]$ (with the standard Borel-$\sigma$ field). 
		\end{remark}

		Next we need a suitable topology on $\mathcal M_f([0,1]\times \TYPE)$:
		
		\begin{definition}\label{def.weak.atomic}
		  A sequence of measures $\mu_1,\mu_2\ldots \in \mathcal M_f([0,1]\times \TYPE)$ converges to some measure $\mu \in \mathcal M_f([0,1]\times \TYPE)$ in the {\it (marked) weak atomic topology}, if $\mu_n \Rightarrow \mu$ in the weak topology, and 
		  \begin{equation}
		   \left(\tilde \mu_n\right)^\ast \Rightarrow \left( \tilde \mu \right)^\ast,
		  \end{equation}
		  where $\left( \tilde \mu\right)^\ast = \sum_{x \in [0,1]} \tilde \mu(\{x\})^2 \delta_x$. 
		\end{definition}
		
		\begin{remark}
		  Using Lemma 2.3 in \cite{EKatomic}, it is straight forward to see that $\mathcal M_f([0,1]\times \TYPE)$ equipped with the (marked) weak atomic topology is Polish.
		\end{remark}

		\begin{assumption}
		  When we consider processes or cadlag functions with values in $\mathcal M_f([0,1]\times \TYPE)$, we will in the following always assume that this space is equipped with the (marked) weak atomic topology.
		\end{assumption}
		
		\begin{definition}\label{def.dMR}
		  Let $(\mfu_T)_{T \ge 0} = ([[0,1],r_T,\mu_T])_{T \ge 0} \in \SEG$ be a (path of an) evolving genealogy. We call a
		  collection $\{\mfx^T: T \ge 0\}$ of cadlag processes 
		  $\mfx^T = (\mfx^T_h)_{h \ge 0} \in D([0,\infty),\mathcal M_f([0,1]\times \TYPE))$ a measure representation of the evolving genealogy $(\mfu_T)_{T \ge 0}$, if for all $T \ge 0$, there is a metric $r_T$ such that 
		  \begin{equation}\label{eq.MR.1}
		   [[0,1],r_T,\mfx^T_h] = \Cut^h(\mfu_{T+h})
		  \end{equation}
		  and 
		  \begin{equation}
		   \supp(\tilde \mfx^T_{h'}) \supset \supp(\tilde \mfx^T_h),
		  \end{equation}
		  for all $0 \le h'<h$. 
			Moreover, we define the {\it space of measure representations} by
		  \begin{equation}
		   \begin{split}
		    \SMR:=\{ (\mfx_t)_{t \ge 0}\in D([&0,\infty),\mathcal M_f([0,1]\times \TYPE)):\\
		    & \supp(\tilde \mfx_t) \subset \supp(\tilde \mfx_{t'}), 0 \le t' \le t,\\
				& \tilde \mfx_t\text{ is purely atomic for all } t > 0\}. 
		   \end{split}
		  \end{equation}
		\end{definition}

		\begin{remark}\label{rem.mr.atomic}
		  (1) Note that by Lemma \ref{lem_representatives}  a measure representation $\{\mfx^T: T \ge 0\}$ always satisfies $\mfx^T_\delta$ is purely atomic for all $\delta > 0$. \par 
		  (2) By construction $ \overline{\mfu}_T = \tilde \mfx^{h}_{T-h}([0,1])$ for all $h \in [0,T]$.
		\end{remark}

		\begin{proposition}\label{prop.MR.EG}
		  Let $(\mfu_T)_{T \ge 0}  \in D([0,\infty),\UMK)$. Then $(\mfu_T)_{T \ge 0} \in \SEG$ implies the existence of a measure representation and the existence 
		  of a collection of cadlag functions $\{\mfx^T: T \ge 0\}$ with $\mfx^T \in \SMR$ for all $T \ge 0$ and with \eqref{eq.MR.1} implies that $(\mfu_T)_{T \ge 0} \in \SEG$. 
		\end{proposition}

		\begin{theorem}\label{thm.compact}
		  Write $(\mfx^{\mfu,T})_{T \ge 0} \in \SMR$ for a measure representation of an evolving genealogy $\mfu := (\mfu_T)_{T \ge 0} \in \SEG$ and take $\Gamma \subset \SEG$. We consider the following conditions: 
		  \begin{itemize}
		    \item[(i)] $\{\mfu_0: (\mfu_T)_{T \ge 0} \in \Gamma\}$ is relatively compact in $\UMK$. 
		    \item[(i')] There is a $C \ge 0$ such that $\nu^{2,\mfu}([C,\infty)) = 0$ for all $\mfu \in \Gamma$. 
		    \item[(ii)] For all $T \ge 0$ and all $\eps > 0$ there is a $H \ge 0$ such that 
		      \begin{equation}
		      \sup_{\mfu \in \Gamma} \left( \left(\int_{[0,1]} \tilde \mfx^{\mfu,T}_H(dx)\right)^2 - \int_{[0,1]} \tilde \mfx^{\mfu,T}_H(\{x\}) \tilde \mfx^{\mfu,T}_H(dx)\right) \le \eps.
		      \end{equation}
		    \item[(iii)] For all $T \ge 0$, there is a compact set 
		      \begin{equation}
		      \begin{split}
			A = A^T\subset \{(\mfx_h)_{h \ge 0} &\in D([0,\infty),\mathcal M_f([0,1]\times \TYPE):\\
			&\tilde \mfx_\delta \text{ is purely atomic for all } \delta > 0\}
		      \end{split}
		      \end{equation}
		      (equipped with the subspace topology), such that $(\mfx^{\mfu,T}_h)_{h \ge 0} \in A$ for all $\mfu \in \Gamma$. 
		    \item[(iii')] For all $T \ge 0$, there is a compact set $A \subset \SMR$ (equipped with the subspace topology), such that $(\mfx^{\mfu,T}_h)_{h \ge 0} \in A$ for all $\mfu \in \Gamma$.
		  \end{itemize}
		  Then the following holds: 
		  \begin{itemize}
		   \item[(1)] Condition (i),(i'),(ii),(iii) imply $\Gamma$ is relatively compact in $D([0,\infty),\UMK)$.
		   \item[(2)] Condition (i),(ii),(iii') imply $\Gamma$ is relatively compact in $\SEG$ (equipped with the subspace topology). Moreover, let $(\mfu^n)_{n \in \mathbb N}$ be a sequence in $\Gamma$
		    with measure representation $\{(\mfx^{n,T}_h)_{h \ge 0}: T \ge 0\}_{n \in \mathbb N}$ and let $D \subset [0,\infty)$ be dense. 
		    Then, convergence of $(\mfx^{n,T}_h)_{h \ge 0} \to (\mfx^{T}_h)_{h \ge 0} \in \SMR$ for all $T \in D$ implies that all limit points of $(\mfu^n)_{n \in \mathbb N}$ have some measure 
		    representation $\{(\hat \mfx^{T}_h)_{h \ge 0}: T \ge 0\}$ with the property that $(\hat \mfx^T_h)_{h \ge 0} =  (\mfx^T_h)_{h \ge 0} $ for all but at most countable many $T \ge 0$. \par
		   \item[(3)] Let $(\mfu^n)_{n \in \mathbb N}$ be a sequence in $\Gamma$ and assume (i),(ii) and (iii'). 
			Moreover, we assume that $(\mfx^{n,T}_h)_{h \ge 0} \to (\mfx^{T}_h)_{h \ge 0} \in \SMR$ for all $T \in D$ as in (2). If 
			for all $T,\delta \in D$, one has 
			\begin{equation}\label{eq.uniq.mvr}
			\sum_{x \in A_1} \tilde \mfx^{T}_\delta(\{x\}) \neq \sum_{x \in A_2} \tilde \mfx^{T}_\delta(\{x\}),
			\end{equation}
			where $A_1,A_2 \subset \{x \in [0,1]:\ \tilde \mfx^{T}_\delta(\{x\}) > 0\}$ with $A_1 \neq A_2$, then uniqueness of the limits in (2) holds, i.e. $(\mfu^n_T)_{T \ge 0} \rightarrow (\mfu_T)_{T \ge 0} \in \SEG$. 
		  \end{itemize}
		\end{theorem} 

		\begin{remark}
		  (1) Note that (iii') implies (iii). \\
		  (2) Note that in the case of purely atomic measures $\mu$ on $[0,1]$, 
		  \begin{equation}
		    \sum_{x \in [0,1]} \mu(\{x\})^2 = \int_{[0,1]} \mu(\{x\}) \mu(dx). 
		  \end{equation} 
		\end{remark}

		\begin{proposition}\label{prop.rel.comp.SMR}
		  A subset $\Gamma \subset \SMR$ is relatively compact (equipped with the subspace topology), if and only if it is relatively compact in $D([0,\infty),\mathcal M_f([0,1]\times \TYPE))$ and: 
		  \begin{itemize}
		    \item[(i)] For all  $T,\bar \delta > 0$ there is a $C \ge  0$ such that 
		      \begin{equation}
		       \sup_{\mfx \in \Gamma} \inf \{C' \ge 0: \sup_{t \in[\delta,T]} \tilde\mfx_t \le C' \tilde\mfx_\delta \text{ for all } \delta \in [\bar \delta,T]\} \le C,
		      \end{equation}
		      where we say $\mu \le \nu$ for two measures $\mu,\nu$ if $\mu(A) \le \nu(A)$ for all measurable $A$. 
		   \item[(ii)] $\{(\ord(\tilde \mfx_h)):\ \mfx \in \Gamma\}$, is relatively compact in $\mathcal S^\downarrow$ for all $h > 0$, where $\ord(\mu)$ is the non increasing reordering of sizes of atoms of 
		      a measure $\mu \in\mathcal M_f([0,1])$ and 
			\begin{align}
			  \mathcal S^\downarrow&:=\left\{(x_1,x_2,\ldots) \in [0,\infty)^{\mathbb N}:\ 
			  \sum_{i \in \mathbb N} x_i < \infty,\ x_1 \ge x_2 \ge \ldots \right\},
			\end{align}
			is equipped with
			\begin{equation}\label{eq_metric_S}
			  \dSeq(x,y) = \sum_{i = 1}^\infty |x_i - y_i| = \norm{x - y}_1.
			\end{equation}
		  \end{itemize}
		  Moreover, using the notation of Theorem \ref{thm.compact}, relative compactness of a set $\Gamma \subset \SEG$ in $\UMK$, implies (ii) for the corresponding measure representations. 
		\end{proposition}
		
		\begin{definition}
		  We call a stochastic process $\U = (\U_T)_{T \ge 0}$ with values in $\UMK$ an {\it evolving genealogy}, if $P(\U \in \SEG) = 1$. \par
		  Let $\{(\hat {\mathcal X}^T_h)_{h \ge 0}:\ T \ge 0\}$ and $\{({\mathcal X}^T_h)_{h \ge 0}:\ T \ge 0\}$ be collections of processes with $P((\hat {\mathcal X}^T_h)_{h \ge 0} \in \SMR) = P(({\mathcal X}^T_h)_{h \ge 0} \in \SMR) = 1$ for all $T\ge 0$. We call 
		  \begin{itemize}
		    \item[(i)] $\{(\hat {\mathcal X}^T_h)_{h \ge 0}:\ T \ge 0\}$ a  {\it strong measure representation} for $\U$, if both are defined on the same probability space and for all $T \in D$, where $D \subset [0,\infty)$ is dense the following holds: Almost surely we find a representative, such that \eqref{eq.MR.1} holds. 
		    \item[(ii)] $\{({\mathcal X}^T_h)_{h \ge 0}:\ T \ge 0\}$ a  {\it measure representation} for $\U$, if $\{({\mathcal X}^T_h)_{h \ge 0}:\ T \ge 0\} \stackrel{\text{f.d.d.}}{=}{\{(\hat {\mathcal X}^T_h)_{h \ge 0}:\ T \ge 0\}}$, where $\stackrel{\text{f.d.d.}}{=}$ means equality of the finite dimensional distributions and $\{(\hat {\mathcal X}^T_h)_{h \ge 0}:\ T \ge 0\}$ is a strong representation defined in (i).
		    \item[(iii)] $\{({\mathcal X}^T_h)_{h \ge 0}:\ T \ge 0\}$ a  {\it weak measure representation} for $\U$, if $({\mathcal X}^T_h)_{h \ge 0} \stackrel{\text{d}}{=} (\hat {\mathcal X}^T_h)_{h\ge 0}$, for all $T \ge 0$. 
		  \end{itemize}
		\end{definition}

		\begin{remark}
		  (1) Note that by Proposition \ref{prop.MR.EG}, an evolving genealogy always has a strong measure representation. \par 
		  (2) We do not want to discuss measurability questions here. If a reader feels uncomfortable, one should have in mind that we could also work with modifications and replace the ``almost surely'' by a ``surely''. This modifications would not change the results. 
		\end{remark}
		
		\begin{remark}\label{rem.MCRA1}
		Let  $\{({\mathcal X}^T_h)_{h \ge 0}:\ T \ge 0\}$ be a weak measure representation of some evolving genealogy $\U$, where we assume for simplicity that $|\TYPE| = 1$. We want to sketch how one can use the measure representation to obtain genealogical information and there are two applications that follow directly by construction:  
		\begin{itemize}
		\item[(i)] {\it (Time to the MRCA)} One can use the measure representation to calculate the time to the most recent common ancestor. Namely, the absorption time
		\begin{equation}
		\tau:= \inf\{h > 0: \mathcal X^T_h \in \{\delta_x: x\in [0,1]\}\}
		\end{equation}
		is exactly the first time such that the population, described by $\U_{T+\tau}$, can be covered by exactly one ball, which in terms of genealogical distance means that this is the first time such that there is one single ancestor that gave birth to the whole population. \par 
		As an example, we will show how one can calculate the expected time to the most recent common ancestor in equilibrium of a Fleming-Viot population (see Remark \ref{rem.MCRA2} below). 
		\item[(ii)] {\it (Distance of two randomly chosen individuals)} Another application is to calculate the typical time to the most recent common ancestor of two randomly chosen individuals, i.e. the first moment measure $E[\nu^{2,\U_T}(\cdot)]$. In terms of the measure representation, this measure is given by
		\begin{equation}
			E[\nu^{2,\U_T}([0,h])] = \sum_{x \in [0,1]} E[\mathcal X^T_h(\{x\})^2]. 
		\end{equation}
		Again, we give an example in Remark \ref{rem.MCRA2}. 
		\end{itemize}
		\end{remark}

		We use the notion measure representation and evolving genealogies in both contexts; stochastic and deterministic. 

		\begin{theorem}\label{thm.tightness}
		  Let $\U^n = (\U_T^n)_{T \ge 0}$, $n = 1,2,\ldots$ be a sequence of evolving genealogies with 
		  \underline{weak} measure representation $\{(\mathcal X^{n,T}_h)_{h \ge 0}:\ T \ge 0\}$, $n = 1,2,\ldots$. Then $(\U^n)_{n \in \mathbb N}$ is tight if
		  \begin{itemize}
		    \item[(i)] $(\U^n_0)_{n \in \mathbb N}$ is tight. 
		    \item[(ii)] For all $T \ge 0$ and all $\eps > 0$, there is a compact set 
		      \begin{equation}
			\begin{split}
			  A \subset \{(\mfx_h)_{h \ge 0} &\in D([0,\infty),\mathcal M_f([0,1]\times \TYPE):\\
			  &\tilde \mfx_\delta \text{ is purely atomic for all } \delta > 0\}
			\end{split}
		      \end{equation}
		      (equipped with the subspace topology), such that 
		      \begin{equation}
			\limsup_{n \rightarrow \infty} P\left(({\mathcal X}^{n,T}_h)_{h \ge 0} \in A^c \right) \le \eps.
		      \end{equation}
		    \item[(iii)] For all $T \ge 0$ and all $\eps > 0$ there is a $H \ge 0$ such that 
		      \begin{equation}
			\limsup_{n \rightarrow \infty} P\left( \left(\int_{[0,1]} \tilde {\mathcal X}^{n,T}_H(dx)\right)^2 - \int_{[0,1]} \tilde {\mathcal X}^{n,T}_h(\{x\}) \tilde {\mathcal X}^{n,T}_H(dx)\ge \eps \right) \le \eps.
		      \end{equation}
		  \end{itemize}
		\end{theorem}
		
		\begin{remark}\label{rem.cond.ii}
		 Note that condition (ii) is satisfied if $\mathcal X^{n,T}$ converges for all $T \ge 0$  to a measure-valued process 
		 $\mathcal X^{T}$ with the property that  $\tilde {\mathcal X}^{n,T}_\delta$ is purely atomic for all $\delta > 0$. 
		\end{remark}

		\begin{theorem}\label{thm.convergence.II}
		  Assume that $\U^n = (\U_T^n)_{T \ge 0}$, $n = 1,2,\ldots$ is a tight sequence of evolving genealogies with measure representation $\{(\mathcal X^{n,T}_h)_{h \ge 0}:\ T \ge 0\}$. 
		  Assume further, that 
		  \begin{itemize}
		    \item[(i)] $\U_0^n \Rightarrow \U_0$,  where $\U_0 = [\{0\},0,\mu]$ and $\mu$ is a random measure with values in $\mathcal M_f(\{0\} \times \TYPE)$.   
		    \item[(ii)] For all $T \ge 0$, the processes $(\mathcal X^{n,T}_h)_{h \ge 0}$ are tight. 
		    \item[(iii)] For all $\tau > 0$ and $\delta > 0$ and all $\eps > 0$, there is a $C \in \mathbb N$ such that 
		    \begin{equation}
			\begin{split}
			  \limsup_{n \rightarrow \infty} P\Big(\exists l =0,\ldots,&\lceil \tau/\delta \rceil, \exists  A \subset [0,1] \text{ measurable}:\\
				&\sup_{h \in [\delta,2\delta)} \tilde{\mathcal X}^{n,l\cdot \delta}_h(A) \ge C \tilde{\mathcal X}^{n,l\cdot \delta}_\delta(A)\Big) \le \eps.
			\end{split}
		    \end{equation}    
		    \item[(iv)] There is a collection $\{(\mathcal X^T_h)_{h \ge 0}:\ T \ge 0\}$, with $P((\mathcal X^T_h)_{h \ge 0} \in \SMR) = 1$ for all $T \ge 0$ and 
		      \begin{equation}
			\begin{split}
			  \left(\mathcal X^{n,T+(K-k)\cdot h}_{l\cdot h}\right)_{1 \le l \le k \le K}
			  \Rightarrow \left(\mathcal X^{T+(K-k)\cdot h}_{l\cdot h}\right)_{1 \le l \le k \le K}
			\end{split}
		      \end{equation}
		      for all $K \in \mathbb N$, $h,T  \in D$ where $D \subset (0,\infty)$ is  dense. 
		    \item[(v)] Let $\{(\mathcal X^T_h)_{h \ge 0}:\ T \ge 0\}$ be the limit given in (iv), let $\delta > 0$ and denote by 
		      $a:= a_\delta := \ord(\tilde {\mathcal X}^T_\delta)$ the reordered sizes of atoms of $\tilde{\mathcal X}^T_\delta$ (recall Proposition \ref{prop.rel.comp.SMR}). 
		      For all $\delta > 0$ the following holds: Given $\tilde{\mathcal X}^T_\delta$ has $K \in \mathbb N\cup \infty$ many atoms, then 
		      $\hat a := \left(\frac{a_j}{\sum_i a_i}\right)_{j = 1,\ldots,K}$ is absolutely continuous to the Lebesgue measure $\lambda|_{[0,1]}^K$. 
		  \end{itemize}
		  Then there is an evolving genealogy $\U$ such that $\U^N \Rightarrow \U$ as processes and $\U$ has measure representation $\{(\mathcal X^T_h)_{h \ge 0}:\ T \ge 0\}$.
		\end{theorem}

		\begin{remark}\label{rem.conv.fdd}
		In order to prove convergence of evolving genealogies as processes in the Skorohod space, one needs to prove tightness and convergence of the finite dimensional distributions 
		(see Theorem 3.7.8 in \cite{EK86}). Assume we are in the situation without marks, i.e. $|\TYPE| = 1$. While tightness can be proven using Theorem \ref{thm.tightness}, 
		the f.d.d convergence follows, when condition (v) of Theorem \ref{thm.convergence.II} holds and the reordered sizes of atoms converge, i.e. 
		      \begin{equation}
			\begin{split}
			  \left(\ord\left(\mathcal X^{n,T+(K-k)\cdot h}_{l\cdot h}\right)\right)_{1 \le l \le k \le K}
			  \Rightarrow \left(\ord\left(\mathcal X^{T+(K-k)\cdot h}_{l\cdot h}\right)\right)_{1 \le l \le k \le K}.
			\end{split}
		      \end{equation}
		This is part of the proof of Theorem \ref{thm.convergence.II}. 
		\end{remark}
		
		\begin{remark}
		One could generalize the initial condition in (i) to so called identifiable (by family sizes) ultra-metric measure spaces, $\UMI$. This space was studied in \cite{decomposition} (see also Lemma \ref{lem.res.fam} below).
		\end{remark}

\section{Application - A finite system scheme result for the genealogy in a spatial Fleming-Viot population}\label{sec.applicationFV}

		Here we apply our theory to the tree-valued interacting Fleming-Viot process. Our goal is to study the behavior of the genealogy in a spatial Fleming-Viot population when the size of the geographical space tends to infinity.   \par 
		In section \ref{sec.TVIMM} we introduce the tree-valued interacting Moran model. We show that this model is an evolving genealogy. In section \ref{sec.TVIFV} we show that its limit, the tree-valued interacting Fleming-Viot process, is again an evolving genealogy and give the main result on the behavior of the genealogy in a spatial Fleming-Viot population.

    \subsection{Tree-valued interacting Moran models}\label{sec.TVIMM}
    
		We start by defining the basic model and then use a graphical construction to define the corresponding tree-valued interacting Moran models, where we follow the approach in \cite{DGP12} and \cite{GPW13}. \\

		\noindent {\bf The model:} We are in the following situation: We want to describe a population that lives in a geographical space $G$, 
		    where we assume that 
		    \begin{equation}\label{eq.geoSpace}
			    G = \{g_1,\ldots,g_m\}\textrm{ is a finite abelian group}.
		    \end{equation}
		    
		    For a fixed $N \in \mathbb N$, we assume that our population consists of $N\cdot |G|$ individuals, where  
		    the initial spatial configuration of the individuals is independent uniformly on $G$, i.e.
		    \begin{equation}
			    (\zeta_i(0))_{i \in\{1,2,\ldots,N\cdot|G|\}}\text{ are i.i.d. with } P(\zeta_1(0) = g) = \frac{1}{|G|},\quad \forall g \in G.
		    \end{equation}
		    The population evolves according to the following dynamics:
		    \begin{itemize}
		      \item[(1)] {\it Resampling:} Every pair $i \not=j$, which is located at the same site, is replaced with the resampling rate $\gamma > 0$.
				    If such an event occurs, $i$ is replaced by an offspring of $j$, or $j$ is replaced by an
				    offspring of $i$, each with probability $\frac{1}{2}$ if their locations coincide.
		      \item[(2)] {\it Migration:} Every individual migrates (independently) according to a random walk kernel $a(\cdot,\cdot)$ on $G$, where we assume for all $\xi,\xi' \in G$,
				    \begin{equation}
					    a(\xi,\xi') \in [0,1],\quad a(\xi,\xi') = a(0,\xi-\xi'), \quad \sum_{\xi \in G} a(0,\xi) = 1.
				    \end{equation}
		    \end{itemize}

		\noindent {\bf Graphical construction:} For $N \in \mathbb N$ and $G$ as above we set 
		    \begin{equation}
			    I_N := \{1,2,\ldots,N\cdot |G|\}. 
		    \end{equation}
		    
		    Let $\zeta = (\zeta(t))_{t \ge 0}$  be a continuous time random walk on $G$ with transition rate $a(\zeta,\zeta')$  and $\{\zeta_k\}_{k \in I_N}$ be a 
		    family of independent copies of $\zeta$, where we assume that $(\zeta_i(0))_{i \in I_N}$ are independent and uniformly distributed on $G$. \\
		    
		    Let $\{\xi_k:\ k \in I_N\}$be a realization of  $\{\zeta_k \}_{k \in I_N}$ and let $\{\pois^{i,j}:\ i,j \in I_N,\ i \not=j\}$
		    be a realization of a family of independent rate $\frac{\gamma}{2}$ Poisson processes, defined on the same probability 
		    space as $\{\zeta_k\}_{k \in I_N}$, where we assume that both are independent.

		    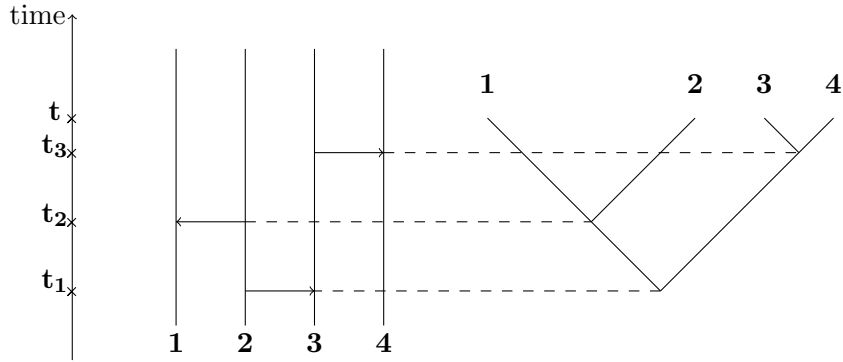
\begin{figure}[ht]
		      \centering
		      \begin{tikzpicture}[scale = 0.46]
			\draw (1,4)-- (1,-4);
			\draw (3,4)-- (3,-4);
			\draw (5,4)-- (5,-4);
			\draw (7,4)-- (7,-4);
			\draw [->] (-2,-5) -- (-2,5);
			\draw (1,-4.5) node {$\mathbf{1}$};
			\draw (3,-4.5) node {$\mathbf{2}$};
			\draw (5,-4.5) node {$\mathbf{3}$};
			\draw (7,-4.5) node {$\mathbf{4}$};
			\draw (10,2)-- (15,-3);
			\draw (15,-3)-- (20,2);
			\draw (19,1)-- (18,2);
			\draw (13,-1)-- (16,2);
			\draw (10,3) node {$\mathbf{1}$};
			\draw (16,3) node {$\mathbf{2}$};
			\draw (18,3) node {$\mathbf{3}$};
			\draw (20,3) node {$\mathbf{4}$};
			\draw [dash pattern=on 4pt off 4pt] (7,1)-- (19,1);
			\draw [dash pattern=on 4pt off 4pt] (15,-3)-- (5,-3);
			\draw [->] (5,1) -- (7,1);
			\draw (-2.5,2.2) node {$\mathbf{t}$};
			\draw (-2.5,-2.7) node {$\mathbf{t_1}$};
			\draw (-2.5,-0.8) node {$\mathbf{t_2}$};
			\draw (-2.5,1.2) node {$\mathbf{t_3}$};
			\draw (-3,5) node {$\text{time}$};
			\draw [->] (3,-1) -- (1,-1);
			\draw [->] (3,-3) -- (5,-3);
			\draw [dash pattern=on 4pt off 4pt] (3,-1)-- (13,-1);
			\draw [color=black] (-2,2)-- ++(-4pt,-4pt) -- ++(7pt,7pt) ++(-7pt,0) -- ++(7pt,-7pt);
			\draw [color=black] (-2,1)-- ++(-4pt,-4pt) -- ++(7pt,7pt) ++(-7pt,0) -- ++(7pt,-7pt);
			\draw [color=black] (-2,-1)-- ++(-4pt,-4pt) -- ++(7pt,7pt) ++(-7pt,0) -- ++(7pt,-7pt);
			\draw [color=black] (-2,-3)-- ++(-4pt,-4pt) -- ++(7pt,7pt) ++(-7pt,0) -- ++(7pt,-7pt);
		      \end{tikzpicture}
		      \caption{\label{f.graph.constr} \footnotesize Graphical construction of the TVIMM (i.e. the tree on the right side), where we assumed for simplicity that $|G| = 1$, i.e. all individuals are located on one single site, and $N = 4$. At times $t_1, t_2, t_3$ we sample the individuals $(x_1^1,x_2^1) =(2,3)$, $(x_1^2,x_2^2)=(2,1)$, $(x_1^3,x_2^3)=(3,4)$ and draw an arrow from $x_1^j$ to $x_2^j$. At time $t$ the ancestors $A_h(i,t) $, $i =1,\ldots,4$ at time $h \in (t_1,t_2]$, for example,  are $A_h(1,t) = 2$, $A_h(2,t) = 2$, $A_h(3,t) = 3$ and $A_h(4,t) = 3$. In the case where $|G| = 1$, $\mu^N_t$ is the uniform distribution on $\{1,2,3,4\}$.}
		    \end{figure}
				
		    \begin{definition}
		      Let $i,i' \in I_N$, $0\le h\le t < \infty$ we say that there is a {\it path} from $(i,h)$ to $(i',t)$ if there is a  $n \in \mathbb N$, 
		      $h \le u_1 < u_2 < \cdots < u_n \le t$ and
		      $j_1,\ldots,j_n \in I_N$ such that for all $k \in \{1,\ldots,n+1\}$ ($j_0:= i, j_{n+1}:=i'$)
		      $\pois^{j_{k-1},j_k}\{u_k\} = 1$, $\xi_{j_{k-1}}(u_k)= \xi(u_k)_{j_{k}}$ and $\pois^{x,j_{k-1}}\{s\}= 0$ for all $x \in I_N$ 
		      with $\xi_x(s) = \xi_{j_{k-1}}(s)$, $s \in (u_{k-1},u_k)$. 
		    \end{definition}

		    Note that for all $i \in I_N$ and $0 \le h \le t$ there exists an unique element
		    \begin{equation}
		    A_h(i,t) \in I_N
		    \end{equation}    
		    with the property that there is a path from $(A_h(i,t),h)$ to $(i,t)$. We call $A_h(i,t)$ the 
		    {\it ancestor of  $(i,t)$ at time $h$} (measured backwards); see Figure \ref{f.graph.constr}. \\
				
		    Let $r_0$ be a pseudo-ultra-metric on $I_N$. We define the pseudo-ultra-metric ($i,j \in I_N$):			
		    \begin{equation}\label{eq.metric.Moran}
			    r_t(i,j):=\left\{ \begin{array}{ll}
					    t-\sup\{h \in [0,t]:A_h(i,t)=A_h(j,t)\},& \textrm{ if } A_0(i,t) = A_0(j,t),\\[0.2cm]
					    t + r_0(A_0(i,t), A_0(j,t)) ,& \textrm{ if } A_0(i,t) \not= A_0(j,t).
			    \end{array}\right.
		    \end{equation}

		    For $t \ge 0$, we define  $\mu^N_t\in \M_f(I_N\times G)$ by
		    \begin{equation}\label{eq.MM.measure}
			    \mu^N_t(A\times\{g\}) = \frac{1}{|\{i:\ \xi_i(t) = g\}| \wedge 1} \sum_{k \in A}  1(\xi_k(t) = g),  
		    \end{equation}		
				$A \subset I_N,\ g \in G$.

		    Now, since $r_t$ is only a pseudo-metric, we  consider the following equivalence relation $\approx_t$  on $I_N$: 
		    \begin{equation}
			    x \approx_t y \Leftrightarrow r_t(x,y) = 0.
		    \end{equation}
		    
		    We define the set $\tilde I_N^t:= I_N\! /\!\!\approx_t$ of equivalence classes and note that we can find a set of representatives
		    $\overline{I}_N^t$ such that $\overline{I}_N^t \to \tilde I_N^t$, $x \mapsto [x]_{\approx_t}$ is a bijection. Let $\bar i,\bar j \in \bar I_N^t, \ g \in G$ and define
		    \begin{align}\label{eq.def.moran}
			    \bar r_t(\bar  i,\bar  j) =  r_t(\bar i,\bar j), \quad 
			    \bar \mu^N_t(\{\bar i\} \times \{g\}) =\mu^N([\bar i]_{\approx_t} \times \{g\}).
		    \end{align}

		    Then the {\it tree-valued interacting Moran model} (TVIMM) of size $N \in \mathbb N$ is defined as 
		    \begin{equation}
			    \U_t^N :=  [\bar I_N^t,\bar  r_t, \bar  \mu^N_t].
		    \end{equation}

		    \begin{remark}
		    In the situation where $|G| = 1$  we can identify $\U^N$ with an $\UM$-valued process and we call this process {\it (non spatial) tree-valued Moran model}. 
		    \end{remark}
	    
		    \begin{assumption}\label{A.init.metric}
		      In the following we will always assume that 
		      \begin{equation}
			r_0 \equiv 0,
		      \end{equation}
		      i.e. at time $0$ all individuals are related.
		    \end{assumption}
		    
		    \begin{remark}\label{r.init.IMM}
		      Note that this assumptions implies
		      \begin{equation}
			\U^N_0 = [\{1\},0,\nu],
		      \end{equation}
		      with $\nu(\{1\} \times B) = |B|$ for all $B \subset G$.
		    \end{remark}
		    
		    As the main result of this section we have the following: 
		    
		    \begin{theorem}\label{thm.TVIMMM.EG}
		      The tree-valued interacting Moran model $\U^N$ is an evolving genealogy. 
		    \end{theorem}

		    In order to work with this result it is necessary to define a suitable measure representation. We will do this in the proof section. But to get an idea, we remark the following: 
		    
		    \begin{remark}
		      The measure representation $\{(\mathcal X^{N,T}_h)_{h \ge 0}:\ T\ge 0\}$ can be described as coupled family of measure-valued interacting Moran models,  
		      $\mathcal X^{N,T}$, for all $T \ge 0$ and the coupling for different $T$ is given in terms of spatial Kingman coalescents. 
		    \end{remark}

	\subsection{A finite system scheme result for the tree-valued interacting Fleming-Viot processes}\label{sec.TVIFV}
		
			In order to analyze properties of a finite (but large) population, it sometimes is useful to study an ``infinite'' population model. 
			In the case of measure-valued interacting Moran models the corresponding process is the so called system of measure-valued interacting Fleming-Viot processes (or measure-valued interacting Fleming-Viot process) and in the case of tree-valued interacting Moran models the resulting process is called the tree-valued interacting Fleming-Viot process (TVIFV): 

			\begin{theorem}\label{t.convMM} Let $\UMG$ be equipped with the Gromov-weak topology. Let $\U^{N}$ be the tree-valued interacting Moran model on the finite geographical space $G$ defined in Section \ref{sec.TVIMM}. Then
			\begin{equation}
			(\U^{N}_t)_{t \ge 0} \stackrel{N \rightarrow \infty}{\Rightarrow} (\U_t)_{t \ge 0}
			\end{equation}
			weakly in the Skorohod topology on $D_{\UMG}([0,\infty))$, where $\L(\U_0) = [\{1\},0,\nu]$ (recall Remark \ref{r.init.IMM}) and $(\U_t)_{t \ge 0}$ is an evolving genealogy.
			\end{theorem}
			
			\begin{remark}\label{rem.MCRA2}
			  Assume $|G| = 1$.  As we will see in the proof section (see Remark \ref{rem.mvfv.MVR}), the measure-valued Fleming-Viot process that starts in the Lebesgue measure on 
			  $[0,1]$ is a weak measure representation for the tree-valued Fleming-Viot process. Having this in mind, we continue Remark \ref{rem.MCRA1}.
			  \begin{itemize}
				  \item[(i)] {\it (Time to the MRCA)} One can show that the tree-valued Fleming-Viot process converges as $T \rightarrow \infty$ to its unique equilibrium $\U_\infty$ (see \cite{GPW13}). When we now want to calculate the expected time to the most recent common ancestor, $E[T_{\text{MRCA}}]$, in this equilibrium we can for example calculate the expected time it takes for a Kingman coalescent to have just one partition element left (this is classical). Since the time it takes for two partition elements to coalesce, given there are $n$ partition elements, is exponentially distributed with rate $\binom{n}{2}$, the expected time to the most recent common ancestor is given by 
			  \begin{equation}
			  E[T_{\text{MRCA}}] = \sum_{i = 2}^\infty \frac{1}{\binom{i}{2}} = 2 \sum_{i = 2}^\infty  \frac{1}{i(i-1)} = 2.
			  \end{equation}
			  On the other hand, when we apply our idea to this question, we need to calculate the mean absorption time of a $n$-type Wright-Fisher diffusion that starts in $(1/n,\ldots,1/n)$ and then we need to take the limit $n \to \infty$. The mean absorption time can be calculated using standard techniques for diffusions (see for example Section 3.4.1.3 in \cite{baxter} - with $M = n$ and $r= n-1$) and is given by 
			  \begin{equation}
			  E[T_{\text{MRCA}}] = \lim_{n \to \infty} -2 \binom{n}{n-1} \frac{n-1}{n}\ln\left(\frac{n-1}{n}\right) = 2. 
			  \end{equation}
			  \item[(ii)] {\it (Distance of two randomly chosen individuals)} Similar to (i), we go to the equilibrium. Then, 
			  \begin{equation}
				  E[\nu^{2,\U_\infty}([0,h])] = \sum_{x \in [0,1]} E[\mathcal X^\infty_h(\{x\})^2] = \lim_{n \rightarrow \infty} \sum_{i = 1}^n E[X_h(i)^2], 
			  \end{equation}
			  where $X_h$ is a $n$-type Wright-Fisher diffusion that starts in $(1/n,\ldots,1/n)$. It is not hard to see that $h \mapsto  \sum_{i = 1}^n E[X_h(i)^2]$ is given in terms of an ODE which can be solved (see \cite{selection}): 
			  \begin{equation}
			  \sum_{i = 1}^n E[X_h(i)^2] = 1-e^{-h}+\frac{1}{n} e^{-h}. 
			  \end{equation}
			  Hence the first moment measure 
			  \begin{equation}
				  E[\nu^{2,\U_\infty}([0,h])] = \lim_{n \rightarrow \infty} \sum_{i = 1}^n E[X_h(i)^2] = 1-e^{-h}
			  \end{equation}
			  is the exponential distribution, which is exactly what one would expect in view of the Kingman coalescent. 
			  \end{itemize}
			  It is possible to include marks, i.e. it is for example possible to include selection. 
			  The approach works similar and we refer to \cite{selection} to get an idea (the proof of the main result there, Theorem 3.1, contains an error, but to learn more about selection and the corresponding measure representation, this paper is the right reference). 
			\end{remark}
			
			Now, we are in the situation where the geographical space $G$ is large (i.e. $|G| \rightarrow \infty$) and the expected meeting time of two random walks (with respect to $a(\cdot,\cdot)$) is also large (i.e. $a$ is ``almost'' transient). In this situation the genealogical distance of two randomly chosen individuals  grows to infinity and the goal is on the one hand to identify the rate of divergence and on the other to determine how the genealogy looks like in this critical time scale. Note that for a fixed finite space $G$ the kernel is recurrent and hence, this question is related to the question in which time scale the finite systems notices that it is finite. This comparison of a large finite and an infinite system is called {\it finite system scheme}. We want to make this more precise: \\
			Let $G = \mathbb Z^d$ and $G_N:= [-N,N)^d \cap \mathbb Z^d$, $N \in \mathbb N$. Moreover, we assume that the migration kernels $a^N(\cdot,\cdot)$ are given by
			\begin{equation}\label{random_walk_kernel_finite}
			  a^N(i,j) = \sum_{k \in \mathbb Z^d} a(i,j+2N k), \qquad i,j \in G_N,
			\end{equation} 
			where $j+2N:=(j_1+2N,\ldots,j_d+2N)$ and $a(\cdot,\cdot)$ is a transient migration kernel on $G$ with the properties that for all $\xi,\xi' \in G$
			\begin{equation}
			  a(\xi,\xi') \in [0,1],\quad a(\xi,\xi') = a(0,\xi-\xi'), \quad \sum_{\xi \in G} a(0,\xi) = 1,
			\end{equation}
			and
			\begin{equation}
			  \sum_{n \in \mathbb N}(a^{(n)}(0,\xi) + a^{(n)}(\xi,0)) > 0,\quad  \sum_{\xi \in G} |\xi|^{d+2} a(0,\xi) < \infty.
			\end{equation}
			We denote by $(\U^N_T)_{T \ge 0}$ the tree-valued interacting Fleming-Viot processes on the geographical spaces $G_N$ (defined in Theorem \ref{t.convMM}). \par 
			
			Observe that given two independent continuous time random walks $Z_1(t),\ Z_2(t)$ on 
			$G$ with transition kernel $ a(\cdot,\cdot)$ the distance process $(Z_1(t) - Z_2(t))_{t \ge 0}$ is a random walk on $G$ with transition kernel $\hat a(\cdot,\cdot)$ and jump rate $2$, where for $\xi,\xi' \in G$:
			\begin{equation} \label{eq_distkernel}
			  \hat a(\xi,\xi') = \frac{1}{2}(a(\xi,\xi') + a(\xi',\xi)).
			\end{equation} 
			Define 
			\begin{equation}\label{dif_const}
			  D = \frac{\gamma}{1+\gamma \int_0^\infty \hat a_{2s}(0,0)ds},
			\end{equation}
			where for $\xi,\xi' \in G$, $\hat a_t(\xi,\xi')$ is given by
			\begin{equation}
			  \hat a_t(\xi,\xi') = e^{-t}\sum_{k = 0}^\infty \frac{t^k}{k!}\hat a^{(k)}(\xi,\xi').
			\end{equation}
			We consider the following functions
			\begin{equation}
			\begin{split}
			  h_N :\UM^{G_N} &\rightarrow \UM:\ [X,r,\mu] \mapsto  \left[X,\frac{1}{|G_N|} r,\frac{1}{|G_N|}\mu(\cdot \times G_N)\right] 
				\end{split}
			\end{equation}
			and 
			\begin{equation}
			  \theta_N :\mathcal M_{|G_N|}([0,1] \times G_N) \rightarrow \mathcal M_1([0,1]), \quad \mu \mapsto \frac{1}{|G_N|} \mu(\cdot \times G_N),
			\end{equation}
			where $\mu \in \mathcal M_{|G_N|}([0,1] \times G_N)\ : \Leftrightarrow \ \mu([0,1]\times G_N) = |G_N|$.  
			\begin{remark}
				Note that 
				\begin{equation}
					\mu(\cdot \times G_N) = \sum_{g \in G_N} \mu(\cdot \times \{g\}),
				\end{equation}
				i.e. 
				\begin{align}
					h_N([X,r,\mu]) &= \left[X,\frac{1}{|G_N|} r,\frac{1}{|G_N|}\sum_{g \in G_N}\mu(\cdot \times\{g\})\right],\\
					\theta_N(\mu) &= \frac{1}{|G_N|} \sum_{g \in G_N} \mu(\cdot \times \{g\}).
				\end{align}
			\end{remark}
			We get the following finite system scheme result from a global perspective: 

			\begin{theorem}\label{thm.FSS}
			  Assume that $\hat a(\cdot,\cdot)$ is transient. Then 
			  \begin{equation}
				  (h_N(\U_{T|G_N|}^N))_{T \ge 0} \Rightarrow (\bar \U_T)_{T \ge 0},
			  \end{equation}
			  where $(\bar \U_T)_{T \ge 0}$ is the (non spatial) tree-valued Fleming-Viot process with resampling rate $D$ that starts in $\bar \U_0 = [\{1\},0,\delta_1]$. 
			\end{theorem}

			We close this section with a remark on a generalization of this result. This remark is based on the observations in the proof section and we note that even though it seems to be straight forward, we did not prove this yet. 

			\begin{remark}
			One can generalize the result to arbitrary abelian groups $G_N$. The only thing needed is on the one hand that $(\theta_N(\MVX^{N,T}_{h \beta_N}))_{h \ge 0}$ (for a suitable scaling $\beta_N$) converges as a processes with values in $\mathcal M_1([0,1])$ equipped with the weak atomic topology to the (non-spatial) Fleming-Viot process with 
			diffusion rate $D$, and on the other hand that the block (or block counting) process of a spatial Kingman coalescent converges in the same time scale to the non-spatial Kingman coalescent with coalescing rate $D$. 
			\end{remark}

\section{Proofs of the main results}\label{sec.p.res}

	Here we give the proofs of our results. 

	\subsection{Preparations}\label{sec.prep}
		
		In the first section we define the notion of family size decomposition, a concept, introduced in \cite{decomposition}. In the second section, we give a Lemma that is used in several proofs. 
		
		\subsubsection{Family size decomposition}\label{sec.fam.size}
		
			\begin{definition}(Definition of $\f$) \label{defin_basic}
			  Let $\mfu \in \UM$ and recall $(\mathcal S^\downarrow,\dSeq(x,y))$ from Proposition \ref{prop.rel.comp.SMR}. .  We define the map $\f(\mfu,\cdot): (0,\infty) \rightarrow \mathcal S^\downarrow$,
			  \begin{equation}
			    \f(\mfu,h) = (a_1(h),a_2(h),\ldots),
			  \end{equation}
			  where  $a_k(h) \ge a_{k+1}(h)$ is the non-increasing reordering of the sequence $(\mu(\bar B(\rep_i^h,h)))_{i = 1,\ldots,\num(h)}$ 
			  (see Lemma \ref{lem_representatives}). Moreover, we define 
			  \begin{equation}
			    \F:\UM \rightarrow D([0,\infty),\mathcal S^\downarrow), \quad \mfu \mapsto \f(\mfu,\cdot).
			  \end{equation}
			  and call $\F(\mfu)$ the {\it family size decomposition} of $\mfu$.
			\end{definition}
			
			\begin{remark}
			  The function $\F$ was studied in \cite{decomposition} in great detail and we refer all interested readers to this paper. 
			\end{remark}
			
			\begin{lemma}\label{lem.res.fam}
			 Assume that $\mfu,\mfu_1,\mfu_2 ,\ldots \in \UM$ with $\mfu_n \rightarrow \mfu$ in the Gromov weak topology, then $\dSeq(\f(\mfu_n,h),\f(\mfu,h)) \rightarrow 0$ for all continuity points $h > 0$
			 of $\f(\mfu,\cdot)$. Moreover, $\F(\mfu) = \F(\mfu')$ implies $\mfu = \mfu'$, whenever $\mfu,\mfu' \in \UMI$, where 
			$\mfu \in \UMI$ if and only if 
			\begin{equation}\label{eq.uniq.fam.size}
			\sum_{i \in I_1} \f(\mfu,\delta)_i \neq \sum_{i \in I_2} \f(\mfu,\delta)_i
			\end{equation}
			for all $\delta > 0$ and all $I_1,I_2 \in \{i \in \mathbb N:\ \f(\mfu,\delta)_i > 0\}$ with $I_1 \neq I_2$. We call the space $\UMI$ {\it identifiable (by family size decomposition)} (compare also \cite{decomposition}).   
			\end{lemma}
			
			\begin{proof}
			 This is Proposition 8.1 and Theorem 3.11 in \cite{decomposition}.
			\end{proof}

			The connection to our situation is given in the following lemma, 
			
			\begin{lemma}\label{l.ord.f}
			  Let $(\mfu_T)_{T \ge 0} = ([[0,1],r_T,\mu_T])_{t \ge 0} \in D([0,\infty),\UMK)$ a (path of an) evolving 
			  genealogy, with measure representation $\{\mfx^T: T \ge 0\}$. Then $\ord(\tilde \mfx^T_h) = \f(\tilde \mfu_{T+h},h)$ for all $0 \le T,h$, where $\ord(\mu)$ is the non increasing reordering of atoms of 
			  a measure $\mu \in\mathcal M_f([0,1])$. 
			\end{lemma}
			
			\begin{proof}
			  This follows directly by definition of the measure representation and the function $\f$. 
			\end{proof}
			
			\begin{lemma}\label{l.cut.f}
			  Recall Definition \ref{def_Phi} and take $\mfu \in \UM$. Then, one has $\f(\mfu,h) = \f(\Cut^{h'}(\mfu),h-h')$ for all $0 \le h' \le h$. Moreover, observe that 
			  $\nu^{2,\mfu}([0,T]) = \sum_{i = 1}^\infty \f(\mfu,T)^2_i$. 
			\end{lemma}
			
			\begin{proof}
			  This follows directly by definition and Lemma \ref{lem_representatives}.
			\end{proof}

			Next we need the following tightness result: 
			
			\begin{lemma}\label{l.tightness.F}
			  Let $\U^n$, $n = 1,2,\ldots$ be a sequence of $\UM$-valued random variables and let $\UM$ be equipped with the Gromov-weak topology. 
			  Assume that for all $\delta > 0$ and all $\eps > 0$ 
			  \begin{itemize}
			    \item[(i)] there is a compact set $\Gamma \subset \mathcal S^\downarrow$ such that 
			      \begin{equation}
				\limsup_{n \rightarrow \infty} P\left(\f(\U_n,\delta) \in \Gamma^c\right)  \le \eps,
			      \end{equation}
			    \item[(ii)] there is an $H \ge  0$ such that 
			      \begin{equation}
				\limsup_{n \rightarrow \infty} P\left(\left(\sum_{i = 1}^\infty \f(\U^n,H)_i\right)^2-\sum_{i = 1}^\infty \f(\U^n,H)^2_i \ge \eps\right) \le \eps
			      \end{equation}
			  \end{itemize}
			  and that the total mass $(\nu^{1,\U^n})_{n \in \mathbb N}$ is tight. Then, $(\U^n)_{n \in \mathbb N}$ is tight. 
			\end{lemma}
			
			\begin{proof}
			  This is Proposition 3.13 in \cite{decomposition}.
			\end{proof}
			
			As a direct corollary we get
			
			\begin{lemma}\label{l.rel.comp.f}
			  A set $\Gamma \subset \UM$ is relatively compact in the Gromov-weak topology, if 
			  \begin{itemize}
			    \item[(i)]  $\sup_{\mfu \in \Gamma} \overline{\mfu} < \infty$, 
			    \item[(ii)] for all $\eps > 0$, there is an $H \ge  0$ such that 
			      \begin{equation}
				\sup_{\mfu \in \Gamma} \left(\left(\sum_{i = 1}^\infty \f(\mfu,H)_i\right)^2-\sum_{i = 1}^\infty \f(\mfu,H)^2_i\right) \le \eps,
			      \end{equation}
			    \item[(iii)]  for all $\bar \delta > 0$, there is a compact set $K \subset \mathcal S^\downarrow$ such that $\f(\mfu,\bar \delta) \in K $ for all $\mfu \in \Gamma$.
			  \end{itemize}
			\end{lemma}

		\subsubsection{Properties of the \texorpdfstring{$h$}{h}-trunk and the Gromov-Prohorov metric}\label{sec.GPM}
		
			Before we start, we cite a result in \cite{DGP11} (see Section 3.2), namely that the marked Gromov-weak topology can be metricized by the so called {\it marked Gromov-Prohorov distance} $\dmG$.
			Recall  the definition of the Prohorov distance of two finite measures $\mu_1$ and $\mu_2$ on a metric space $(E,r)$ with Borel $\sigma$-field $\B(E)$ 
			\begin{equation}\label{eq_def_Pr}
			\begin{split}
				\dPr(\mu_1,\mu_2):= \inf \Big\{\eps>0:\ &\mu_1(A) \le \mu_2(A^\e) + \eps,\\
				& \mu_2(A) \le \mu_1(A^\e) + \eps \ \textrm{for all } A \textrm{ closed} \Big\},
			\end{split}
			\end{equation}
			where 
			\begin{equation}
				A^\e:=\Big\{x \in E:\ r(x,x') < \eps, \textrm{ for some } x' \in A\Big\}. 
			\end{equation}
			For two marked metric measure spaces $[X,r_X,\mu_X]$ and $[Y,r_Y,\mu_Y]$, the marked Gromov-Prohorov distance is defined as 
			\begin{equation}
				\dmG([X,r_X,\mu_X],[Y,r_Y,\mu_Y]):= \inf_{(\tilde \varphi_X,\tilde \varphi_Y,Z)} d^{(Z,r_Z)}_{\textrm{Pr}} 
				\big(\mu_X\circ \varphi_X^{-1},\mu_Y\circ \varphi_Y^{-1} \big),
			\end{equation}
			where $\varphi_X(x,\type):= (\tilde \varphi_X(x),\type)$ and $\varphi_Y(y,\type):= (\tilde \varphi_Y(y),\type)$ and the infimum is
			taken over all isometric embeddings $\tilde \varphi_X$ and $\tilde \varphi_Y$ from $\supp(\tilde \mu_X)$ and $\supp(\tilde \mu_Y)$
			into some complete separable metric space $(Z,r_Z)$ and $d^{(Z,r_Z)}_{\textrm{Pr}}$
			denotes the Prohorov distance on $\mathcal M_f(Z \times \TYPE)$.

			Recall the functions given in Definition \ref{def_Phi}, then the following holds: 
			\begin{lemma}\label{lem_restr}
			  Let $0<\h$ and $\mfu =  [X,r,\mu]\in \UMK$.
			  \begin{itemize}
			    \item[(i)] If $A \subset X$ is measurable, and $\mu_A(\cdot\times \cdot):= \mu(\cdot \cap A\times \cdot)$ then 
			    \begin{equation}
			     \dmG ([A,r,\mu_A], [X,r,\mu]) \le \mu(X\backslash A\times \TYPE).
			    \end{equation}
			    \item[(ii)] If $\mfu' = [X,r,\mu'] \in \UMK$, then 
			    \begin{equation}
			     \dmG(\mfu,\mfu') \le \dPr(\mu,\mu'),
			    \end{equation}
			    where the Prohorov distance is taken on the set of Borel-measures on $X \times \TYPE$ (with the product metric). 
			    \item[(iii)]  Let $\Cut_h$ and $\CutTwo_h$ be the functions given in Definition \ref{def_Phi}. Then
			      \begin{equation}
			      \dmG(\mfu, \CutTwo_h(\mfu)) \le h,\qquad \dmG(\Cut_h(\mfu),\CutTwo_h(\mfu)) \le h.
			      \end{equation}
			    \item[(iv)] The functions $h \mapsto \Cut_h(\mfu)$ and $h \mapsto \CutTwo_h(\mfu)$ are cadlag. 
			    \item[(v)]  Assume that $\mfu_n$ is a sequence in $\UMK$ with $\mfu_n \rightarrow \mfu \in \UMK$. 
			    Then $\Cut^h(\mfu_n) \rightarrow \Cut^h(\mfu)$ for all continuity points $h \mapsto \Cut^h(\tilde \mfu)$.
			  \end{itemize}
			\end{lemma}

			\begin{proof}
			  (i) Since the identity $id:X \rightarrow X$ is an isometric embedding from $A$ to $X$, it is enough to bound the Prohorov distance
			  of $\mu$ and $\mu_A$. Choosing $\epsilon = \mu(X\backslash A\times \TYPE)$ in  \eqref{eq_def_Pr} gives the result(note that $\mu_A \le \mu$). \\
		
			  (ii) As in (i) one can use the identity as an isometric embedding. Now this is obvious by definition. \\

			  (iii) We use the notation of Definition \ref{def_Phi} and note that $id$ is an isometric embedding from $\{ \rep_i^\h, \ i\in \mathbb N\}$ to $X$. 
			  Let 
			  \begin{equation}
			  \bar \mu ((A_1 \times T_1) \times (A_2\times T_2)):= \sum_{i \in \mathbb N} \int_{A_1\cap \bar B(\rep_i^\h,\h)} \kappa(x,T_1\times T_2) \tilde \mu(dx) \delta_{\rep_i^\h}(A_2)
			  \end{equation}
			  for all measurable sets $A_1,A_2 \subset X$ and $T_1, T_2 \subset \TYPE$ and observe that $\bar \mu$ is a coupling of $\mu$ and $\mu_h$. 
			  Since $\mu_h(\{ \rep_i^\h , \ i\in \mathbb N\}\times \TYPE) = \mu(X\times \TYPE)$, we can use the coupling characterization of the Prohorov distance (see Theorem 3.1.2 in \cite{EK86})
			  and get
			  \begin{equation}
			  \begin{split}
			    \dGP(&\mfu,\CutTwo_\h(\mfu))\\ 
			    &\le \inf\{ \epsilon > 0:\ \bar \mu (\{((x,\type),(x',\type')) \in (X \times \TYPE)^2: r(x,x')+d^\TYPE(\type',\type)\ge  \eps\}) \le \eps\} \\
			    &=  \inf\{ \epsilon > 0:\ \bar \mu (\{((x,\type),(x',\type')) \in (X \times \TYPE)^2: r(x,x') \ge  \eps\}) \le \eps\}
			  \end{split}
			  \end{equation}
			  and if we choose $\eps>\h$ then 
			  \begin{equation}
			  \begin{split}
			    \bar \mu (\{((x,\type),&(x',\type')) \in (X \times \TYPE)^2: r(x,x') \ge  \eps\}) \\
			    &\le \sum_{i,j \in \mathbb N, \ i\not=j} \mu(\bar B(\rep_i^\h,\h)\cap \bar B(\rep_i^\h,\h)\times \TYPE) \delta_{\rep_i^\h}(\bar B(\rep_j^\h,\h))= 0.
			  \end{split}
			  \end{equation}
			  For the second part, we use the same argument as in Section 3 in \cite{Loehr13}: Let $Y:=\{\rep_i^h:\ i \in \{1,\ldots,\num(h)\}$, $r^1 = r$, $r^2 = r - h 1(x\neq y)$ and $\mu^1 = \mu^2 = \mu_h$.
			  We denote by $Y \sqcup Y$ the disjoint union of $Y$ and $Y$ and let $\tilde \varphi_i: Y\to Y \sqcup Y$ the canonical embeddings, $i = 1,2$. Moreover, let $\varphi_i(y,\type) := (\tilde \varphi_i(y),\type)$ for $(y,\type) \in Y \times \TYPE$ and define the metric $\tilde d$ on $Y \sqcup Y$ by 
			  \begin{align}
			  \tilde d(\tilde \varphi_1(x),\tilde \varphi_1(y)) &= r^1(x,y) , \\
			  \tilde d(\tilde \varphi_2(x),\tilde \varphi_2(y)) &= r^2(x,y) ,  \\
			  \tilde d(\tilde \varphi_1(x),\tilde \varphi_2(y)) &= \inf_{z \in Y} (r^1(x,z)+r^2(y,z))+ h, 
			  \end{align}
			  where $x,y \in Y$. Then, as in \cite{Loehr13} it is easy to see that this is a metric on $Y\sqcup Y$ that extends the metrics $r^1 $ and $r^2$ (i.e. $\tilde \varphi_i$ is an isometry for $i = 1,2$) and we have 
			  \begin{equation}
			  \varphi_2(\varphi_1^{-1}(F)) \subset F^{h_0}:= \{(x,\type) \in Y \sqcup Y \times \TYPE:\ \exists (x',\type') \in F \text{ s.t. } 
			  d((x,\type),(x',\type')) <h_0\}, 
			  \end{equation}
			  for all $h_0 > h$, where $d = \tilde d + \dTYPE$. Since $\mu^1 = \mu^2$ this gives: 
			  \begin{equation}
			  \mu^1 \circ \varphi_1^{-1}(F) = \mu^2 \circ \varphi_1^{-1}(F) \le \mu^2 \circ \varphi_2^{-1} (\varphi_2(\varphi_1^{-1}(F)) 
			  \le \mu^2 \circ \varphi_2^{-1} (F^{h_0}) + h_0, 
			  \end{equation}
			  for all $h_0 > h$ and the result follows. \\
			  
			  (iv) A similar argument as in (iii) shows that $\CutTwo_{h'}(\mfu) \rightarrow\CutTwo_{h}(\mfu)$ for $h' \downarrow h$.  By definition, we have $\Cut_{h + \delta} (\mfu) = \Cut_\delta(\Cut_h(\mfu))$ and hence by (iii) $\Cut_{h'}(\mfu) \rightarrow\Cut_{h}(\mfu)$ for $h' \downarrow h$. This shows the right continuity. For the left continuity set  
			  \begin{align}
			  \mfu_h &:= \left[\{\rep_i^h:\ i \in \{1,\ldots,\num(h)\}\},\ r- h \cdot 1(\tau_i \neq \tau_j), \ \mu_h^\circ\right], \\ 
			  \hat \mfu_h &:= \left[\{\rep_i^h:\ i \in \{1,\ldots,\num(h)\}\},\ r, \ \mu_h^\circ\right], 
			  \end{align}
			  where 
			  \begin{equation}
			  \mu_h^\circ(A \times B):= \sum_{i = \in \{1,\ldots,\num(h)\}}\int_{B(\rep_i^h,h)} \KERNEL(x,B) \tilde \mu(dx) \delta_{\rep_i^h}(A)
			  \end{equation} 
			  is given in terms of open balls with radius $< h$ (instead of $\le h$ - see  Remark \ref{rem_pos} (ii)). 
			  By a similar argument as in (iii) we get 
			  \begin{equation}
			  \dmG(\Cut_{h'}(\mfu),\mfu_h) \vee \dmG(\CutTwo_{h'}(\mfu),\hat \mfu_h) \rightarrow 0,\qquad h' \uparrow h.
			  \end{equation}
			  
			  (v)  It is not hard to see that $\nu^{k,\Cut^h(\tilde \mfu)}([0,h']) = \nu^{k,\tilde \mfu}([0,h+h'])$ and that a continuity point of $h \mapsto \Cut^h(\tilde \mfu)$ is a continuity point 
			  of $h \mapsto \nu^{k,\tilde \mfu}([0,h])$ for all $k \in \mathbb N$ and all $h,h' > 0$.  The result follows by the definition of convergence in the marked Gromov-weak topology 
			  (see Definition \ref{def_GS}) and the continuous mapping theorem (see for example Theorem 8.4.1 in \cite{Bog}). 
			\end{proof}
				
		\subsection{Proof of Proposition \ref{prop.MR.EG}}
		
			By the definition of an evolving genealogy, $([0,1],r_T,\mu^h_{T+h}) \in \Cut^h(\mfu_{T+h})$ for all $0 \le T,h$. Moreover, 
			since there is an isometry $\tilde \tau_{h',h}^T:\supp(\tilde \mu^h_{T+h}) \to \supp(\tilde \mu^{h'}_{T+h'})$ for all $0 \le h' \le h$ and the composition 
			$\tilde \tau_{h'',h'}^T\circ \tilde \tau_{h',h}^T$ is an isometry $\supp(\tilde \mu^h_{T+h}) \to \supp(\tilde \mu^{h''}_{T+h''})$ for all $0 \le h'' \le h' \le h$, 
			it is possible to find representatives $([0,1],r_T,\mfx^T_{h}) \in \Cut^h(\mfu_{T+h})$ , such that the identity is an isometric embedding 
			$(\supp(\mfx^T_{h}),r_T) \to (\supp(\mfx^T_{h'}),r_T) $ for all $0 \le h' \le h$. 
			By Lemma \ref{lem_restr} and Lemma 5.8 in \cite{GPW09}, the map $h \mapsto \mfx^T_{h}$ is cadlag and, since $\mfx^T_{h}$ is purely atomic for all $h > 0$ (recall Remark \ref{rem.mr.atomic}), the first direction follows. \par
			The other direction is straight forward.

		\subsection{Proof of Theorem \ref{thm.compact}}

			${}$\\
			
			{\it (1) + (2) part 1: (relative compactness in the Skorohod space)} According to Theorem 3.6.3 in \cite{EK86}, we need to prove two things, (recall the definition of the marked Gromov-Prohorov distance - Section \ref{sec.GPM})
			\begin {itemize}
			\item[(a)] For all $T \ge 0$ there is a compact set $K_T \subset \UMK$ such that 
			  $\mfu_T \in K_T$ for all $\mfu \in \Gamma$.
			\item[(b)] For each $\mathcal T > 0$, 
			  \begin{equation}
			    \lim_{\delta \downarrow 0} \sup_{\mfu \in \Gamma} w'(\mfu,\delta,\mathcal T) = 0,
			  \end{equation}
			  with
			  \begin{equation}\label{p.thm.compact.eq.1}
			    w'(\mfu,\delta,\mathcal T) = \inf_{t_i}\max_i \sup_{s,t \in [t_{i-1}, t_i)} \dmG(\mfu_t,\mfu_s),
			  \end{equation} 
			  where $\{t_i\}$ is a finite partition of $[0,\mathcal T]$ with $\min_i(t_i-t_{i-1}) \ge \delta$ (see (3.6.2) in \cite{EK86}).
			\end {itemize}
			In order to prove (a) we first note that, since $\TYPE$ is compact and according to Theorem 3 in \cite{DGP11}, we only need to verify
			\begin {itemize}
			  \item[(a')] For all $T \ge 0$ there is a compact set $K_T \subset \UM$ such that 
			    $\tilde \mfu_T \in K_T$ for all $\mfu \in \Gamma$.
			\end {itemize}
			In view of this observation, we assume for the following w.l.o.g. $|\TYPE| = 1$, i.e. we assume $\mfu_T \in \UM$ for $T \ge 0$. 
			By the assumptions, Lemma \ref{l.ord.f} (connection of $\f$ and the reordered sizes of atoms of the measure representation) and Lemma 2.5 (c) in \cite{EKatomic} (characterization of the weak atomic convergence in terms of $\ell^1$ convergence of the reordered sizes), we first note, that condition (iii) of Lemma \ref{l.rel.comp.f} (compactness in $\UM$ via conditions on $\f$) is satisfied: 
			\begin{itemize}
			  \item[] for all $0 < \delta < T$, there is a compact set $K \subset \mathcal S^\downarrow$ such that $\f(\mfu_T,\delta) \in K $ for all $\mfu \in \Gamma$.
			\end{itemize}
			Since the initial conditions are relatively compact together with the fact that the map to the total mass is continuous 
			(see Remark \ref{r.totalMassCont}) combined with Remark \ref{rem.mr.atomic} gives that the total mass is uniformly bounded, 
			i.e. (i) of Lemma \ref{l.rel.comp.f} is also satisfied. It remains to verify (ii) of that Lemma, i.e. we need to prove 
			that for all $T \ge 0$ and all $\eps > 0$, there is an $H \ge  0$ such that 
			\begin{equation}\label{thm.comp.proof.1}
			  \sup_{\mfu \in \Gamma} \left(\left(\sum_{i = 1}^\infty \f(\mfu_T,H)_i\right)^2-\sum_{i = 1}^\infty \f(\mfu_T,H)^2_i\right) \le \eps.
			\end{equation}
			We will prove this under condition (i),(ii),(iii'). The result under condition (i),(i'),(ii), (iii), then follow by a simple modification. \par 
			We may assume w.l.o.g. that $\overline{\mfu}_T = 1 $ for all $T \ge 0$ and all $\mfu \in \Gamma$, i.e.
			\begin{equation}
			    \left(\sum_{i = 1}^\infty \f(\mfu_{T},H)_i\right)^2-\sum_{i = 1}^\infty \f(\mfu_T,H)^2_i 
			    = 1-\sum_{i = 1}^\infty \f(\mfu_T,H)^2_i. 
			\end{equation}
			We apply Lemma \ref{l.cut.f} and get for 
			$\mfu_0 = [X,r,\mfx^{\mfu,0}_0]$, by the definition of a measure representation and evolving genealogy, 
			\begin{equation}
			  \begin{split}
			    1-&\sum_{i = 1}^\infty \f(\mfu_T,H+T)^2_i 
			    = 1-\sum_{i = 1}^\infty \f(\Cut^{T}(\mfu_T),H)^2_i \\
			    &= 1-\sum_{i = 1}^\infty \tilde \mfx^0_T\big(\bar B(\rep_i^{H},H)\big)^2
			    = \sum_{i = 1}^\infty \tilde \mfx^0_T\big(\bar B(\rep_i^{H},H)\big)\tilde \mfx^0_T\big(\bar B(\rep_i^{H},H)^c\big)
			  \end{split}
			\end{equation}
			where $\rep_i^{H}$ are given in Lemma \ref{lem_representatives}  and $\bar B(x,h)$ is the closed ball of radius $h> 0$ around $x$ with respect to $r$. Now, by
			Proposition \ref{prop.rel.comp.SMR}, 
			\begin{equation}
			  \begin{split}
			    1-&\sum_{i = 1}^\infty \f(\mfu_T,H+T)^2_i 
			    = \sum_{i = 1}^\infty \tilde \mfx^0_T\big(\bar B(\rep_i^{H},H)\big)\tilde \mfx^0_T\big(\bar B(\rep_i^{H},H)^c\big) \\
			    &\le (C')^2 \sum_{i = 1}^\infty \tilde \mfx^0_\delta\big(\bar B(\rep_i^{H},H)\big)\tilde \mfx^0_\delta\big(\bar B(\rep_i^{H},H)^c\big) 
			    = (C')^2 \nu^{2,\mfu_\delta}\left((H,\infty)\right).
			  \end{split}
			\end{equation}
			Since this holds for all $\delta > 0$ and since the map $\mfu \mapsto \nu^{2,\mfu}$ is continuous and $T \mapsto \mfu_T$ is cadlag, 
			we get for all continuity points $H$ of $H \mapsto \nu^{2,\mfu_0}((H,\infty))$,
			\begin{equation}
			    1-\sum_{i = 1}^\infty \f(\mfu_T,H+T)^2_i 
			    \le  (C')^2 \nu^{2,\mfu_0}\left((H,\infty)\right).
			\end{equation}
			By the relative compactness of the initial condition, (a') follows. \\
			
			In order to prove (b), we take $\eps, \delta > 0$ and assume w.l.o.g. that $\frac{\e}{8} K = \mathcal T$ and $L \delta = \mathcal T$ for some $K,L \in \mathbb N$. Since 
			$\{(\mfx^{\mfu,T}_h)_{h \ge 0}: \ \mfu \in \Gamma\}$ is relatively compact in the Skorohod topology for all $T \ge 0$, we apply Theorem 3.6.3 in \cite{EK86} and choose $\delta$ small enough, 
			such that 
			\begin{equation}
			  \sup_{\mfu \in \Gamma} w'(\mfx^{\mfu,s_j},2\delta,\mathcal T) \le \frac{\eps}{2(K+1)},\qquad \text{for all } j = 1,\ldots,K, 
			\end{equation}
			where $w'(\mfx^{\mfu,T},\delta,\mathcal T)$ is defined as in \eqref{p.thm.compact.eq.1} but with $\dmG$ replaced by $\dPr$. \par 
			Now, let $s_j = (j-1) \frac{\e}{8}$, $j = 1,\ldots,K+1$ and $t_i = (i-1) \delta$, $i = 1,\ldots,L+1$. We define $I(i) := j$ iff $s_{j} \le t_i < s_{j+1}$. 
			Note that for $s,t \in [t_{i-1},t_i)$ ($s_{-1}:=0$) 
			\begin{equation}
			  \begin{split}
			    \dmG(\mfu_t,\mfu_s) \le  \dmG(\mfu_s,&\Cut^{s-s_{I(i-1)-1}}(\mfu_s)) \\
			    &+ \dmG(\Cut^{s-s_{I(i-1)-1}}(\mfu_s),\Cut^{t-s_{I(i-1)-1}}(\mfu_t)) \\
			    &+ \dmG(\mfu_t,\Cut_{t-s_{I(i-1)-1}}(\mfu_t)). 
			  \end{split}
			\end{equation}
			Now, by Lemma \ref{lem_restr}, 
			\begin{equation}
			  \dmG(\mfu_s,\Cut^{s-s_{I(i-1)-1}}(\mfu_s)) + \dmG(\mfu_{t},\Cut^{t-s_{I(i-1)-1}}(\mfu_{t})) \le  2\frac{\e}{4} = \frac{\e}{2}
			\end{equation}
			and 
			\begin{equation}
			  \begin{split}
			    \max_i \sup_{s,t \in [t_{i-1}, t_i)} &\dmG(\Cut^{s-s_{I(i-1)-1}}(\mfu_s),\Cut^{t-s_{I(i-1)-1}}(\mfu_t)) +  \frac{\e}{2}\\
			    &\le \max_i \sup_{s,t \in [t_{i-1}, t_i)} \dPr(\mfx^{\mfu,s_{I(i-1)-1}}_{s-s_{I(i-1)-1}},\mfx^{\mfu,s_{I(i-1)-1}}_{t-s_{I(i-1)-1}}) +  \frac{\e}{2} \\
			    &\le \max_{j \in \{1,\ldots,K+1\}} \max_{i}\sup_{s,t \in [t_{i-1}, t_i)} \dPr(\mfx^{\mfu,s_{j}}_{s},\mfx^{\mfu,s_{j}}_{t})  +  \frac{\e}{2}\\
			    &\le \sum_{j = 1}^{K+1} \max_{i}\sup_{s,t \in [t_{i-1}, t_i)} \dPr(\mfx^{\mfu,s_{j}}_{s},\mfx^{\mfu,s_{j}}_{t}) + \frac{\e}{2},
			  \end{split}
			\end{equation}
			where we again used Lemma \ref{lem_restr}. Since $t_{i}-t_{i-1} =\delta$, we  have
			\begin{equation}
			    \max_{i} \sup_{s,t \in [t_{i-1}, t_i)} \dPr(\mfx^{\mfu,s_{j}}_{s},\mfx^{\mfu,s_{j}}_{t})  
			    \le \max_{i} \sup_{s,t \in [\tilde t_{i-1}, \tilde t_i)} \dPr(\mfx^{\mfu,s_{j}}_{s},\mfx^{\mfu,s_{j}}_{t}) 
			\end{equation}
			for all partitions $\{\tilde t_i\}$ of $[0,\mathcal T]$ with $\min_{i}|\tilde t_{i}-\tilde t_{i-1}| \ge 2\delta$. Finally, 
			\begin{equation}
			  \begin{split}
			    \sup_{\mfu \in \Gamma}\ & w'(\mfu,\delta,\mathcal T) 
			    \le \sup_{\mfu \in \Gamma} \sum_{j = 0}^{K}  w'(\mfx^{\mfu,s_j},2\delta,\mathcal T) + \frac{\e}{2}\\
			    &\le \sum_{j = 1}^{K+1} \sup_{\mfu \in \Gamma} w'(\mfx^{\mfu,s_j},2\delta,\mathcal T) + \frac{\e}{2}
			    \le (K+1) \frac{\e}{2(K+1)} + \frac{\e}{2} = \eps.
			  \end{split}
			\end{equation}
			  
			{\it (2) - part 2}  In order to prove relative compactness in $\SEG$, we pick a sequence $\mfu^n \in \Gamma$ and assume that $\mfu^n \rightarrow \mfu$ in $D([0,\infty),\UMK)$ 
			for some $\mfu \in \Gamma$ (this is true along some subsequence where we suppress the dependence).
			Fix a continuity point $T \ge 0$ of $T \mapsto \mfu_T$ and let $\mfx^{n,T}$ be a measure representation of $\mfu^n$ then, along some subsequence, $\mfx^{n_k,T} \rightarrow \mfx^{T} \in \SMR$.
			By Proposition 3.5.2 in \cite{EK86} (pointwise convergence for all continuity points of the limit function), $[[0,1],r_T^{n_k},\mfx^{n_k,T}_0] = \mfu^n_T \rightarrow \mfu_T = [[0,1],r_T,\mu] $ and, by Lemma 2.4 in \cite{DGP11}, there is a 
			complete separable metric space $(Z,r_Z) $ and isometric embeddings 
			$\tilde \varphi,\tilde \varphi_1,\tilde \varphi_2,\ldots: [0,1] \to Z$, with $\mfx^{n_k,T}_0 \circ \varphi_{n_k}^{-1} \Rightarrow \mu \circ \varphi^{-1}$ (independent of the subsequence). 
			By the continuous mapping theorem (see for example Theorem 8.4.1 in \cite{Bog}) this is enough to prove that for all continuity points $h \ge 0$ of $h \mapsto \mfx^{T}_h$,  
			\begin{equation}
			 \Cut^h(\mfu^{n_k}_{T+h}) = [[0,1],r_T^{n_k},\mfx^{n_k,T}_h] \rightarrow [[0,1],r_T,\mfx^{T}_h].
			\end{equation}
			In particular $\mfx^T_0 = \mu$. It remains to verify $\Cut^h(\mfu^{n}_{T+h}) \rightarrow \Cut^h(\mfu_{T+h})$ for all but countable many $h \ge 0$. 
			Since $h \mapsto \Cut^h(\mfu)$ is cadlag (see Lemma \ref{lem_restr} (iv)) and the set of discontinuity points of a cadlag function is at most countable (see Lemma 3.5.1. in \cite{EK86}), 
			this is Lemma \ref{lem_restr} (v). It follows that 
			\begin{equation}
			 \Cut^h(\mfu^{n_k}_{T+h}) = [[0,1],r_T^{n_k},\mfx^{{n_k},T}_h] \rightarrow [[0,1],r_T,\mfx^{T}_h] = \Cut^h(\mfu_{T+h}),
			\end{equation}
			for all but countably many $h \ge 0$. By the cadlag property of $h \mapsto \mfx^T_h$ in combination with Lemma \ref{lem_restr} (ii), this extends to all $h \ge 0$ and since $T\mapsto \mfu_T$
			is cadlag also to all $T \ge 0$. Hence $\mfu$ is an evolving genealogy, by Proposition \ref{prop.MR.EG}. The last part of (2) is a reformulation of this observation. \\
			  
			{\it (3)} We take a subsequence and assume that $\mfu^{n_k} \rightarrow \mfu \in D([0,\infty),\UMK)$. Then, by Proposition 3.5.2 in \cite{EK86}, 
			$\mfu^{n_k}_T \rightarrow \mfu_T$ for all continuity points of $T \mapsto\mfu_T$. Note that by Lemma 3.5.1 in \cite{EK86}, the number of discontinuity points is at most countable. \par 
			Now observe, that by the connection of the family size decomposition and the measure representation (see Lemma \ref{l.ord.f}), the property of the measure representation given in \eqref{eq.uniq.mvr} is exactly the property \eqref{eq.uniq.fam.size} in Lemma \ref{lem.res.fam}, i.e. the corresponding ultra-metric measure space is identifiable. We can now use Lemma 2.5 in \cite{EKatomic}  and the fact that $\F(\tilde \mfu)$ is cadlag, to deduce that $\tilde \mfu_T$ is unique for 
			all continuity points $T$, and therefore, since $T\mapsto \tilde \mfu_T$ 
			is cadlag, for all $T \ge 0$. This implies that if $\mfu'$ is another limit point, then $\tilde \mfu_T = [X,r,\tilde \mfx^T_0] = [X,r,(\tilde \mfx^T_0)'] = \tilde \mfu_T'$. But $\mfx^T$
			does not depend on the subsequence and therefore, we also have  $\mfu_T = [X,r,\mfx^T_0] = [X,r,(\mfx^T_0)'] = \mfu_T'$.

			\subsection{Proof of Proposition \ref{prop.rel.comp.SMR}}		
			
			First note that, according to Lemma 2.5 in \cite{EKatomic}, condition (ii) is equivalent to 
			\begin{itemize}
			 \item[(ii')] For all  $h > 0$, $\{\tilde \mfx_h:\ \mfx \in \Gamma\}$ is relatively compact in the space of purely atomic finite measures, 
			    $\mathcal M_f^a([0,1])$ (equipped with the subspace topology). 
			\end{itemize}

			``$\Leftarrow$'' Let $(\mfx^n)_{n \in \mathbb N}$ be a sequence in $\Gamma$. Then there is a subsequence, where we suppress the dependence, and a $\mfx \in D([0,\infty),\mathcal M_f([0,1]\times \TYPE))$ 
			such that $\mfx^n \rightarrow \mfx$. Since $(ii)$ holds, $\mfx_\delta$ is purely atomic for all $\delta > 0$. Assume now, that there is a $x$ such that $\tilde \mfx_t(\{x\})> 0$ but 
			$\tilde \mfx_{t'}(\{x\}) = 0$ for some $0 <t' < t$. Let $t_n, t_n'$ be two sequences in $(0,\infty)$ with $t_n \rightarrow t$ and $t_n' \rightarrow t'$ and
			$\tilde\mfx_{t_n}^n \rightarrow \tilde\mfx_t$ and $\tilde\mfx_{t_n'}^n\rightarrow \tilde\mfx_{t'} $ in the weak atomic topology. Then, according to Lemma 2.5 in \cite{EKatomic}, 
			there is exactly one sequence  $x_n \in \mathcal A(\tilde \mfx_{t_n}^n)$, where $\mathcal A$ denotes the set of atoms, such that $x_n \to x$ and $\tilde \mfx^n_{t_n}(\{x_n\}) \rightarrow \mfx_{t}(\{x_n\})$. Moreover, by the same 
			Lemma in \cite{EKatomic}, $\limsup_{n \rightarrow \infty} \mfx_{t_n'}(\{x_n\}) = 0$, which contradicts (i).  \par 
			``$\Rightarrow$'' Clearly relative compactness in the Skorohod space as well as (ii) holds. To verify (i), one can use a similar argument as for the $\Leftarrow$ direction, we leave that to the reader.\\
			
			The last part is a consequence of (ii'), Lemma \ref{lem.res.fam} ($\mfu$ is - under some conditions - uniquely determined by its family size decomposition) and Lemma \ref{l.ord.f} (connection of the family size decomposition and the measure representation).
			
		\subsection{Proof of Theorem \ref{thm.tightness}}
			
			To show relative compactness recall that we have to prove two things (see Corollary 3.7.4 in \cite{EK86}): 
			\begin{itemize}
			  \item[(1)] For all $\eps> 0$, $T \ge 0$, there is a compact set $K\subset \UMK$ such that $\liminf_{n \rightarrow \infty} P(\U^n_T \in K) \ge 1-\e$. 
			  \item[(2)] For every $\eps > 0$ and $\mathcal T > 0$ there is a $\delta > 0$ such that 
			    \begin{equation}
			      \limsup_{n \rightarrow \infty} P(w'(\U^n,\delta,\mathcal T) \ge \eps) \le \eps,
			    \end{equation} 
			    with
			    \begin{equation}\label{eq.mod.cont}
			      w'(\U^n,\delta,\mathcal T) = \inf_{t_i}\max_i \sup_{s,t \in [t_{i-1}, t_i)} \dmG(\U^n_t,\U^n_s),
			    \end{equation} 
			    where $\{t_i\}$ is a finite partition of $[0,\mathcal T]$ with $\min_i(t_i-t_{i-1}) \ge \delta$ (see (3.6.2) in \cite{EK86}) and the Gromov-Prohorov distance is defined in Section \ref{sec.GPM}).
			\end{itemize}

			(1)  First note that $\f(\tilde \U^n_{T+h},h) \stackrel{d}{=} \ord(\tilde {\mathcal X}^{T,n}_{h})$ (see Lemma \ref{l.ord.f}) and 
			that $\tilde {\mathcal X}^{T,n}_{\delta}$ is tight for all $\delta > 0$, 
			where $\mathcal M_f([0,1])$ is equipped with the weak atomic topology, as $n \rightarrow \infty$.
			By Lemma 2.5 (c) in \cite{EKatomic}, this is enough to prove tightness of $\ord(\tilde {\mathcal X}^{T,n}_{h})$ for all $h > 0$, 
			as random variables with values in $(\mathcal S^\downarrow,d^1)$. 
			We can now apply Theorem 4 in \cite{DGP11}, i.e. we need to prove tightness of $\tilde \U^n_T$ and we use Lemma \ref{l.tightness.F} to do that. 
			To check the conditions of this Lemma works analogue to the proof of Theorem \ref{thm.compact}. \\
			
			(2) Let $\eps, \delta > 0$ and assume w.l.o.g. that $\frac{\e}{8} K = \mathcal T$ and $L \delta = \mathcal T$ 
			for some $K,L \in \mathbb N$. Moreover, since 
			$\{(\MVX^{\mfu,T}_h)_{h \ge 0}: \ \mfu \in \Gamma\}$ is relatively compact in the Skorohod topology for 
			all $T \ge 0$, we apply Theorem 3.7.2 in \cite{EK86} and choose $\delta$ small enough, 
			such that 
			\begin{equation}
			  \sup_{n \in \mathbb N} P\left(w'(\mathcal X^{n,s_j},2\delta,\mathcal T) \ge \frac{\eps}{K+1}\right)\le \frac{\eps}{K+1},\qquad  j = 1,\ldots,K, 
			\end{equation}
			where $w'(\mathcal X^{n,t},\delta,\mathcal T)$ is defined as in \eqref{eq.mod.cont} but with $\dmG$ replaced by $\dPr$. \par 
			Now, analogue to the proof of Theorem \ref{thm.compact}, we get
			\begin{equation}
			  \begin{split}
			    P(&w'(\U^n,\delta,T) \ge \eps) 
			      \le \sum_{j = 0}^{K} P(w'(\MVX^{n,s_j},2\delta,T) \ge \frac{\e}{2}) \\
			      &\le \sum_{j = 0}^{K} \sup_{n}P(w'(\MVX^{n,s_j},2\delta,T) \ge \frac{\e}{K+1}
			      \le (K+1) \frac{\e}{K+1} = \eps.
			  \end{split}
			\end{equation}
		
		\subsection{Proof of Theorem \ref{thm.convergence.II}}
		
			Before we start, we need the following: 
			\begin{lemma}\label{lem.coupl.C}
			Let $\mfu \in \SEG$ with measure representation $\{\mfx^{\mfu,T}:\ T \ge 0\}$. If $\tilde \mfx^{\mfu,T+t}_{h} \le C \tilde \mfx^{\mfu,T+t}_{\delta}$ for some 
			$T,C,t > 0$ and $0 < \delta \le h$, then $\tilde \mfx^{\mfu,T}_{t+h} \le C\tilde \mfx^{\mfu,T}_{t+\delta}$. 
			Moreover, let $A_\delta^T := \{\mfu \in \SEG:\ \sup_{h \ge \delta}\tilde \mfx^{\mfu,T}_{h} \le C\tilde \mfx^{\mfu,T}_{\delta}\}$, then 
			\begin{equation}
			 A_{\delta}^{T+\delta'} \subset A_{\delta+\delta'}^{T},\qquad \text{for all } \delta,\delta' >0.
			\end{equation}
			\end{lemma}
			
			\begin{proof}
			  First note that by definition of $\Cut^h$ (see Definition \ref{def_Phi}) and the definition of a measure representation, 
			  there is on the one hand a measure preserving map $\psi^{T,t}_{\delta}: \supp(\tilde \mfx^{\mfu,T+t}_{\delta}) \to \supp(\tilde \mfx^{\mfu,T}_{t + \delta})$ for all $T,t \ge 0$ and $\delta > 0$.
			  On the other hand, this maps can be chosen in such a way, that $\psi^{T,t}_{\delta}\big|_{\supp(\tilde \mfx^{T+t}_{h})} = \psi^{T,t}_{h}$. It follows that 
			  \begin{equation}
			      \tilde \mfx^{\mfu,T}_{t+h}(A) = \tilde \mfx^{\mfu,T+t}_{h} \circ (\psi^{T,t}_{\delta})^{-1}(A) \le C\tilde \mfx^{\mfu,T}_{t+\delta}(\psi^{T,t}_{\delta})^{-1}(A) = \tilde \mfx^{\mfu,T}_{t+\delta}(A)
			  \end{equation}
			\end{proof}
			
			We are now ready to prove the theorem. Let $\U$ be a weak limit point along some subsequence of $\U^n$, where we suppress the dependence and assume for the moment 
			$\U^n \Rightarrow \U$ (as processes). By Skorohod's representation theorem (see for example Theorem 8.5.4 in \cite{Bog}), we can 
			find a probability space and random variables $\hat \U^n \stackrel{d}{=} \U^n$, $n = 1,2,\ldots$ and $\hat \U \stackrel{d}{=} \U$, 
			such that $\hat \U^n \rightarrow \hat \U$ almost surely (see Proposition 3.5.2 in \cite{EK86}). By convergence in the Skorohod topology, this implies that for continuity points $T$ of $T \mapsto \hat \U_T$ we have
			$\hat \U^n_T \rightarrow \hat \U_T$ almost surely. We can now apply Lemma \ref{lem.res.fam} to get
			$ \f(\hat \U^n_T,h) \rightarrow \f(\hat \U_T,h)$ almost surely, for all but at most countably many $T \ge 0$ and $h > 0$ (see Lemma 3.5.1 in \cite{EK86}). We now use the fact that $\f$ is also given in terms of the reordered sizes of atoms of the measure representation (see Lemma \ref{l.ord.f}) and the fact that that convergence in the weak atomic topology implies $\ell^1$ convergence of the reordered sizes of atoms (see Lemma 2.5 (c) in \cite{EKatomic}). This gives that the f.d.d. law of $\{\f(\tilde{\hat \U}_T,h): T > 0, h > 0\}$ is completely determined by the law of
			$\{(\tilde{\mathcal X}^{T+k\cdot h}_{l\cdot h})_{1\le l\le k \le K}:\ K \in \mathbb N, h \in D_1,T\in D_2\}$. Combining this with the fact that $\f(\mfu,\cdot)$ uniquely determines $\mfu$ 
			as long as $\mfu$ is identifiable, i.e. $\mfu \in \UMIK$ (see Lemma \ref{lem.res.fam} and compare the proof of Theorem \ref{thm.compact}), which 
			is satisfied by assumption (v) and Proposition 3.16 in \cite{decomposition} gives uniqueness of the finite dimensional distributions of $\U$ and hence the convergence, $\U^n \Rightarrow \U$, follows. \\
			
			Next we prove that the limit $\U$ is again an evolving genealogy. To do that, we need to prove that tightness of the processes holds in $\mathcal M_1(\SEG)$. Note that the above construction immediately gives 
			for all $\eps > 0$ the existence of a compact set $\Gamma^0 \subset D([0,\infty),\UMK)$  such that 
			\begin{align}
			 \limsup_{n \rightarrow \infty} P\left(\U^n \in (\Gamma^0)^c\right) &\le \frac{\eps}{3}.
			\end{align}
			Moreover, applying Corollary 3.12 in \cite{decomposition}, where we note that the result stays valid when one considers the Gromov-weak topology (this follows by a careful read of the proof of Theorem 3.11 in \cite{decomposition} - compare also Lemma \ref{lem.res.fam}) together with the fact that the family size decomposition can be reformulated in terms of the reordered vector of sizes of atoms of the measure representation (see Lemma \ref{l.ord.f}), 
			gives the existence of a compact set  $K \subset \mathcal S^\downarrow$  such that the following compact containment condition (see Remark 3.7.3 in \cite{EK86}) holds:
			\begin{align}
			 \limsup_{n \rightarrow \infty} P\left(\U^n \in (\Gamma^1)^c\right) &\le \frac{\eps}{3},
			\end{align}
			where
			\begin{equation}
			 \Gamma^1 := \{\mfu \in \SEG: \ \ord(\tilde\mfx^{\mfu,T}_h) \in K \text{ for all } h > 0, T \ge 0\}.
			\end{equation}
			Let $\tau > 0$, $\bar \delta > 0$, $C  \in \mathbb N$ and $l = 1,\ldots,\lceil \tau/\bar \delta\rceil$. Define
			\begin{equation}
			 \Gamma^2_{\tau,\bar \delta,C,l} :=  \left\{\mfu \in \SEG: \ \sup_{h \in [(l+1)\bar \delta - T,(l+2)\bar \delta - T)} \tilde \mfx^{\mfu,T}_h \le C\tilde \mfx^{\mfu,T}_{(l+1)\bar \delta - T} \text{ for all } 
			 T \in [(l-1)\bar \delta,l\bar \delta)\right\}. 
			\end{equation}
			By  Lemma \ref{lem.coupl.C}, 
			\begin{equation}
			  \sup_{h \in [(l+1)\bar \delta - T,\tau]} \tilde \mfx^{\mfu,T+\bar \delta}_h \le C\tilde \mfx^{\mfu,T+\bar \delta}_{(l+1)\bar \delta - T}.
			\end{equation}
			implies 
			\begin{equation}
			  \sup_{h \in [(l+2)\bar \delta - T,\tau]} \tilde \mfx^{\mfu,T}_h \le C\tilde \mfx^{\mfu,T}_{(l+2)\bar \delta - T}
			\end{equation}
			Let $\bar \delta_k \downarrow 0$ for $k \to \infty$ and define
			\begin{equation}
			\begin{split}
			 \Gamma^2&:=  \bigcap_{\tau \in \mathbb N}\bigcap_{k \in \mathbb N} \bigcup_{C \in \mathbb N} \bigcap_{l = 1}^{\lceil \tau/\bar \delta_k\rceil}\Gamma^2_{\tau,\bar \delta,C,l}
			 \end{split} 
			\end{equation}
			Now, according to Proposition \ref{prop.rel.comp.SMR} and Theorem \ref{thm.compact}, the set $\bar \Gamma = \Gamma^0 \cap \Gamma^1 \cap \Gamma^2$ is relatively compact in $\SEG$ and we have 
			\begin{equation}
			 \limsup_{n \rightarrow \infty} P(\U^n \in \bar \Gamma^c) \le \frac{2\eps}{3} + \limsup_{n \rightarrow \infty} P(\U^n \in (\Gamma^2)^c).  
			\end{equation}
			Next observe that by Lemma \ref{lem.coupl.C} we have for a strong measure representation $\{\hat{\mathcal X}^{\U^n,T}: \ T \ge 0\}$, 
			\begin{equation}
			 \begin{split}
			    P(\U^n &\in \Gamma^2_{\tau,\bar \delta,C,l}) \\
			    &= P\left(\forall T \in [(l-1)\bar \delta,l \bar \delta): 
			    \sup_{h \in [(l+1)\bar \delta - T,(l+2)\bar \delta - T)} \tilde{\hat{\mathcal X}}^{\U^n,T}_h \le C\tilde{\hat{\mathcal X}}^{\U^n,T}_{(l+1)\bar \delta - T}\right)\\
			    &\ge P\left(  \sup_{h \in [2\bar \delta,3\bar \delta)} \tilde{\hat{\mathcal X}}^{\U^n,(l-1) \bar \delta}_h \le C\tilde{\hat{\mathcal X}}^{\U^n,(l-1) \bar \delta}_{2 \bar \delta}\right)\\
			    &= P\left(\sup_{h \in [2 \bar \delta,3\bar \delta)} \tilde {\mathcal X}^{\U^n,(l-1) \bar \delta}_h \le C\tilde{\mathcal X}^{\U^n,(l-1) \bar \delta}_{2\bar \delta}\right).
			 \end{split}
			\end{equation}
			Take  $\bar \delta_k = 2^{-k-1}$, for $k \in \mathbb N$, then
			\begin{equation} 
			\bigcap_{l = 1}^{\lceil \tau/\bar \delta_k\rceil}\Gamma^2_{\bar \delta,l} = \bigcap_{l = 1}^{\lceil 2^k \tau\rceil}\Gamma^2_{k,l} 
			\supset \bigcap_{l = 1}^{\lceil 2^k \tau'\rceil}\Gamma^2_{k+1,l}, 
			\end{equation}
			for all $k \in \mathbb N$ and all $\tau' \ge \tau$. It follows, by assumption, that 
			\begin{equation}
			 \begin{split}
			    \inf_{n \in \mathbb N}P(\U^n &\in \Gamma^2) \\
			    &\ge \inf_{n \in \mathbb N} \inf_{\tau \in \mathbb N} \inf_{k \in \mathbb N} \sup_{C \in \mathbb N}P(\forall l = 0,\ldots,2^k \tau:  
			    \sup_{h \in[2^{-k},2^{-k+1})} \mathcal X^{\U^n,l 2^{-k}}_h \le C{\mathcal X}^{\U^n,l 2^{-k}}_{2^{-k}}) \\
			    &\ge  \inf_{\tau,k \in \mathbb N} \sup_{C \in \mathbb N} \inf_{n \in \mathbb N} P(\forall l = 0,\ldots,2^k \tau:  
			    \sup_{h \in[2^{-k},2^{-k+1})} \mathcal X^{\U^n,l 2^{-k}}_h \le C{\mathcal X}^{\U^n,l 2^{-k}}_{2^{-k}})  \\
			    &\ge 1-\eps.
			 \end{split}
			\end{equation}

	\section{Proof of the application}\label{sec.p.app}
		
		Here we prove the applications. 
		
		\subsection{Proofs for section \ref{sec.TVIMM}}\label{sec.p.TVIMM}
		
			In the following we will show how to construct a measure representation for the tree-valued interacting Moran models.
			In section \ref{sec.MM.King} we recall the definition of the measure-valued Moran model and show how one can construct the genealogy of a Moran model given the Kingman coalescent.  
			In section \ref{sec.MVR.TVIMM} we use these observations to give a measure representation, i.e. we prove Theorem \ref{thm.TVIMMM.EG}. 
			
		\subsubsection{The measure-valued interacting Moran models and the spatial Kingman coalescent}\label{sec.MM.King}
			
			Firstly, we define the measure-valued interacting Moran models, i.e. the model that describes the evolution of relative frequencies of types in a population that 
			lives on some geographical space. Secondly, we define the Kingman coalescent, which can be used to model the genealogy of the population for a fixed time $T > 0$. \\
			
			{\bf The measure-valued model:} We denote by $(\hat \eta_i(t),\xi_i(t))_{i \in I_N}$ the type and location (in $G$) of an individual $i$, where the type space is $[0,1]$. Then this process 
			  (with values in $([0,1]\times G)^{I_N}$) has the following dynamic:  
			  \begin{itemize}
			    \item (resampling) At rate $\gamma/2$ we pick a pair $(i,j) \in I_N\times I_N$, $i \neq j$. If $\xi_i(t) = \xi_j(t)$ we have the following transition: 
			      \begin{equation}
				(\hat \eta_k(t),\xi_k(t)) \rightarrow \left\{\begin{array}{ll}
				  (\hat \eta_i(t),\xi_i(t)),&\quad \text{if } k = j,\\
				  (\hat \eta_k(t),\xi_k(t)),&\quad \text{if } k\neq j. 
				\end{array}\right.
			      \end{equation}
			    \item (migration) At rate $1$ we pick an individual $i \in I_N$ and we have the following transition:  
			      \begin{equation}
				(\hat \eta_k(t),\xi_k(t)) \rightarrow \left\{\begin{array}{ll}
				  (\hat \eta_i(t),g),&\quad \text{if } k = i, \text{ with rate } a(\xi_i(t),g),\\
				  (\hat \eta_k(t),\xi_k(t)),&\quad \text{if } k\neq i. 
				\end{array}\right.
			      \end{equation}
			  \end{itemize}
			  
			  Now, we are interested in the frequency of types at each colony. Namely, we define 
			  \begin{equation}
			    X^{N,g}_t:= \frac{1}{|\{i:\ \xi_i(t) = g\}|\wedge 1} \cdot \sum_{i \in I_N} \delta_{\hat \eta_i(t)} 1(\xi_i(t) = g) 
			  \end{equation}
			  for $g \in G$ and we call the process $((X_t^{N,g})_{g \in G})_{t \ge 0}$ with values in $\mathcal M_{\le 1}([0,1])^G$ (note that $X^{N,g}_t$ is a probability measure except for the case where no individuals are located at some site $g \in G$) a {\it system of interacting measure-valued Moran models}. \\
			  
			  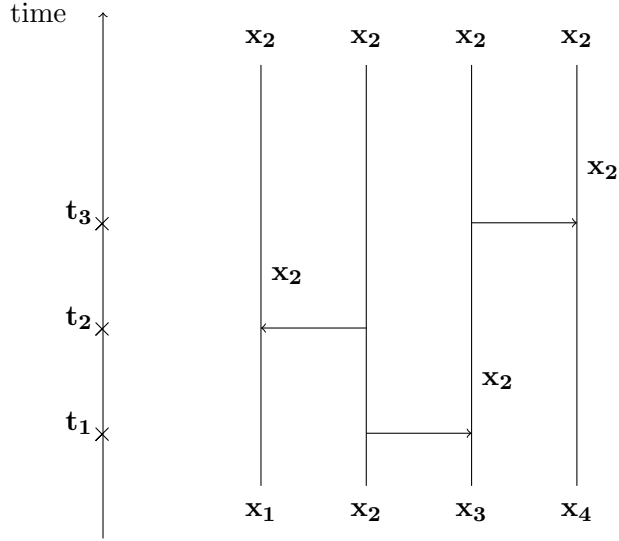
\begin{figure}[ht]
			    \centering
			    \begin{tikzpicture}[scale=0.7]
			      \draw (1,4)-- (1,-4);
			      \draw (3,4)-- (3,-4);
			      \draw (5,4)-- (5,-4);
			      \draw (7,4)-- (7,-4);
			      \draw [->] (-2.,-5.) -- (-2.,5.);
			      \draw (1,-4.5) node {$\mathbf{x_1}$};
			      \draw (3,-4.5) node {$\mathbf{x_2}$};
			      \draw (5,-4.5) node {$\mathbf{x_3}$};
			      \draw (7,-4.5) node{$\mathbf{x_4}$};
			      \draw [->] (5,1) -- (7,1);
			      \draw (-2,-2.8) node[anchor=east] {$\mathbf{t_1}$};
			      \draw (-2,-0.8) node[anchor=east] {$\mathbf{t_2}$};
			      \draw (-2,1.2) node[anchor=east] {$\mathbf{t_3}$};
			      \draw (-2.5,5) node[anchor=east] {$\text{time}$};
			      \draw [->] (3.,-1.) -- (1.,-1.);
			      \draw [->] (3.,-3.) -- (5.,-3.);
			      \draw (5.5,-2) node {$\mathbf{x_2}$};
			      \draw (1.5,0) node {$\mathbf{x_2}$};
			      \draw (7.5,2) node {$\mathbf{x_2}$};
			      \draw (1,4.5) node {$\mathbf{x_2}$};
			      \draw (3,4.5) node {$\mathbf{x_2}$};
			      \draw (5,4.5) node {$\mathbf{x_2}$};
			      \draw (7,4.5) node {$\mathbf{x_2}$};
			      \draw [color=black] (-2,1)-- ++(-4pt,-4pt) -- ++(7pt,7pt) ++(-7pt,0) -- ++(7pt,-7pt);
			      \draw [color=black] (-2,-1)-- ++(-4pt,-4pt) -- ++(7pt,7pt) ++(-7pt,0) -- ++(7pt,-7pt);
			      \draw [color=black] (-2,-3)-- ++(-4pt,-4pt) -- ++(7pt,7pt) ++(-7pt,0) -- ++(7pt,-7pt);
			    \end{tikzpicture}
			    \caption{\label{f.graph.constr.measure} \footnotesize
			    We are in the situation of Figure \ref{f.graph.constr}. Here $(\hat \eta_i(t_0))_{i = 1,\ldots,4} = (x_i)_{i = 1,\ldots,4}$, where 
			    $x_i \in [0,1]$ (not necessarily pairwise different). The process $t \mapsto (\hat \eta_i(t))_{i = 1,\ldots,4}$ is constant up to the three times $t_1,t_2,t_3$, where $\hat \eta(t_1) = (x_1,x_2,x_2,x_4)$, $\hat \eta(t_2) = (x_2,x_2,x_2,x_4)$ and $\hat \eta(t_3) = (x_2,x_2,x_2,x_2)$. The corresponding frequency of types $X^{N}_t$ is also constant up to the three times $t_1,t_2,t_3$, where 
			    $X^{N}_{t_0} = \frac{1}{4} \sum_{i = 1}^4 \delta_{x_i}$,  $X^{N}_{t_1} = \frac{1}{4} \delta_{x_1}+\frac{1}{2} \delta_{x_2}+\frac{1}{4} \delta_{x_4}$, $X^{N}_{t_2} =  \frac{3}{4} \delta_{x_2}+\frac{1}{4} \delta_{x_4}$ and $X^{N}_{t_3} =  \delta_{x_2}$.
			    }
			  \end{figure}	
			  
			  \begin{remark}\label{r.Generator.MM}
			    Alternatively, one can also define 
			    \begin{equation}
			     \hat X^{N,g}_t:= \frac{1}{N} \cdot \sum_{i \in I_N} \delta_{\hat \eta_i(t)} 1(\xi_i(t) = g)
			    \end{equation}
			    and we note that this process can be characterized as the solution of a well-posed martingale problem (see section 4 in \cite{EK86}), where the linear operator 
			    $G^N$, acting on continuous functions $F \in C_b(\mathcal M_{f}([0,1])^G)$, is given by 
			    \begin{equation}
			      \begin{split} 
			      G^N F((y_g)_{g \in G}) =&  \sum_{g \in G} \gamma \cdot \binom{N y_g([0,1])}{2}\int \int \big(F(y_g^{u,v}) - F(y)\big) y_{g}(du)y_{g}(dv)\\
				&+ \sum_{g,q \in G} N \cdot y_g([0,1]) a(g,q) \int \big(F(\hat y_{g,q}^{u}) - F(y)\big) y_{g}(du),
			      \end{split}
			    \end{equation}
			    where
			    \begin{equation}
			      (y_g^{u,v})_\xi = \left\{\begin{array}{cl}
				y_\xi & \textrm{if } \xi \not=g,\\[0.2cm]
				(y_g + \frac{1}{N}\delta_{u} - \frac{1}{N}\delta_{v})& \textrm{if } \xi  =g
			      \end{array}\right.
			    \end{equation} 
			    and 
			    \begin{equation}
			      (\hat y_{g,q}^{u})_\xi = \left\{\begin{array}{cl}
				y_\xi & \textrm{if } \xi \not=g,q,\\[0.2cm]
				(y_g  - \frac{1}{N}\delta_{u})& \textrm{if } \xi  =g, \\[0.2cm]
				(y_q  + \frac{1}{N}\delta_{u})& \textrm{if } \xi  =q.
			      \end{array}\right.
			    \end{equation}
			    If we set  $M^{N,g}_t := |\{i \in I_N:\ \xi_i(t) = g\}|$, where $(\xi_i)_{i \in \mathbb N}$ is a sequence of independent random walks on $G$ such that $\xi_1(0)$ is uniformly distributed on $G$,  then it is not hard to see that $((M^{N,g}_t /N)_{g \in G})_{t \ge 0} \Rightarrow I$ weak in $\mathcal M_1( D([0,\infty),[0,1]^G))$, where $I \equiv 1$. This follows since the one (and hence the finite) dimensional distributions converge almost surely by the strong law of large numbers together with a generator calculation in the sense of Lemma 4.5.1 and Remark 4.5.2 in \cite{EK86}. 
			    We note that this also gives convergence in probability with respect to the Skorohod metric $\dSK$, i.e. for all $\eps> 0$: 
			    \begin{equation}
			      \lim_{N \rightarrow \infty }P\left(\dSK\left(\left(\left(\frac{M^{N,g}_t }{N}\right)_{g \in G}\right)_{t \ge 0},I\right) \ge \eps\right) = 0. 
			    \end{equation}
			    If we now set $\hat \MVX^{N}_t(A \times B) := \sum_{g \in B}\hat X^{N,g}_t(A)$ and $\MVX^{N}_t(A \times B) := \sum_{g \in B} X^{N,g}_t(A)$, then
			    \begin{equation}
			      \begin{split}
				\dPr(\hat \MVX^{N}_t, \MVX^{N}_t) &\le \sum_{g \in G} \dPr(\hat X^{N,g}_t, X^{N,g}_t)  \\
				  &\le \sum_{g \in G} \sum_{x \in [0,1]} |\hat X^{N,g}_t(\{x\}) - X^{N,g}_t(\{x\})| \le \sum_{g \in G} \left|1-\frac{M^{N,g}_t }{N}\right| 
			      \end{split}
			    \end{equation}
			    and hence for all $\eps> 0$:
			    \begin{equation}
			      \begin{split}
				P(\dSK(\hat \MVX^{N}, \MVX^{N}) \ge \eps) &\le \sum_{g \in G} P\left( \int_0^\infty \sup_{0 \le t \le u} \left|1-\frac{M^{N,g}_{t} }{N}\right| du\ge \frac{\e}{|G|} \right) \\
				  &= \sum_{g \in G}P\left(\dSK\left(\left(\frac{M^{N,g}_t }{N}\right)_{t \ge 0},I\right) \ge \frac{\e}{|G|} \right) \\
				  &\stackrel{N \rightarrow \infty}{\longrightarrow }0.
			      \end{split}
			    \end{equation}
			  \end{remark}

			{\bf The spatial Kingman coalescent:} Next we introduce the spatial Kingman coalescent and show how one can construct the genealogy for a fixed time $T > 0$ in a Moran model. We construct the coalescent as in \cite{LS} - see also \cite{Ber}. \par 
			  
			  Let $N \in \mathbb N$ and write $\mathcal P([N])$ for the set of partitions of $[N]:=\{1,\ldots,N\}$, 
			  where we assume that the partition elements are ordered by their least elements, i.e. given $\pi = \{\pi_1,\ldots,\pi_m\} \in \mathcal P([N])$ we assume  $\min(\pi_i) \le \min(\pi_j)$, 
			  whenever $i \le j$. Analogue $\mathcal P(\mathbb N)$ denotes the set of partitions of $\mathbb N$, where we assume the same order as before. For $\pi \in \mathcal P([N])$
			  or $\pi \in \mathcal P(\mathbb N)$ we denote by $|\pi|$ the number of blocks (or elements) of $\pi$. \par 
			  Let $G= \{g_1,\ldots,g_m\}$ be the geographical space, defined as in section \ref{sec.TVIMM} (see \eqref{eq.geoSpace}). 
			  Then the spatial Kingman coalescent takes values in 
			  \begin{equation}\label{eq.King.StateSpace}
				\begin{split}
			  \mathcal P(\mathbb N)^G:= \{\{(\pi_i,&\xi_i): i = 1,\ldots,|\pi|\}:\  \\
				&\pi_i \in  \pi,\ \xi_i \in G,\ i = 1,\ldots,|\pi|, \ \pi \in \mathcal P(\mathbb N)\}. 
				\end{split}
			  \end{equation}
			  Analogue, we define $\mathcal P([N])^G$ as above but with $\mathcal P(\mathbb N)$ replaced by $\mathcal P([N])$. For $\pi = \{(\pi_1,\xi_1),(\pi_2,\xi_2),\ldots\} \in \mathcal P([N])^G$ or in $\mathcal P(\mathbb N)^G$
			  we define $\pi\big|_n \in \mathcal P([n])$ for $n \le N$ as the element induced by $(\pi_1\cap[n],\xi_1),(\pi_2\cap[n],\xi_2),\ldots$. 
			  Moreover we define the map $\pi:\mathcal P(\mathbb N) \to \mathbb N, \pi(i) = \min(\pi_i)$  (analogue $\mathbb N$ replaced by $ [N]$). 
			  We equip $\mathcal P(\mathbb N)^G$ with the following distance: 
			  
			  \begin{equation}
			    d^\Pi(\pi,\pi'):= \sup_{m \in \mathbb N} 2^{-m} 1_{\{\pi|_m \neq \pi'|_m\}}, \qquad \pi,\pi' \in \mathcal P(\mathbb N)^G,
			  \end{equation}
			  and $\mathcal P([N])^G$ with 
			  \begin{equation}
			    d^\Pi_N(\pi,\pi'):= \sup_{m  \le N} 2^{-m} 1_{\{\pi|_m \neq \pi'|_m\}}, \qquad \pi,\pi' \in \mathcal P([N])^G.
			  \end{equation}
			  We are now able to define the spatial Kingman coalescent.

			  \begin{proposition}\label{p.prop.King}
			    For each $\pi \in \mathcal P(\mathbb N)^G$ there is a cadlag Feller and strong Markov process $\King^G$ on $\mathcal P(\mathbb N)^G$ called the
			    spatial Kingman coalescent such that $\King^G(0) = \pi$ and 
			    \begin{itemize}
			      \item[(i)] if we write $\King^G = \{(\king^G_i,\xi_i):\ i = 1,\ldots,n\}$ for some $n \in \mathbb N$, then
				    two blocks $\king^G_i$ and $\king^G_j$ with the same label (i.e. $\xi_i = \xi_j$) coalesce according to a (non spatial) Kingman coalescent with rate $\gamma$, i.e. at rate $\gamma$ we pick independently two  elements $\King^G_i(t)$ and $\King^G_j(t)$ (at some time $t$) with the same label and have the transition $\King^G(t) \rightarrow \King^{G,i,j}(t)$, where (assume $i < j$)
				    \begin{equation}
				    \King^{G,i,j}_k(t) = \left\{ \begin{array}{ll}
				    (\king_k^G(t),\xi_k(t)),&\quad k < i,\\
				    (\king_i^G(t) \cup \king_j^G(t),\xi_i(t)),&\quad k = i,\\
				    (\king_{k}^G(t),\xi_{k}(t)),&\quad i < k < j, \\ 
				    (\king_{k+1}^G(t),\xi_{k+1}(t)),&\quad j \le k \le n-1. 
				    \end{array}\right.
				    \end{equation}
			      \item[(ii)] independently, each block with label $g_i \in G$ changes its label to $g_j \in G$ at rate $a(g_i,g_j)$. 
			    \end{itemize}
			    This process also satisfies 
			    \begin{itemize}
			      \item[(iii)] $(\King^G(t)\big|_n)_{t \ge 0}$ is a spatial Kingman coalescent started from $\King^G(0)\big|_n$, 
			    \end{itemize}
			    and its law is characterized by (iii) and the initial distribution $\pi$. 
			  \end{proposition}

			  \begin{proof}
			    This is Theorem 1 and the first Remark below this Theorem in \cite{LS}. 
			  \end{proof}

			  In the following we will always write $\king^G$ for the partition process of the spatial Kingman-coalescent $\King^G$. \\ 
			  
			  Now we give the connection to the tree-valued interacting Moran model. In order to do this recall the notations from Section \ref{sec.TVIMM}. \\

			  Let $T > 0$ be fixed. We set $A_i(h):=A_{T-h}(i,T)$ and $\hat \xi_i(h) := \xi_i((T-h)-)$, $0 \le h \le T$, $i \in  I_N$. 
			  Then $(A_i(h),\hat \xi_i(h))_{ i \in I_N}$ is a processes with values in 
			  $(I_N\times G)^{|I_N|}$ that starts in
			  $A_0(i) = i,\hat \xi_i(0) = \xi_i(T)$, $i \in I_N$ and has the following dynamic: Whenever $\pois^{i,j}(\{T-h\}) = 1$ for some $i \neq j,\ i,j \in I_N$ and 
			  $\hat \xi_i(h) = \hat \xi_j(h) $, then 
			  \begin{equation}\label{eq.anc.trans}
			    A_{k}(h) = i, \qquad \forall k \in \{l \in I_N: A_{l}(h-) = j\}.
			  \end{equation}
			  and $\hat \xi_1,\hat \xi_2,\ldots$ are independent random walks with transition kernel $\bar a(g,q) := a(q,g)$. If we define
			  \begin{equation}
			    \king_i^G(h) = \{j\in I_N:\ A_{j}(h) = A_{i}(h)\},
			  \end{equation}
			  then it is straight forward to see that $\King^G(h) := \{(\king_1^G(h),\hat\xi_1(h)),\ldots\}$ is a spatial Kingman coalescent up to time $T$ that starts in 
			  $\{(\{1\},\xi_1(T)),\ldots,(\{|I_N|\},\xi_{|I_N|}(T))\}$ (compare also the construction in \cite{Ber}, section 2.1).
			  
			  \begin{remark}\label{r.marginals.King.finite}
			    Even though the above process is only defined up to time $T$, we will in the following always assume that there is a spatial Kingman coalescent $\hat \King$ 
			    (on an extension of the probability space) defined for all times $t \in \mathbb R_+$ such that $(\hat \King_{t})_{0 \le t \le T} = (\King^G_t)_{0 \le t \le T}$.
			  \end{remark}

			  Since $\King^G$ depends on $N$ and $T$ we will write in the following $\King^{G,N,T}$ in order to indicate this dependence or $\King$ when it is clear from the context what $N,T$ and $G$ are. 
			  As a direct consequence of the construction, we get the following: 
			  
			  \begin{lemma}\label{lem.rep.Moran}
			    Let $h \ge 0$. We say $i \sim_h j$ iff there is a $ (B,g) \in \King^{G,N,T}(h)$ such that $i,j \in B$. If we define
			    \begin{equation}
			      \hat r^\king_T(i,j) := \inf\{0 \le h \le T:\ i \sim_h j\}\wedge T,\qquad i,j \in I_N,
			    \end{equation}	 
			    then (see \eqref{eq.MM.measure} for the definition of $\mu^N_T$)
			    \begin{equation}
			      \U^N_T = [I_N,\hat r^\kappa_T,\mu^N_T]. 
			    \end{equation}
			  \end{lemma}
			  
			\begin{figure}[ht]
			  \centering
			  \begin{tikzpicture}[scale=0.6]
			    \draw (-1,5)-- (-1,-3);
			    \draw (1,5)-- (1,-3);
			    \draw (3,5)-- (3,-3);
			    \draw (5,5)-- (5,-3);
			    \draw [->] (3,2) -- (5,2);
			    \draw [->] (1,-2) -- (3,-2);
			    \draw [->] (-2,-4) -- (-2,6);
			    \draw [->] (1,0) -- (-1,0);
			    \draw (-1,-4) node {\textbf{1}};
			    \draw (1,-4) node {\textbf{2}};
			    \draw (3,-4) node {\textbf{3}};
			    \draw (5,-4) node {\textbf{4}};
			    \draw (-2.5,2.2) node {$\mathbf{t_3}$};
			    \draw (-2.5,0.2) node {$\mathbf{t_2}$};
			    \draw (-2.5,-1.8) node {$\mathbf{t_1}$};
			    \draw (-2.5,3.2) node {$\mathbf{T}$};
			    \draw (-3,6) node {$\text{time}$};
			    \draw (12,3)-- (12,0);
			    \draw (12,0)-- (14,0);
			    \draw (14,0)-- (14,3);
			    \draw (13,0)-- (13,-2);
			    \draw (13,-2)-- (17,-2);
			    \draw (15,-2)-- (15,-3);
			    \draw (17,-2)-- (17,2);
			    \draw (17,2)-- (18,2);
			    \draw (17,2)-- (16,2);
			    \draw (16,3)-- (16,2);
			    \draw (18,3)-- (18,2);
			    \draw [->] (20,3) -- (20,-3);
			    \draw (12,3.7) node {$\mathbf{1}$};
			    \draw (14,3.7) node {$\mathbf{2}$};
			    \draw (16,3.7) node {$\mathbf{3}$};
			    \draw (18,3.7) node {$\mathbf{4}$};
			    \draw [dash pattern=on 5pt off 5pt] (9.5,3.)-- (18.,3.);
			    \draw [dash pattern=on 5pt off 5pt] (9.5,2.)-- (17.,2.);
			    \draw [dash pattern=on 5pt off 5pt] (9.5,0.)-- (13.,0.);
			    \draw [dash pattern=on 5pt off 5pt] (9.5,-2.)-- (15.,-2.);
			    \draw (9.5,3.5) node {\tiny $\{\{1\},\{2\},\{3\},\{4\}\}$};
			    \draw (9.5,2.5) node {\tiny $\{\{1\},\{2\},\{3,4\}\}$};
			    \draw (9.5,0.5) node {\tiny $\{\{1,2\},\{3,4\}\}$};
			    \draw (9.5,-1.5) node{\tiny $\{\{1,2,3,4\}\}$};
			    \draw [color=black] (-2,2)-- ++(-4pt,-4pt) -- ++(7pt,7pt) ++(-7pt,0) -- ++(7pt,-7pt);
			    \draw [color=black] (-2,0)-- ++(-4pt,-4pt) -- ++(7pt,7pt) ++(-7pt,0) -- ++(7pt,-7pt);
			    \draw [color=black] (-2,-2)-- ++(-4pt,-4pt) -- ++(7pt,7pt) ++(-7pt,0) -- ++(7pt,-7pt);
			    \draw [color=black] (-2,3)-- ++(-4pt,-4pt) -- ++(7pt,7pt) ++(-7pt,0) -- ++(7pt,-7pt);
			    \draw [color=black] (20,2)-- ++(-4pt,-4pt) -- ++(7pt,7pt) ++(-7pt,0) -- ++(7pt,-7pt);
			    \draw [color=black] (20,0)-- ++(-4pt,-4pt) -- ++(7pt,7pt) ++(-7pt,0) -- ++(7pt,-7pt);
			    \draw [color=black] (20,-2)-- ++(-4pt,-4pt) -- ++(7pt,7pt) ++(-7pt,0) -- ++(7pt,-7pt);
			  \end{tikzpicture}
			  \caption{\label{f.graph.constr.King}\footnotesize We are in the situation of Figure \ref{f.graph.constr}. 
			  On the left side we see the graphical construction for the TVIMM and on the right side a realization of the corresponding Kingman coalescent.}
			\end{figure}
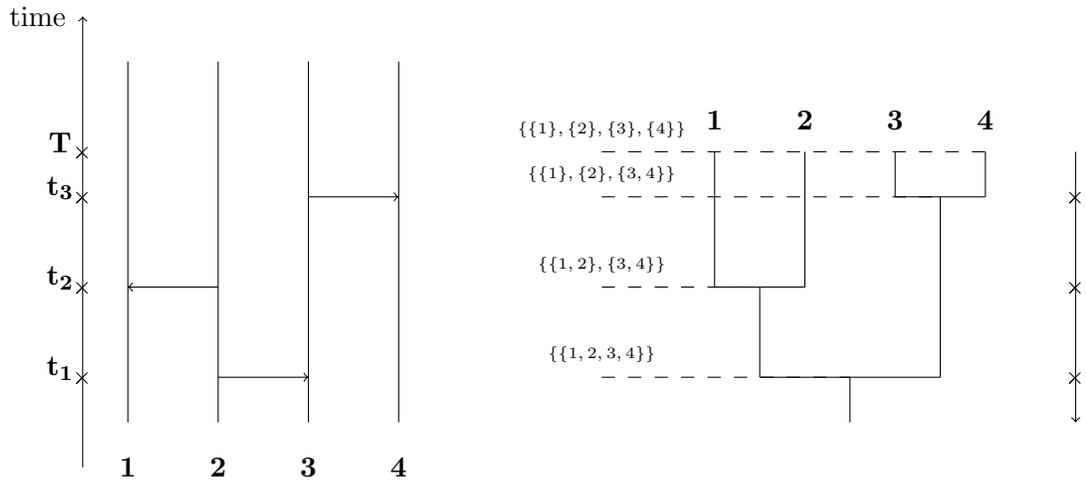
			
			Another consequence of the above construction, together with the fact that the Poisson processes used in the graphical construction have independent increments, is the following Lemma: 
			
			\begin{lemma}\label{lem.indep.increments.king}
			  The processes $(\King^{G,N,T}(h))_{0 \le h \le T}$ and $(\King^{G,N,T+T'}(h))_{0 \le h \le T'}$ are independent conditioned on $\xi_T$. 	 
			\end{lemma}	
			
			We close this section by the following assumption, which is, due to the consistency of the Kingman-coalescent, no loss of generality: 
			
			\begin{assumption}\label{a.coupl.King}
			  The Kingman coalescents are coupled for different $N \in \mathbb N$ in such a way that if we denote by $\King^{G,T}$ the spatial Kingman coalescent with $\King^{G,T}(0) = \{(\{1\},\xi_1(T)),(\{2\},\xi_2(T)),\ldots\}$, 
			  where $((\xi_i(t))_{t \ge 0})_{i \in \mathbb N}$ are independent random walks on $G$ (defined analogue as in Section \ref{sec.TVIMM}), 
			  then $\King^{G,T}\big|_{I_N} = \King^{G,N,T}$ (compare Proposition \ref{p.prop.King}) for all $N \in \mathbb N$. 
			\end{assumption}
	
		\subsubsection{A measure representation for the tree-valued interacting Moran models - Proof of Theorem \ref{thm.TVIMMM.EG}}\label{sec.MVR.TVIMM}
		  
		  Let $M \subset I_N$, $0 \le t \le T$ and define 
		  \begin{equation}
		  D_{t,T}(M):= \{i \in I_N:\ A_t(i,T) \in M\}.
		  \end{equation}
		  We call $D_{t,T}(M)$ the set of descendants at time $T$ of the individuals in $M$ that lived at time $t$. 
		  We abbreviate $D_{t,T}(x):=D_{t,T}(\{x\})$ and write $\DES_{t,T}:= \{i \in I_N: D_{t,T}(i) \neq \emptyset\} $ for the set of ancestors at time $T-t$ (measured backward).  
		  We note that $ D_{t,T}(i) = B^{r_{T}}_{T-t}(j) = \{i \in I_N:\ r_T(i,j) \le T-t\}$ is a closed ball of radius $T-t$ for all $i \in I_N$ and $j \in  D_{t,T}(i)$. It follows that 
		  \begin{equation}\label{eq.decomp}
		  I_N = \biguplus_{i \in \DES_{t,T}}  D_{t,T}(i)
		  \end{equation}
		  is the disjoint union of closed  balls (with respect to $r_{T}$) with radius $T-t$. \\ 
		  
		  Now let $T \ge 0$ and $(V^{T}_i)_{i \in \mathbb N}$ be independent (also independent of the random mechanisms in Section \ref{sec.TVIMM}) uniformly $[0,1]$-distributed random variables, where we assume w.l.o.g. that the underlying probability space is large enough. \\

		  We define for $0\le t,T$, $g \in G$
		  \begin{equation}
		    \begin{split}
		      X_t^{N,T,g}&:= \sum_{i \in \DES_{T,T+t}} \mu^N_{T+t}(D_{T,T+t}(i) \times \{g\}) \delta_{V^T_i} \\
		      &= \sum_{i \in I_N} \mu^N_{T+t}(D_{T,T+t}(i) \times \{g\})  \delta_{V^T_i}. 
		    \end{split}
		  \end{equation}
		  
		  Now  $((X_t^{N,T,g})_{g \in G})_{t \ge 0}$ is a system of interacting measure-valued Moran models that starts in 
		  \begin{equation}
		    X^{N,T,g}_0  = \frac{1}{|\{i:\ \xi_i(T) = g\}| \wedge 1}\sum_{i \in I_N}\delta_{V_i^T} 1(\xi_i(T) = g).
		  \end{equation}
		  To see this define $\hat \eta_i(t):= V_{A_T(t+T,i)}^T$ and $\hat \xi_i(t):=\xi_i(t+T)$, $i \in I_N$. Then $t \mapsto (\hat \eta_i(t),\hat \xi_i(t))_{i \in I_N}$ has the following dynamic. 
		  If $\pois^{i,j}(\{t+T\}) = 1$ and $\hat \xi_i(t) =\hat  \xi_j(t)$, then $A_T(j,T+t) = A_{T}(i,T+t)$ and hence 
		  \begin{equation}\label{eq.MM.particle.MVR}
		    (\hat \eta_k(t),\hat \xi_k(t)) \rightarrow \left\{\begin{array}{ll}
		      (\hat \eta_i(t),\hat \xi_i(t)),&\quad \text{if } k = j,\\
		      (\hat \eta_k(t),\hat \xi_k(t)),&\quad \text{if } k\neq j. 
		    \end{array}\right.
		  \end{equation}
		  Moreover, at rate $1 $ each random walk $\hat \xi_i(t)$ moves (independently of the others) from $\hat \xi_i(t)$ to $g \in G$ at rate $a(\hat \xi_i(t),g)$, i.e. we have the transition 
		  \begin{equation}
		    (\hat \eta_k(t),\hat \xi_k(t)) \rightarrow \left\{\begin{array}{ll}
		      (\hat \eta_i(t),g),&\quad \text{if } k = i, \text{ at rate } a(\hat \xi_i(t),g),\\
		      (\hat \eta_k(t),\hat \xi_k(t)),&\quad \text{if } k\neq i. 
		    \end{array}\right.
		  \end{equation}
		  It follows that 
		  \begin{equation}\label{eq.MVX.MM.1}
		    \begin{split}
		      X^{N,T,g}_t&= \frac{1}{|\{i \in I_N:\ \hat \xi_i(t) = g\}|\wedge 1} \cdot \sum_{i \in I_N} \delta_{\hat \eta_i(t)} 1(\hat \xi_i(t) = g) \\
			&=\sum_{i \in I_N} \sum_{j \in D_{T,T+t}(i)} \mu^N_{T+t}(\{j\} \times \{g\}) \delta_{\hat \eta_j(t)} \\
			&=\sum_{i \in I_N} \sum_{j \in D_{T,T+t}(i)} \mu^N_{T+t}(\{j\} \times \{g\}) \delta_{V^T_i} \\
			&=\sum_{i \in I_N} \mu^N_{T+t}(D_{T,T+t}(i) \times \{g\}) \delta_{V_{i}^T}
		    \end{split}
		  \end{equation}
		  satisfies the definition of a system of interacting Moran models. 
		  
		  \begin{remark}
		    If we look at Figure \ref{f.graph.constr.measure} (see also Figure \ref{f.graph.constr}) and define $x_i:= V_i^T$ then the above observation is clear: The mass of descendents at some time $t> T$ of the individual $x_i$ can be interpreted as the mass of the "type" $x_i$ at time $t$.  
		  \end{remark}

		  \begin{proof}(Theorem \ref{thm.TVIMMM.EG})
		   Recall the definition of $\Cut^h$ (see Definition \ref{def_Phi}), then, as a direct consequence of the above construction,  $\{(\mathcal X^{N,T}_h)_{h \ge 0}:\ T\ge 0\}$, with 
		   $\mathcal X^{N,T}_h(A \times B) = \sum_{g \in B} X^{N,T,g}_h(A)$, is a strong measure representation and hence, the TVIMM is an evolving genealogy (see Proposition \ref{prop.MR.EG}). 
		  \end{proof}
			
			We close this section with the following remark: 
			
			\begin{remark}\label{rem.ind.atoms}
				Recall that the atoms of $X^{N,T,g}$ are given by $V_{A_T(t+T,i)}^T$, where $A_T(t+T,i)$ is the ancestor of $(t+T,i)$. 
				Since the set of ancestor are sampled without replacement from the whole population, we get that, for $N \to \infty$, the 
				atoms of $X^{N,T,g}$  and $X^{N,T',g}$ for $T \neq T'$ become independent. 
			\end{remark}

		\subsubsection{Connection to the Kingman-coalescent}\label{sec.Connection.Kingman}
		  
		  Let $T \ge 0$ be fixed, $\King^{G,N,T}$ be the spatial Kingman coalescent given in Section \ref{sec.MM.King} and $V^T:= (V_i^T)_{i \in \mathbb N}$ be the sequence of independent uniformly $[0,1]$-distributed random variables, given as in the previous section. 			
		  Let $\#_h:= \#_h^T :=\#_h^{N,T}= |\king^{G,N,T}(h)|$ be the number of blocks of $\King^{G,N,T}_h$ at time $h \in [0,T]$ and $M^{N,T,g}:= \{i \in I_N:\ \xi_i(T)=g\}$ be the individuals located at $g$. Moreover, we denote by $\tau_k:= \inf \{h \ge 0: \#_h = k\}$, $k = \#_T,\ldots,|I_N|$ the coalescing times of the partition process $\king^{G,N,T}$ and note that $\tau_{|I_N|} = 0$ almost surely. \\
	
		  Let  $T,h > 0$. We define  $\bar b^{T}_h:I_N \rightarrow I_N$, which maps all elements $j$ of a partition element $\king^{G,N,T}_i(h)$ to a single element $A_{h}(j,T) $, where $A_h(j,T)$ is the ancestor of $(j,T)$ at time $h$ (measured backwards). Note that by construction of the Kingman coalescent, one has $j_1,j_2 \in \king^{G,N,T}_i(h)$ for some $i$ if and only if
		  $A_h(j_1,T) = A_h(j_2,T)$. \par 
		  Now, by the construction  via the ancestors (compare also the section before), we observe the following
		  \begin{itemize}
		  \item[(i)] The maps $\bar b^{T_1}$ and $\bar b^{T_2}$ for $T_1,T_2 >0$ satisfy the following consistency condition: 
		  For $h'>0$ and $T  > h' > 0$ one has $\bar b^{T+h}_{h+h'} = \bar b^{T}_{h'}\circ \bar b^{T+h}_h$.
		  \item[(ii)] The Moran model from the previous section satisfies 
		    \begin{equation}
		      X_h^{N,T,g} = \frac{1}{|M^{N,T+h,g}|\wedge 1}\sum_{i \in M^{N,T+h,g}} \delta_{V_{\bar b^{T+h}_{h}(i)}^T},
		    \end{equation}
		    for $g \in G$, $T,h > 0$. 
		  \end{itemize}
		
		  Now, as a consequence of Remark \ref{rem.ind.atoms}, i.e. the independence of the atoms of $X^{N,T,g}$ and $X^{N,T',g}$ for $T\neq T'$ as $N \to \infty$, 
		  and by  definition of the Kingman coalescent, we can rewrite (ii) as (note that this formula is only asymptotically true, which is ok for our purpose): 
		  For $1 \le l\le k\le K$, $g \in G$ and $T \ge 0$
		  {
		  \begin{equation}
		      \begin{split}
			&\mathcal X^{N,T+(K-k+1)\cdot h,g}_{l\cdot h}\\
			  &\phantom{xx} = \sum_{i_1 = 1}^{\#_h^{N,T + (K-k+2)\cdot h}} \sum_{i_2 \in \king^{G,N,T+(K-k+2)\cdot h}_{i_1}(h)} \cdots  \sum_{i_l \in \king^{G,N,T+(K-k+l)\cdot h}_{i_{l-1}}(h)}\\
			  & \phantom{xx}\phantom{xxxxxxxx}
			  \frac{|M^{N,T+(K-k+l+1)\cdot h,g} \cap \king^{G,N,T+(K-k+l+1)\cdot h}_{i_l}(h)|}{|M^{N,T+(K-k+l+1)\cdot h,g}|} \delta_{V_{i_1}^{T+(K-k+1)\cdot h}},
		      \end{split}
		  \end{equation}
					}
		  where $\king^{G,N,T}_i(h) := \emptyset$, whenever $i > \#^{N,T}_h$. \\

		  We will use this observations to construct a limit object. In order to do that recall that if $B \subset \mathbb N$ satisfies: 
		  \begin{equation}
		  \lim_{N \rightarrow \infty} \frac{|B \cap [N]|}{ N} =:\freq{B}
		  \end{equation}
		  exists, then we call $\freq{B}$ the {\it asymptotic frequency} of $B$. We start with the following important observation: 
		  
		  \begin{lemma}\label{l.equal.freq}
		  Let $\pi = (\pi(t))_{t \ge 0}$ be the non spatial Kingman coalescent that starts in $\{\{1\},\{2\},\ldots\}$ and let $\tau^\pi_k:=\inf \{t \ge 0:\ |\pi(t)| = k\}$. 
		  Moreover, let $\King^{G,T}$ be the spatial Kingman-coalescent and $\tau_k:= \inf\{t \ge 0: |\king^{G,T}| = k\}$. 
		  Then there is a coupling of $\pi$ and $\King^{G,T}$ such that  for all $k \in \mathbb N$ and $i = 1,\ldots,k$
		  \begin{equation}
		  \freq{\pi_i(\tau_k^\pi)} = \freq{\king^{G,T}_i(\tau_k)}\qquad \text{almost surely.}
		  \end{equation}  
		  As a consequence, 
		  $(\freq{\king^{G,T}_i(\tau_k)})_{i =1,\ldots,k}$ has the law of 
		  a variable that is uniformly distributed on the simplex 
		  \begin{equation}
		  \Delta_{k-1}:=\left\{(x_1,\ldots,x_k) \in [0,1]^k:\ \sum_{i= 1}^k x_k = 1\right\},
		  \end{equation}
		  i.e. it has the density $(k-1)! 1_{\Delta_{k-1}}$ with respect to the Lebesgue measure (see Definition 2.4 in \cite{B}).  
		  \end{lemma}
		  
		  \begin{proof}
		  Fix a $N \in \mathbb N$ and let $\pi^N:=\pi|_{I_N}$ be the non spatial coalescing that starts in $\pi^{N}(0) = \{\{1\},\ldots,\{|I_N|\}\}$. Since the locations $(\xi_i)_{i = 1,\ldots,|I_N|}$ of the partition elements of $\King^{G,N,T}$ evolve as i.i.d. random walks on $G$, with 
		  $\xi_1(0)$ being uniformly distributed on $G$, at each coalescing event we uniformly choose two partition elements over all existing partition elements and then merge them to one big block (see also the proof of Proposition 14 in \cite{LS}).  But the same holds for the non spatial process $\pi^N$. Hence we can couple both processes such that (we use the same Assumption \ref{a.coupl.King} for the processes $\pi^N$)
		  \begin{equation}\label{eq.coupl.neutr.spatial}
		  \pi^N(\tau_k^\pi) = \king^{G,N,T}(\tau_k),\qquad \forall k = 1,\ldots,|I_N|. 
		  \end{equation}
		  The result is now a consequence of Kingman's paintbox construction; see section 4.1.3 in \cite{B} (see in particular Corollary 4.1).
		  \end{proof}
		  As a consequence we get the following: 
		  
		  \begin{lemma}\label{lem.conv.fdd.MVH}
		  Assume we are in the situation of Lemma \ref{l.equal.freq}. Then for all $k \in \mathbb N$, $i = 1,\ldots,k$ and $g \in G$
		  \begin{equation}\label{eq.con.King.AS.1}
		  \frac{1}{|M^{N,T,g}|\vee 1}  |M^{N,T,g} \cap \king^{G,N,T}_i(h)| \rightarrow  a_i^{g,T}(h) =: a_i^g, 
		  \end{equation} 
		  almost surely. The family $\{(a_i^g)_{i=1,\ldots\#^{T}_{h}}:\ g \in G\}$ is exchangeable and $ (a_i)_{i = 1,\ldots,\#^{T}_{h}}$ with $a_i:= \frac{1}{|G|}\sum_{g \in G} a_i^g$ is uniformly distributed over the simplex $\Delta_{\#^{T}_{h}-1}$.
		  \end{lemma}

		  \begin{proof}
		    We define $M^{T,g}:= \{i \in \mathbb N:\ \xi_i(T) = g\}$, then $\{M^{T,g} \cap \king^{G,T}_i(\tau_k): i \in \mathbb N, g \in G\}$ is an exchangeable random partition of $\mathbb N$. Therefore, by Theorem 2.1 in \cite{B} it posses asymptotic frequency:
		    \begin{equation}
		      \frac{1}{|I_N|} |M^{N,T,g} \cap \king^{G,N,T}_i(\tau_k)| = \frac{1}{|I_N|} |M^{T,g} \cap \king^{G,T}_i(\tau_k) \cap I_N| \rightarrow \hat a_i^g
		    \end{equation}
		    almost surely. Since 
		    \begin{equation}
		      \frac{|M^{N,T,g}|}{|I_N|} = \frac{|M^{T,g} \cap I_N|}{|I_N|} \rightarrow \frac{1}{|G|}
		    \end{equation}
		    almost surely, by the strong law of large numbers for all $g \in G$, we finally get
		    \begin{equation}
		      \begin{split}
			\lim_{N \rightarrow \infty} \frac{1}{|M^{N,T,g}|\vee 1} |M^{N,T,g} \cap \king^{G,N,T}_i(\tau_k)| = |G| \hat a_i^g =: a_i^g.
		      \end{split} 
		    \end{equation}
		    The last part is a consequence of Lemma \ref{l.equal.freq}.
		  \end{proof} 
		  
		  As a consequence we get
		  
		  \begin{proposition}\label{prop.fdd.MVR}
		   Let $\{\mathcal X^{N,T}:\ T\ge 0\}$ be the measure representation of the tree-valued interacting Moran models. Then, for all $K \in \mathbb N$ and $T,h > 0$, 
		   \begin{equation}
		    \begin{split}
			&\left(\left(X^{N,T+(K-k+1)\cdot h,g}_{l\cdot h}\right)_{g \in G}\right)_{1 \le l \le k \le K}\\
			&\phantom{xx}\stackrel{N\rightarrow \infty}{\Longrightarrow} \left(\left(X^{T+(K-k+1)\cdot h,g}_{l\cdot h}\right)_{g \in G}\right)_{1 \le l \le k \le K},
		    \end{split}
		  \end{equation}
		  where for $1 \le l \le k \le K$ and $g \in G$, 
		  
		  \begin{equation}
		    \begin{split}
		      &\mathcal X^{T+(K-k+1)\cdot h,g}_{l\cdot h}\\
			&\phantom{xx} = \sum_{i_1 = 1}^{\#_h^{T + (K-k+2)\cdot h}} \sum_{i_2 \in \king^{G,T+(K-k+2)\cdot h}_{i_1}(h)} \cdots  \sum_{i_l \in \king^{G,T+(K-k+l)\cdot h}_{i_{l-1}}(h)}a_{i_l}^g \delta_{V_{i_1}^{T+(K-k+1)\cdot h}},
		    \end{split}
		\end{equation}   
		  with $a_i^g = a_i^{g,T+(K-k+l+1)h}(h)$ is given in Lemma \ref{lem.conv.fdd.MVH}.
		  \end{proposition}

		  \begin{proof}
		  This is a summary of the above observations, where 
			convergence in the weak atomic topology follows by Lemma 2.5 in \cite{EKatomic} together with our Assumption \ref{a.coupl.King}.
		  \end{proof}

	 \subsection{Proofs for Section \ref{sec.TVIFV}}
		
	    We first introduce the measure-valued Fleming-Viot process. We then prove Theorem \ref{t.convMM} and Theorem \ref{thm.FSS} by using our results on evolving genealogies. 

		    \subsubsection{Measure-valued interacting Fleming-Viot processes} \label{sec.measureFV}
		    We recall the definition of a system of interacting Fleming-Viot processes on the geographical space $G$. This process can be characterized as the solution of a well-posed martingale problem. We will only sketch the 
		    properties here, for more details see \cite{DGV} (compare also \cite{D} for a general introduction to measure-valued processes and \cite{EK86}, section 4, for the definition and properties of solutions of martingale problems).  \par
		    We denote by $\mathcal D \subset C((\mathcal M_1[0,1])^G)$ the set of {\it polynomials}, i.e. the span of functions $F$ of the following form: 
		    \begin{equation} 
		    F(x) = \prod_{i=1}^m \left(\int_{[0,1]} f_i ~d x_{\xi_i}\right), \qquad f_i \in C_b([0,1]) \ i = 1,\ldots,m,  
		    \end{equation} 
		    where $m \in \mathbb N$ and $(\xi_i)_{i = 1,\ldots,m}$ is a finite collection of elements in $G$. We note that $\mathcal D$ is an algebra that separates points and hence is dense in $C((\mathcal M_1[0,1])^G)$ by the Theorem of Stone-Weierstrass, when $\mathcal M_1[0,1]$ is equipped with the weak topology. Moreover, we can define partial derivations for elements $F \in \mathcal A$: 
		    \begin{equation}
		    \frac{\partial F(x)}{\partial x_\xi} (u) := \lim\limits_{\e\searrow 0} \frac{1}{\e} \big(F(x^{\e}(\xi,u)) -F(x)\big),
		    \end{equation} 
		    with
		    \[\big(x^{\e}(\xi,u)\big)_{\xi'}:= \left \{
		    \begin{array}{cl}
		    x_{\xi}\quad& \mbox{falls }\ \xi' \not= \xi, \\
		    x_{\xi} + \eps \delta_u \quad & \mbox{falls }\ \xi' = \xi.
		    \end{array} \right.\]
		    The second derivations are given by
		    \begin{equation}
		    \frac{\partial^2 F(x)}{\partial x_\xi x_\xi} (u,v) = \frac{\partial }{\partial x_\xi}\left(\frac{\partial F(x)}{\partial x_\xi}(u)\right)(v).
		    \end{equation}

		    Let 
		    \begin{equation}
		    Q_{x_\xi} (du,dv) = x_\xi (du) \delta_u (dv) - x_{\xi}(du) x_\xi(dv)
		    \end{equation} 
		    and define the linear operator $\Omega: \mathcal D\to C((\mathcal M_1[0,1])^G)$ by 
		    
		    \begin{equation}\label{eq.Generator.FV}
		      \begin{split} 
			\Omega F(x) = \sum_{\xi,\xi' \in G} a(\xi,\xi') &\int_{[0,1]} \frac{\partial F(x)}{\partial x_\xi} (u) \big(x_{\xi'}(du)-x_{\xi}(du)\big) \\
			  &+ \frac{\gamma}{2}\sum_{\xi \in G} \int_{[0,1]}\int_{[0,1]} \frac{\partial^2 F(x)}{\partial x_\xi \partial x_\xi} (u,v)  Q_{x_\xi}(du,dv),
		      \end{split}
		    \end{equation}
		    where $\gamma > 0$ is the resampling rate. 
		    
		    \begin{proposition}\label{p.properties.FV} (Characterization and properties of the measure-valued Fleming-Viot process)
		      The system of interacting measure-valued Fleming-Viot processes $((X_t^g)_{g \in G})_{t \ge 0}$  with $(X_0^g)_{g \in G} = (\mu_g)_{g \in G}=:\bar \mu \in (\mathcal M_1[0,1])^G$ can be characterized  by the well-posed martingale problem $(\Omega,\mathcal D,\bar \mu)$. 
		      It has the following properties: 
		      \begin{itemize}
		      \item[i)] $((X_t^g)_{g \in G})_{t \ge 0}$ has continuous paths.
		      \item[ii)] $(X_t^g)_{g \in G}$ is purely atomic for all $t > 0$ almost surely. Note that this property holds for all (even continuous) initial measures $\bar \mu$. 
		      \item[iii)] $(X_t^g)_{g \in G}$ is Markov and Feller.
		      \end{itemize}
		    \end{proposition}
		    
		    \begin{proof}
		      The above properties can be found in Theorem 0.0 in \cite{DGV}.
		    \end{proof}
		    
		    \begin{remark}
		       $(X_t^g)_{g \in G}$ has only finitely many atoms for all $t > 0$ almost surely. This follows immediately 
		      by the fact that the number of atoms is given by the number of blocks of a spatial Kingman coalescent (see Section \ref{sec.Connection.Kingman}), 
		      and that the spatial Kingman coalescent on a finite geographical space comes down from infinity (see \cite{LS}). 
		    \end{remark}
		    
		    \begin{lemma}\label{l.conv.MVR}
		      Recall that the measure representation satisfies $\MVX^{N,T}_t(\cdot \times \{g\})   =  X^{N,T,g}_t(\cdot)$ for $g \in G$, 
		      where $(X^{N,T,g})_{g \in G}$ is a system of interacting Moran models for all $T \ge 0$ (see section \ref{sec.MM.King}).
		      Let $(X^{T,g})_{g \in G}$ be the system of interacting Fleming-Viot processes with 
		      $X^{T,g}_0 = \lambda$ (the Lebesgue-measure on $[0,1]$) for all $g \in G$ and set  $\MVX^T_t(\cdot \times \{g\}) := X^{T,g}_t$, $t,T \ge 0$.
		      Then 
		      \begin{equation}
			  \MVX^{N,T} \Rightarrow \MVX^T
		      \end{equation}
		      weakly as cadlag processes, where $\mathcal M_f([0,1]\times G)$ is equipped with the topology given in 
			Definition \ref{def.weak.atomic}. 
		    \end{lemma}
		    
		    \begin{proof}
		    Since 
		    \begin{equation}
		    \MVX^{N,T}_0(\cdot \times \{g\}) = \frac{1}{|M^{N,T,g}| \wedge 1} \sum_{i \in M^{N,T,g}} \delta_{V_i^T},
		    \end{equation}
		    where $M^{N,T,g}$ is the index set of random walks that are located at $g$ at time $T$,  
		    and $|M^{N,T,g}|/N \rightarrow 1$ almost surely for all $T \ge 0$, by the strong law of large numbers, 
		    we can apply the Glivenko-Cantelli theorem (see \cite{parthasarathy}) and Lemma 2.1 in \cite{EKatomic} ($\lambda$ has no atoms) to get 
		    \begin{equation}
		    \MVX^{N,T}_0(\cdot \times \{g\}) \Rightarrow \lambda 
		    \end{equation}
		    in the weak atomic topology almost surely. \par
		    We know that the initial distributions converge, the space $\mathcal M_1([0,1])^G$ is compact, equipped with the weak topology, 
		    and the martingale problem for the measure-valued interacting Fleming-Viot processes is well-posed (see Proposition \ref{p.properties.FV}). 
		    Now the convergence $(X^{N,T,g})_{g \in G} \Rightarrow (X^{T,g})_{g \in G}$ follows by the convergence of the generators given in 
		    Remark \ref{r.Generator.MM} and \eqref{eq.Generator.FV} together with Lemma 4.5.1 and Remark 4.5.2 in \cite{EK86}, 
		    where the convergence of the generators is  a straight forward calculation and hence we skip it. 
		    As a consequence we get that 
		    
		    \begin{equation}
		      \MVX^{N,T} \Rightarrow \MVX^T
		    \end{equation}
		    weakly as cadlag processes, where $\mathcal M_f([0,1]\times G)$ is equipped with the weak topology. \par
		    It remains to prove the result, when $\mathcal M_f([0,1]\times G)$ is equipped with the (marked) weak atomic topology. 
		    According to Theorem 2.12 and Remark 2.13 in \cite{EKatomic} we need to prove that for all $S > 0$, $\delta > 0$ there is an $\eps> 0$ such that 
		    \begin{equation}
		      \liminf_{n \rightarrow \infty} P\left(\sup_{0 \le t \le S} \left( \int \int f_\eps(x,y) \MVX^{N,T}_t(dx \times G) \MVX^{N,T}_t(dy \times G)\right) \le \delta \right) \ge 1-\delta,
		    \end{equation}
		    where 
		    \begin{equation}\label{eq.def.feps}
		    f_\eps(x,y) = \Psi\left(\frac{|x-y|}{\e}\right) - 1(x = y), 
		    \end{equation}
		    and $\Psi:[0,\infty) \to [0,1]$ is an arbitrary non increasing continuous function with $\Psi(0) = 1$, $\Psi(1) = 0$. Recall Remark \ref{r.Generator.MM} (also for the definition of $\hat{\cdot}$) and observe that, 
		    by the martingale problem characterization, 
		    \begin{equation}\label{eq.Martingale.Property}
		      \begin{split}
			\int \int f_\eps(x,y) \ \hat \MVX^{N,T}_t(dx  &\times G) \hat \MVX^{N,T}_t(dy  \times G) \\
			  & =  \sum_{g,q \in G}\int \int f_\eps(x,y) \ \hat X^{N,T,g}_t(dx)\hat X^{N,T,q}_t(dy)
		      \end{split}
		    \end{equation} 
		    is a non negative martingale with cadlag paths. We can now apply Doob's martingale inequality 
		    (see for example  Proposition 2.2.16 in \cite{EK86}) and get that   
		    \begin{equation}
		      \begin{split}
			P\Big(\sup_{t \le S} &\Big( \int \int f_\eps(x,y) \hat \MVX^{N,T}_t(dx \times G) \hat \MVX^{N,T}_t(dy \times G)\Big) \ge \delta \Big) \\
			  &\le \sum_{g,q \in G} \frac{1}{\delta} E\left[ \int \int f_\eps(x,y) \ \hat X^{N,T,g}_S (dx)\hat X^{N,T,q}_S(dy) \right] \\
			  &= \sum_{g,q \in G} \frac{1}{\delta} E\left[ \int \int f_\eps(x,y) \ \hat X^{N,T,g}_0 (dx)\hat X^{N,T,q}_0(dy) \right].
		      \end{split}
		    \end{equation}
		    As we have seen above $X^{N,T,g}_0 \rightarrow \lambda$ almost surely and hence    
		    \begin{equation}
		      \begin{split}
			E&\left[ \int \int f_\eps(x,y) \ \hat X^{N,T,g}_0 (dx)\hat X^{N,T,q}_0(dy) \right] \\
			&\rightarrow \int \int \Psi\left(\frac{|x-y|}{\e}\right) \lambda(dx) \lambda(dy) 
			  \le \int_{[0,1]} \int_{y-\e}^{y+\e} \ \lambda(dx) \lambda(dy) = 2\eps. 
		      \end{split}
		    \end{equation}
		    If we choose $\eps= \frac{\delta}{2|G|^2}$ and combine this with the convergence result in Remark \ref{r.Generator.MM}, 
		    where we note that convergence in the total variation norm implies convergence in the weak atomic topology (this follows by Lemma 2.5 in \cite{EKatomic}),
		    then the above observations give the result in the weak atomic topology.    
		    \end{proof}

	 \subsubsection{The large population limit - proof of Theorem \ref{t.convMM} }

		    In order to prove Theorem \ref{t.convMM} we will use our Theorem \ref{thm.convergence.II}
		    i.e. we need to verify tightness of the tree-valued interacting Moran models and we need to check conditions (i)-(v).  \\

		    \noindent {\bf Tightness:} In order to prove tightness we verify condition (i)-(iii) of Theorem \ref{thm.tightness}. 
			Note that (i) is satisfied by Definition (see Remark \ref{r.init.IMM}) and 
			(ii) is Lemma \ref{l.conv.MVR} combined with Proposition \ref{p.properties.FV} and Remark \ref{rem.cond.ii}. Hence,
			it remains to verify (iii), i.e. for all $\eps> 0$ and $T \ge 0$ there is a $H \ge  0$ such that 
			\begin{equation}
			  \limsup_{N \rightarrow \infty} P\left(|G|^2- \sum_{x \in [0,1]}\MVX_{H}^{N,T}(\{x\}\times G)^2 \ge \varepsilon\right) \le \varepsilon.
			\end{equation}
			In order to see that observe that the result follows, if there is an almost surely finite time $C< \infty$ such that
			 $\MVX_{C}^{N,T}(\cdot \times G)$ has only one atom for all $N \in \mathbb N$. 
			 By construction of $\MVX^{N,T}$, the number of atoms of $\MVX^{N,T}_h$ is given by the number of blocks of the spatial Kingman coalescent $\King^{G,N,T+h}(h)$. 
			 Since there is an almost surely finite time $\tau$ such that $|\king^{G,T}(\tau)| = 1$ and $|\king^{G,N,T}(t)| \le |\king^{G,T}(t)| $ for all $t \ge 0$ 
			 (recall Assumption \ref{a.coupl.King}), the result follows. \\

			 \noindent {\bf Convergence:} Here we check condition (i)-(v) of Theorem \ref{thm.convergence.II}. 
			 Note that (i) is satisfied by Definition (see Remark \ref{r.init.IMM}), (ii) is Lemma \ref{l.conv.MVR},  
			 (iv) is Proposition \ref{prop.fdd.MVR} and (v) follows by Proposition \ref{prop.fdd.MVR} and Lemma \ref{lem.conv.fdd.MVH}. 
			 Hence, it remains to verify (iii). \par 
			 Similar to the proof of Lemma \ref{l.conv.MVR},  we may assume w.l.o.g.  that 
			 $t \mapsto  \hat \MVX^{N,T}_t(\{x\} \times G))$ is a positive martingale and hence we get by Doob's martingale inequality for all $\delta,\tau > 0$, $l \in \mathbb N$ and $C > 0$,
			 where we abbreviate $\mu_t^{N,l}:=\tilde{\mathcal X}^{N,l\cdot \delta}_t$ and denote by $\mathcal A(\mu)$ the set of atoms of a measure $\mu$:

			  \begin{equation}
			    \begin{split}
			      P\Big( \exists x &\in \mathcal A(\mu^{N,l}_\delta):\ \sup_{h \in [\delta,2\delta)} \mu^{N,l}_h(\{x\}) \ge C \mu^{N,l}_\delta(\{x\})| \mu^{N,l}_\delta\Big) \\
				&\le \sum_{x \in \mathcal A(\mu^{N,l}_\delta)} \frac{1}{C \cdot \mu^{N,l}_\delta(\{x\}) } E_{\mu^{N,l}_\delta}[\mu^N_{2\delta}(\{x\}) ] \\
				&= \sum_{x \in \mathcal A(\mu^{N,l}_\delta)} \frac{1}{C \cdot \mu^{N,l}_\delta(\{x\}) } \mu^{N,l}_\delta(\{x\}) \\
				&= \frac{\left|\mathcal A(\tilde{\mathcal X}^{N,l\cdot \delta}_\delta(\cdot \times G))\right|}{C}=: \frac{K^{N,l}(\delta)}{C}. 
			    \end{split}
			  \end{equation}

			  Observe that  $K^{N,l}(\delta)$ is the number of blocks of the spatial Kingman-coalescent 
			  $\King^{G,N,(l+1)\delta}$ at time $\delta$ and that this number satisfies $ K^{N,l}(\delta) \le K^{l}(\delta)  < \infty$ almost surely, where $K^l(\delta)$ denotes the number of blocks in 
			  $\King^{G,(l+1)\delta}$ at time $\delta$ (see Section 3 in \cite{LS} - also recall Assumption \ref{a.coupl.King}).  Hence, if we choose for $\eps> 0$ a $C> 0$ such that  
			  \begin{equation}
			    P(\tilde C > C) \le \frac{\e}{2\lceil \tau/\delta \rceil},\quad \text{where } \tilde C = \frac{2|K^l(\delta)| \lceil \tau/\delta \rceil}{\e}, 
			  \end{equation}
			  then we get  
			  \begin{equation}
			    \begin{split}  
			      \sup_{N \in \mathbb N} P\Big( \exists x &\in \mathcal A(\mu^{N,l}_\delta):\sup_{h \in [\delta,2\delta)} \mu^{N,l}_h(\{x\}) \ge C \mu^{N,l}_\delta(\{x\})\Big) \\  
				&= \sup_{N \in \mathbb N}  P\Big( \exists x \in \mathcal A(\mu^{N,l}_\delta)\sup_{h \in [\delta,2\delta)} \mu^{N,l}_h(\{x\}) \ge C \mu^{N,l}_\delta(\{x\}), (\tilde C \le C \vee \tilde C > C) \Big)  \\ 
				&\le \sup_{N \in \mathbb N}  P\Big( \exists x \in \mathcal A(\mu^{N,l}_\delta)\sup_{h \in [\delta,2\delta)} \mu^{N,l}_h(\{x\}) \ge \tilde C \mu^{N,l}_\delta(\{x\}) \Big)  + P(\tilde C > C) \\ 
				&\le \frac{\e}{2\lceil \tau/\delta \rceil} + \frac{\e}{2\lceil \tau/\delta \rceil} = \frac{\eps}{\lceil \tau/\delta \rceil}. 
			    \end{split}
			  \end{equation}
			  
			  It follows that

			  \begin{equation}
			    \begin{split}  
			      \sup_{N \in \mathbb N} P\Big(\exists l =1,\ldots,\lceil \tau/\delta \rceil, x &\in \mathcal A(\mu^{N,l}_\delta): \sup_{h \in [\delta,2\delta)} \mu^{N,l}_h(\{x\}) \ge C \mu^{N,l}_\delta(\{x\})\Big)\\
				&\le  \lceil \tau/\delta \rceil \cdot \frac{\eps}{\lceil \tau/\delta \rceil}= \eps. 
			    \end{split}
			  \end{equation}

			 \begin{remark}\label{rem.mvfv.MVR}
			  As we have seen in the above proof, the system of interacting measure-valued Fleming-Viot processes that starts in  the Lebesgue measure $\lambda$ on $[0,1]$ at each site, is a weak measure representation for the 
			  tree-valued interacting Fleming-Viot process.
			 \end{remark}

		\subsubsection{Proof of Theorem \ref{thm.FSS}; the finite system scheme}
				  
			  Let $\{\MVX^{N,T}: \ T \ge 0\}$ be a measure representation of the tree-valued interacting Fleming-Viot processes on $G_N$. 
			  Then, by definition, $\{\hat \MVX^{N,T}:\ T \ge 0\}$ given by 
			  \begin{equation}
			  \hat \MVX^{N,T}_h = \theta_N\left(\MVX^{N,T|G_N|}_{h|G_N|}\right),\qquad T,h \ge 0,
			  \end{equation}
			  defines a measure representation of $\hat \U^N:=(h_N(\U_{T |G_N|}^N))_{T \ge 0}$. 
			  We need to check two things, tightness of $\hat \U^N$ and convergence of the finite dimensional distributions (see Remark \ref{rem.conv.fdd}).\\

			  \noindent {\bf Tightness of $\hat \U^N$:} We start by proving tightness of $\hat \U^N:=(h_N(\U_{t|G_N|}^N))_{t \ge 0} $,
			  where we check condition (i)-(iii) of Theorem \ref{thm.tightness}.
			  Note that  (i) is satisfied  by definition and observe that  we can write 
			  \begin{equation}
			    \hat \MVX^{N,T}_h = \frac{1}{|G_N|}\sum_{g \in G_N} X^{N,T|G_N|,g}_{h|G_N|}, 
			  \end{equation}
			  where $(X^{N,T|G_N|,g})_{g \in G_N}$ is a system of measure-valued interacting Fle\-ming-Viot processes 
			  (see Remark \ref{rem.mvfv.MVR}) for all $T \ge 0$. In particular, it is identically distributed for different $T$, i.e. $(X^{N,T|G_N|,g})_{g \in G_N} \stackrel{d}{=} (X^{N,g})_{g \in G_N}$. \par
			  A classical result (see \cite{FSS} and \cite{DGV}) says that $(\frac{1}{|G_N|}\sum_{g \in G_N} X^{N,g}_{h|G_N|})$ converges for $N \rightarrow \infty$, when $\mathcal M_1([0,1])$ is 
			  equipped with the weak topology, i.e.  
			  \begin{equation}
			  \hat \MVX^{N,T} \Rightarrow \hat \MVX^{T}, 
			  \end{equation}
			  where $\hat \MVX^T$ is a (non spatial) Fleming-Viot process with diffusion rate $D$. \par
			  Next we prove, that this result stays valid when $\mathcal M_1([0,1])$ is equipped with the weak atomic topology. 
			  In order to do that we note that $\sum_{g \in G_N} X^{N,g}_t(f)$ is a martingale for all bounded continuous functions $f$ 
			  (this is a direct consequence of the martingale-problem characterization of this process). 
			  Moreover, the class of functions $f$, such that $\sum_{g \in G_N} X^{N,g}_t(f)$ is a martingale, is closed under bounded pointwise convergence. Hence the process 
			  \begin{equation}
			  \sum_{g \in G_N}\sum_{g \in G_N} \int \int f_\eps(x,y)   X^{N,g}_t(dx)  X^{N,g}_t(dy)
			  \end{equation}
			  where $f_\eps \ge 0$ is given in \eqref{eq.def.feps}, is a continuous non-negative martingale. 
			  When we now apply Doob's martingale inequality, then (ii) follows analogue to the proof of Lemma \ref{l.conv.MVR}. \par 
			  In order to prove (iii), we need to check if for $\eps> 0$ and $T \ge 0$ there is a $H \ge  0$ such that  
			  \begin{equation}
				\begin{split}
				  \limsup_{N \rightarrow \infty} &P\left(\Big|1 - \sum_{x \in [0,1]}{\hat \MVX}_{H}^{N,T}(\{x\})^2\Big|> \varepsilon\right)
				  \le \varepsilon.
				 \end{split}
			  \end{equation}
						  
			  Note that the number of atoms of ${\hat \MVX}^{N,T}_h$ is given by the number of blocks of the spatial Kingman coalescent 
			  $\King^{G_N,(T+h)|G_N|}(h|G_N|) \stackrel{d}{=} \King^{G_N,T}(h|G_N|) $ (in the sense of Remark \ref{r.marginals.King.finite}). Since the block counting process of the spatial Kingman coalescent converges in this time scale to a non-spatial Kingman coalescent with coalescing rate $D$ (see \eqref{eq.block.conv} below) we get the result analogue to (iii) in the proof of Theorem \ref{t.convMM} (also compare the coupling argument for the f.d.d. convergence). \\

			  \noindent {\bf F.d.d. convergence of $\hat \U^N$:}  Let $T \ge 0$, $\hat \MS_h^N$ be the non increasing reordering of the block frequencies $(\freq{\king^{G_N,T|G_N|}_i(|G_N|h)})_{i = 1,\ldots,|\king^{G_N,T|G_N|}(|G_N|h)|}$ and observe that the Kingman coalescents 
			  $\King^{G_N,T}$ are identically distributed for different $T$ (in the sense of Remark \ref{r.marginals.King.finite}), i.e. the law of $\hat \MS^N$ does not depend on $T$. \par 
			  By Theorem 19 in \cite{LS}, for all $\delta > 0$: 
			  \begin{equation}\label{eq.block.conv}
			     (\hat \#^N_h)_{h \ge \delta} \Rightarrow (\hat \#_h)_{h \ge \delta}
			  \end{equation}
			  where $\hat \#^N_h := \hat \#^{N,T}_h := |\king^{G_N,T}(|G_N|h)|$ and $\hat \#_h := \hat \#^{T}_h := |\bar \king^{T}(h)|$ is the block counting process and $\bar \king^{T}(h)$
			  is the non spatial Kingman coalescent with coalescing rate $D$. By the convergence in the Skorohod topology (compare also Lemma 3.8.1 in \cite{EK86}), 
			  \begin{equation}
			   (\hat \tau_k^N)_{k = 1,\ldots,\hat \#^N_\delta} \Rightarrow  (\hat \tau_k)_{k = 1,\ldots,\hat \#_\delta},
			  \end{equation}
			  where 
			  \begin{align}
			   \hat \tau^N_k &:= \inf \{h \ge 0: \ \hat \#^N_h = k\},\\
			   \hat \tau_k &:= \inf \{h \ge 0: \ \hat \#_h = k\},
			  \end{align}
			  are  the jump times of the block counting processes. \par 
			  Note that these jump times are exactly the jump times of the corresponding Kingman coalescents and hence, when we apply Lemma \ref{l.equal.freq}, this gives
			  \begin{equation}
			   (\hat \MS_h^N)_{h \ge \delta} \Rightarrow (\hat \MS_h)_{h \ge \delta},
			  \end{equation}
			  where $\hat \MS_h$ is the non increasing reordering of $(\freq{\bar \king^{T}_i(h)})_{i = 1,\ldots,\hat \#_h}$. \par 
			  			  
			  Now, observe that the convergence of the one (or finite) dimensional distributions can be proven in terms of the family size decomposition, when the limit function $f$  is the 
			  family size decomposition of some $\UM$-valued variable that satisfies condition (v) of Theorem \ref{thm.convergence.II}, i.e. in our situation, 
			  for all $\delta > 0$ and given $f(\delta)$ has $K$ non zero entries, $(f(\delta)_i)_{i = 1,\ldots,K}$ is absolutely continuous to the Lebesgue measure $\lambda|_{[0,1]}^K$ 
			  (see Remark \ref{rem.conv.fdd} combined with Lemma \ref{l.ord.f}). By Proposition \ref{prop.fdd.MVR} (compare also Lemma \ref{lem.conv.fdd.MVH}), the reordered sizes of atoms of the measure
			  representation are given by the reordered block frequencies of the Kingman coalescent, i.e. 
			  \begin{equation}
			  \begin{split}
			   (\f(\hat \U^N_T,h|G_N|))_{h \in [\delta,T]} &=  (\hat \MS_h^N)_{h \in [\delta,T]} \\
			   &\Rightarrow (\hat \MS_h)_{h \in [\delta,T]} = (\f(\bar \U_T,h))_{h \in [\delta,T]}
			  \end{split}
			  \end{equation}
			  and convergence of the one dimensional marginals follows.\par
			  In order to complete the proof, we recall that the  Kingman coalescents $\King^{G_N,T|G_N|}(|G_N|h)$ and $\King^{G_N,(T+h)|G_N|}(|G_N|h)$ are independent given their locations 
			  (see Lemma \ref{lem.indep.increments.king}). Following the proof of Proposition 14 in \cite{LS}, the locations of $\King^{G_N,T|G_N|}(|G_N|h)$ and $\King^{G_N,(T+h)|G_N|}(|G_N|h)$
			  are asymptotically independent and therefore, $\king^{G_N,T|G_N|}(|G_N|h)$ and $\king^{G_N,(T+h)|G_N|}(|G_N|h)$ are asymptotically independent. By Proposition \ref{prop.fdd.MVR}, 
			  this gives the result.


\end{document}